\documentclass[letterpaper, 11pt,  reqno]{amsart}

\usepackage{amsmath,amssymb,amscd,amsthm,amsxtra, esint}
\usepackage{graphicx}

\usepackage[implicit=true]{hyperref}

\usepackage[left=32mm, right=32mm, 
bottom=27mm]{geometry}
\allowdisplaybreaks[2]

\usepackage{color}

\definecolor{gr}{rgb}   {0.,   0.69,   0.23 }
\definecolor{bl}{rgb}   {0.,   0.5,   1. }
\definecolor{mg}{rgb}   {0.85,  0.,    0.85}
\definecolor{yl}{rgb}   {0.8,  0.7,   0.}
\definecolor{or}{rgb}  {0.7,0.2,0.2}

\newtheorem{theorem}{Theorem} [section]

\newtheorem{lemma}[theorem]{Lemma}
\newtheorem{proposition}[theorem]{Proposition}
\newtheorem{remark}[theorem]{Remark}

\newtheorem{corollary}[theorem]{Corollary}

\newtheorem*{ack}{Acknowledgments}

\DeclareMathOperator*{\supp}{supp}

\newcommand{\1}{\hspace{0.5mm}\text{I}\hspace{0.5mm}}
\newcommand{\II}{\text{I \hspace{-2.8mm} I} }
\newcommand{\III}{\text{I \hspace{-2.9mm} I \hspace{-2.9mm} I}}

\newcommand{\noi}{\noindent}
\newcommand{\Z}{\mathbb{Z}}
\newcommand{\R}{\mathbb{R}}
\newcommand{\C}{\mathbb{C}}
\newcommand{\T}{\mathbb{T}}

\let\Re=\undefined\DeclareMathOperator*{\Re}{Re}
\let\Im=\undefined\DeclareMathOperator*{\Im}{Im}

\let\P= \undefined
\newcommand{\P}{\mathbf{P}}

\newcommand{\E}{\mathbf{E}}
\newcommand{\D}{\mathbb{D}}
\DeclareMathOperator{\GNS}{GNS}

\newcommand{\M}{\mathcal{M}}

\newcommand{\N}{\mathcal{N}}
\newcommand{\NB}{\mathbb{N}}

\newcommand{\F}{\mathcal{F}}

\def\norm#1{\|#1\|}

\newcommand{\al}{\alpha}
\newcommand{\be}{\beta}
\newcommand{\ga}{\gamma}
\newcommand{\dl}{\delta}

\newcommand{\nb}{\nabla}

\newcommand{\Dl}{\Delta}
\newcommand{\eps}{\varepsilon}
\newcommand{\kk}{\kappa}
\newcommand{\g}{\gamma}

\newcommand{\ld}{\lambda}

\newcommand{\s}{\sigma}

\newcommand{\ft}{\widehat}

\newcommand{\wt}{\widetilde}
\newcommand{\cj}{\overline}
\newcommand{\dx}{\partial_x}
\newcommand{\dt}{\partial_t}
\newcommand{\dd}{\partial}

\newcommand{\ta}{\theta}

\renewcommand{\l}{\ell}
\renewcommand{\o}{\omega}

\newcommand{\les}{\lesssim}
\newcommand{\ges}{\gtrsim}

\newcommand{\jb}[1]
{\langle #1 \rangle}

\newcommand{\ind}{\mathbf{1}}

\newcommand{\ds }{\mathrm{d}}

\DeclareMathOperator{\Id}{id}

\numberwithin{equation}{section}
\numberwithin{theorem}{section}

\newcommand{\too}{\longrightarrow}

\renewcommand{\d}{\ds }

\begin{document}
\baselineskip = 14pt

\title[Optimal mass threshold for the  focusing  Gibbs measure]
{Optimal integrability threshold for Gibbs measures associated with 
focusing NLS on the torus}

\author[T.~Oh, P.~Sosoe,  and L.~Tolomeo]
{Tadahiro Oh, 
Philippe Sosoe,
and Leonardo Tolomeo}

\address{
Tadahiro Oh\\
 School of Mathematics\\
The University of Edinburgh\\
and The Maxwell Institute for the Mathematical Sciences\\
James Clerk Maxwell Building\\
The King's Buildings\\
Peter Guthrie Tait Road\\
Edinburgh\\ 
EH9 3FD\\
 United Kingdom}

\email{hiro.oh@ed.ac.uk}

\address{Philippe Sosoe\\
Department of Mathematics\\
Cornell University\\ 
310 Malott Hall\\ 
Cornell University\\
 Ithaca\\ New York 14853\\ USA}

\email{psosoe@math.cornell.edu}

\address{
Leonardo Tolomeo\\ 
School of Mathematics\\
The University of Edinburgh\\
and The Maxwell Institute for the Mathematical Sciences\\
James Clerk Maxwell Building\\
The King's Buildings\\
Peter Guthrie Tait Road\\
Edinburgh\\ 
EH9 3FD\\
 United Kingdom
 and 
Mathematical Institute, Hausdorff Center for Mathematics, Universit\"at Bonn, Bonn, Germany}

\email{tolomeo@math.uni-bonn.de}



\begin{abstract}

We study an optimal mass threshold for normalizability 
of  the Gibbs measures associated with the focusing mass-critical nonlinear Schr\"odinger equation
on the one-dimensional torus. 
In   an influential paper, Lebowitz, Rose, and Speer (1988)
proposed a critical mass threshold given by the mass of the ground state
on the real line.
We provide a proof for the optimality of this critical mass threshold.
The proof also applies to the two-dimensional radial problem posed on the unit disc. 
In this case, we answer a question posed by  Bourgain and Bulut (2014)
on the optimal mass threshold.

Furthermore, 
in the one-dimensional case, 
we  show that the Gibbs measure is indeed normalizable at 
the optimal mass threshold, thus answering an open question 
posed by Lebowitz, Rose, and Speer (1988).
This normalizability  
at the optimal mass threshold is rather striking in view
of the minimal mass blowup solution
for the focusing quintic  nonlinear Schr\"odinger equation
on the one-dimensional torus. 

\end{abstract}



%
\maketitle

\vspace{-3mm}

\tableofcontents

%
%
%
%

\section{Introduction}

\subsection{Focusing Gibbs measures}

In this paper, we continue the study of the focusing Gibbs measures
for the nonlinear Schr\"odinger equations (NLS), 
initiated in the seminal papers by Lebowitz, Rose, and Speer~\cite{LRS}
and Bourgain \cite{BO94}.
A focusing Gibbs measure $\rho$ 
is a probability measure on 
functions\,/\,distributions with a formal density:
\begin{align*}
\d\rho = Z^{-1} e^{- H(u)} \d u , 
\end{align*}

\noi
where $Z$ denotes the partition function and 
the Hamiltonian functional $H(u)$ is given by 
\begin{align}
H(u) = \frac 12\int |\nb u|^2 \ds x - \frac 1 p\int |u|^p \ds x.
\label{Gibbs2}
\end{align}

\noi
The NLS equation:
\begin{align}
i \dt u +  \Dl u + |u|^{p-2} u = 0, 
\label{NLS1}
\end{align}

\noi
generated by the Hamiltonian functional $H(u)$, 
 has been studied extensively as models for describing  
 various  physical phenomena ranging from Langmuir waves in plasmas to signal propagation in optical fibers \cite{SS, KMM, Agrawal}. 
 Furthermore, 
 the study of the equation~\eqref{NLS1}
 from the point of view of the (non-)equilibrium statistical mechanics 
 has received wide attention; see for example \cite{LRS, BO94, BO96, BO97,Tzv1,  Tzv2, LMW, BBulut, 
CFL, DNY2, Bring3}.
See also \cite{BOP4} for a survey on the subject, more from the dynamical point of view.
Our main goal in this paper is to study the construction of the  focusing Gibbs measures
on   
 the one-dimensional torus 
$\T  = \R/\Z$ and   the two-dimensional unit disc $\D \subset \R^2$
(under the radially symmetric assumption
with the Dirichlet boundary condition)
and 
 determine
optimal mass thresholds of their normalizability 
in the {\it critical} case.
In particular, we 
resolve an issue in the Gibbs measure construction on~$\T$ \cite[Theorem 2.2]{LRS}
and also 
answer a question posed by Bourgain and Bulut \cite[Remark 6.2]{BBulut}
on 
the optimal mass threshold for the focusing Gibbs measure on the unit disc $\D$.
Furthermore, in the case of the one-dimensional torus, 
we prove normalizability at the optimal mass threshold
in spite of the existence of the minimal mass blowup solution
to NLS at this mass, 
thus answering an open question posed by Lebowitz, Rose, and Speer~\cite{LRS}.

We first go over the case of the one-dimensional torus.
Consider the mean-zero  Brownian loop $u$ on $\T$, defined by the Fourier-Wiener series: 
\begin{equation}
\label{bloop}
u(x)=\sum_{n\neq 0}\frac{g_n(\o)}{2\pi n}e^{2\pi inx},
\end{equation}

\noi
where $\{g_n\}_{n \in \Z\setminus\{0\}}$ 
denotes  a sequence of independent standard complex-valued\footnote{Namely, 
$\Re g_n$ and  $\Im g_n$ are independent  real-valued mean-zero 
  Gaussian random variables
  with variance $\frac 12$.}
  Gaussian random variables.
Then, 
the law $\mu_0$ of  
  the mean-zero  Brownian loop $u$ in \eqref{bloop}
  has the formal density 
 given by 
\begin{align*}
\d \mu_0 = Z^{-1} e^{-\frac 12 \int_\T |  u'|^2 \d x} \d u .
\end{align*}

\noi
The main difficulty in constructing the focusing Gibbs measures
comes from the unboundedness-from-below
of the Hamiltonian $H(u)$.
This makes the problem very different from the defocusing case,
which is a well studied subject in constructive Euclidean quantum field theory.
In \cite{LRS}, Lebowitz, Rose, and Speer proposed
to consider the Gibbs measure with an $L^2$-cutoff:\footnote{Recall that the $L^2$-norm
is conserved under the NLS dynamics.}
\begin{align*}
\d\rho & = Z_{p, K}^{-1} 
\, e^{\frac{1}{p}\int_\T |u|^p\,\ds  x}
\, \ind_{\{\|u \|_{L^2(\T)} \leq K\}}\, 
\d \mu_0 
\end{align*}

\noi
and claimed the following results.
\begin{theorem}\label{THM:LRS}

Given $p > 2$ and $K > 0$, define
the partition function $Z_{p, K}$ by 
\begin{equation}\label{Z1}
Z_{p,K}=\E_{\mu_0}\Big[e^{\frac{1}{p}\int_\T |u|^p\,\ds  x}\mathbf{1}_{\{\|u\|_{L^2(\T)}\le K\}}\Big], 
\end{equation}

\noi
where $\E_{\mu_0}$ denotes an expectation with respect
to the law $\mu_0$ of  
  the mean-zero  Brownian loop 
  \eqref{bloop}.
Then, the following statements hold\textup{:}

\smallskip

\noi
\begin{itemize}
\item
[\textup{(i)}] \textup{(subcritical case)}
If $2< p<6$, then $Z_{p,K}<\infty$ for any $K>0$.

\smallskip

\noi
\item[\textup{(ii)}] \textup{(critical case)}
Let $p=6$.
Then, $Z_{6,K}<\infty$ if $K<\|Q\|_{L^2(\R)}$, and $Z_{6,K}=\infty$ if $K>\|Q\|_{L^2(\R)}$. Here, $Q$ is the \textup{(}unique\footnote{Up to the symmetries.}\textup{)} optimizer for the Gagliardo-Nirenberg-Sobolev inequality
on $\R$
such that $\|Q\|_{L^6(\R)}^6 = 3\|Q'\|_{L^2(\R)}^2$.

\end{itemize}

\end{theorem}

There remains a question
of normalizability at the optimal threshold
$K= \|Q\|_{L^2(\R)}$
in the critical case ($p = 6$).
We address this issue in Subsection \ref{SUBSEC:1.2}.

Lebowitz, Rose, and Speer proved
the non-normalizability result for $K > \|Q\|_{L^2(\R)}$
in Theorem \ref{THM:LRS}\,(ii)
by using a Cameron-Martin-type theorem
and the following 
sharp Gagliardo-Nirenberg-Sobolev (GNS) inequality on $\R^d$:
\begin{equation}\label{GNS0}
\|u\|_{L^p(\R^d)}^p \le C_{\GNS}(d,p)\|\nb u\|^{\frac{d(p-2)}{2}}_{L^2(\R^d)}
\|u\|^{2+\frac{p-2}{2}(2-d)}_{L^2(\R^d)}
\end{equation}

\noi
with $d = 1$ and $p = 6$.
See  Section \ref{SEC:GNS}
 for a further discussion on 
 the sharp GNS inequality.

The threshold value $p=6$ and the relevance of the GNS inequality  can be understood at an intuitive level by formally rewriting \eqref{Z1} as a functional integral with respect to the (periodic) Gaussian free field
(= the mean-zero Brownian loop in \eqref{bloop}):
\begin{equation}\label{Z2}
Z_{p,K} \text{ ``=" } \int_{\{\|u\|_{L^2(\T)}\le K\}} e^{-\frac{1}{2}\int_\T |u'(x)|^2\,\ds  x
+\frac{1}{p}\int_\T |u(x)|^p\,\ds  x}\d u .
\end{equation}

\noi
Applying the GNS inequality \eqref{GNS0}, this quantity is bounded by
\[\int_{\{\|u\|_{L^2(\T)}\le K\}} e^{-\frac{1}{2}\int_\T |u'(x)|^2\,\ds  x+\frac{C_{\GNS}(1, p)}{p}
K^{\frac{p+2}2}\left(\int_\T |u'(x)|^2\,\ds  x\right)^{\frac{p-2}4}}\,\d u.\]

\noi
Thus, 
 when $p<6$
 or when $p=6$ and $K$ is sufficiently small, 
  we expect the Gaussian part of the measure to dominate, and 
  hence the partition function to be finite.

Regarding the construction of the focusing Gibbs measure, 
a pleasing probabilistic proof of Theorem \ref{THM:LRS} based on this idea was given in \cite{LRS}, using the explicit joint density of the times that the Brownian path hits certain levels on a grid. 
Unfortunately, 
as pointed out by 
Carlen, Fr\"ohlich, and Lebowitz 
\cite[p.\,315]{CFL}, 
there is a gap 
in the proof of Theorem~2.2 
in~\cite{LRS}.
More precisely, 
%
%
%
 the proof in \cite{LRS} seems to apply only to the case, where the expectation in the definition of $Z_{p,K}$ is taken with respect to a standard (``free'') Brownian motion started at 0, rather than the random periodic function \eqref{bloop}. 

Subsequently, 
a more analytic proof due to Bourgain appeared in \cite{BO94}, 
establishing normalizability of the focusing Gibbs measure
(i.e.~$Z_{p, K} < \infty$)
for (i) $2 < p < 6$ and any $K > 0$
and for (ii) $p = 6$ and sufficiently small $K > 0$.
 His argument combines basic estimates for Gaussian vectors with the Sobolev embedding to identify the tail behavior of the random variable $\int_\T |u|^p\,\ds  x$, subject to the condition $\|u\|_{L^2(\T)}\le K$. 
 It also applies  to the case $p=6$, but shows only that $Z_{6,K}<\infty$ for sufficiently small $K>0$.

As the first main result in this paper, 
we obtain the optimal threshold
when $p = 6$
claimed in Theorem~\ref{THM:LRS}\,(ii)
by proving $Z_{6,K}<\infty$ for any  $K < \|Q\|_{L^2(\R)}$.
In particular, our argument  resolves
the issue in \cite{LRS} mentioned above.
Our proof is closer in spirit to Bourgain's, since it uses the series representation \eqref{bloop} 
 of the Brownian loop, as opposed to the path space approach taken in \cite{LRS}. 
 In Section \ref{SEC:Bourgain}, 
 we go over Bourgain's argument and point out that, 
 in this approach, closing the gap between small $K$ and the optimal threshold seems difficult. 
 We then present our proof of the direct implication of 
 Theorem~\ref{THM:LRS}\,(ii)
 in Subsection \ref{SUBSEC:1d}.
As in~\cite{LRS}, the idea is to make rigorous the computation suggested by \eqref{Z2} by 
a finite dimensional approximation.

\begin{remark} \label{REM:loop}\rm

Theorem \ref{THM:LRS} also applies
 when we replace
the mean-zero  Brownian loop in~\eqref{bloop}
by 
 the Ornstein-Uhlenbeck loop: 
\begin{equation}
\label{bloop3}
u(x)=\sum_{n \in \Z}\frac{g_n(\o)}{\jb{n}}e^{2\pi inx},
\end{equation}

\noi
where $\jb{n} = (1 + 4\pi^2 |n|^2)^\frac{1}{2}$
and  $\{g_n\}_{n \in \Z}$ is a sequence of independent standard complex-valued  Gaussian random variables.
See Remark \ref{REM:gauss}.
The same comment applies to Theorem~\ref{THM:OPT} below.
The law $\mu$ of  the Ornstein-Uhlenbeck loop has 
the formal density 
\begin{align}
\d \mu = Z^{-1} e^{-\frac 12 \|u\|_{H^1(\T)}^2} \d u .
\label{Gibbs3a}
\end{align}

\noi
As seen in \cite{BO94}, 
$\mu$ is a more natural base Gaussian measure to consider for
the nonlinear Schr\"odinger equations, 
due to the lack of the conservation of the spatial mean under the dynamics.

We also point out that 
Theorem \ref{THM:LRS} also holds in the real-valued setting.
The same  comment applies to Theorems \ref{THM:BB} and \ref{THM:OPT}.
For example, this  is relevant to the study of the generalized KdV equation
(gKdV)
on $\T$:
\begin{align}
\dt u + \dx^3 u + \dx (u^{p-1}) = 0.
\label{KdV}
\end{align}

\end{remark}

Our method also applies to 
the focusing Gibbs measures on 
the two-dimensional unit disc $\D \subset \R^2$, 
under the radially symmetric assumption
with the Dirichlet boundary condition.
In the subcritical case ($p < 4$), 
Tzvetkov \cite{Tzv1}
constructed the focusing Gibbs measures,
along with the associated  invariant dynamics.
His analysis was complemented 
in  \cite{BBulut} by a study of the critical  case $p=4$, 
under a small mass assumption.
See also \cite{Tzv1, Tzv2, BBulut} for results in the defocusing case.

Our approach to Theorem \ref{THM:LRS}
allows us to 
establish the optimal mass threshold in the critical case ($p = 4$), 
thus 
answering the question posed by Bourgain and Bulut in \cite[Remark 6.2]{BBulut}.
We first introduce some notations.
Let $\D=\{(x,y)\in \R^2: x^2+y^2<1\}$ be the unit disc. 
Let $J_0(r)$ be the  Bessel function of order zero, defined by
\[J_0(x)=\sum_{j=0}^\infty \frac{(-1)^j}{(j!)^2}\left(\frac{x}{2}\right)^{2j},\]
and $z_n$, $n\ge 1$, be its successive, positive zeros. 
Then, it is known \cite{Tzv1} that $\{e_n\}_{n \in \NB}$ defined by 
\[e_n(r) =\|J_0(z_n\cdot)\|^{-1}_{L^2(\D)}J_0(z_n r), \quad 0\le r\le 1,\]

\noi
forms an orthonormal basis of $L^2_\text{rad}(\D)$, 
consisting of 
the radial eigenfunctions of 
the Dirichlet self-adjoint realization of $-\Dl$ on $\D$.
Here, $L^2_\text{rad}(\D)$ denotes the subspace
of  $L^2(\D)$, 
consisting of radial functions.
Now, consider
 the random series:
\begin{equation}\label{bloop2}
v(r)=\sum_{n=1}^\infty \frac{g_n(\o)}{z_n} e_n( r), \quad r^2=x^2+y^2, 
\end{equation}

\noi
where  $\{g_n\}_{n \in \NB}$ is a sequence of independent standard complex-valued  Gaussian random variables.

\begin{theorem}\label{THM:BB}

Given $p > 2$ and $K > 0$, define
the partition function $\wt Z_{p, K}$ by 
\begin{align}
\wt Z_{p,K}= \E\Big[e^{\frac{1}{p}\int_{\D}|v|^p\,\ds  x}\mathbf{1}_{\{\|v\|_{L^2(\D)}\le K\}}\Big],
\label{Z3}
\end{align}

\noi
where $\E$ denotes an expectation with respect
to the law of the random series \eqref{bloop2}.
Then, the following statements hold\textup{:}

\smallskip

\noi
\begin{itemize}
\item[\textup{(i)}] \textup{(subcritical case)}
If $p<4$, then $\wt Z_{p,K}<\infty$ for any $K>0$.

\smallskip

\noi
\item[\textup{(ii)}] \textup{(critical case)}
Let  $p=4$.  Then,  $\wt Z_{4,K}<\infty$ if $K<\|Q\|_{L^2(\R^2)}$, 
and  $\wt Z_{4,K}= \infty$ if $K> \|Q\|_{L^2(\R^2)}$, where $Q$ is the optimizer for the 
Gagliardo-Nirenberg-Sobolev inequality \eqref{GNS0} on $\R^2$
such that $\|Q\|_{L^4(\R^2)}^4 = 2\|\nb Q\|_{L^2(\R^2)}^2$.

\end{itemize}

\end{theorem}

Part (i) of Theorem \ref{THM:BB} is due to  Tzvetkov \cite{Tzv1}.
In \cite{BBulut}, Bourgain and Bulut 
considered the critical case ($p = 4$)
and proved $\wt Z_{4,K}<\infty$ if $K \ll 1$, 
and  $\wt Z_{4,K}= \infty$ if $K\gg 1$, leaving a gap.

Theorem \ref{THM:BB}\,(ii) answers the question posed by Bourgain and Bulut in \cite{BBulut}.
See Remark~6.2 in~\cite{BBulut}.
We present the proof of Theorem \ref{THM:BB}\,(ii)
in Subsection~\ref{SUBSEC:2d}
and Section~\ref{SEC:non}.
In Subsection~\ref{SUBSEC:2d}, 
we prove 
$\wt Z_{4,K}<\infty$ for the optimal range $K<\|Q\|_{L^2(\R^2)}$
by following our argument for Theorem \ref{THM:LRS}
on $\T$.
We point out that some care is needed here due to the 
growth
of the $L^4$-norm of the eigenfunction $e_n$;
see \eqref{H6}.
In Section~\ref{SEC:non}, 
we prove $\wt Z_{4,K} = \infty$ for  $K> \|Q\|_{L^2(\R^2)}$.
Our argument of the non-normalizability 
follows closely that on $\T$ by Lebowitz, Rose, and Speer \cite{LRS}.

\subsection{Integrability at the optimal mass threshold}
\label{SUBSEC:1.2}

We now consider the normalizability issue of the focusing Gibbs
measure on $\T$ in the critical case ($p = 6$)
 at the optimal threshold
$K= \|Q\|_{L^2(\R)}$, 
at which a phase transition takes place.
Before doing so, let us first 
discuss the situation
for the associated dynamical problem, 
namely, the focusing quintic NLS, \eqref{NLS1} with $p = 6$.
On the real line, the optimizer $Q$ 
for the sharp Gagliardo-Nirenberg-Sobolev inequality \eqref{GNS0}
is the ground state for the associated elliptic problem (see \eqref{elliptic} below).
Then, by applying the pseudo-conformal transform to 
the solitary wave solution $Q(x)e^{2it}$, 
we obtain the minimal mass blowup solution 
to the focusing quintic NLS on $\R$.
Here, the minimality refers to the fact that 
 any solution to 
the focusing quintic NLS on $\R$
with 
$\|u\|_{L^2(\R)} < \|Q\|_{L^2(\R)}$ exists globally in time; see \cite{weinstein}.

In \cite{OgT}, Ogawa and Tsutsumi constructed
an analogous minimal mass blowup solution~$u_*$
with $\|u_*\|_{L^2(\T)} = \|Q\|_{L^2(\R)}$
to the focusing quintic NLS on the one-dimensional torus $\T$.
It was also shown that, as time approaches a blowup time, 
the potential energy 
$\frac 1 6\int |u_*|^6 \ds x$ tends to $\infty$.
In view of the structure of the partition function $Z_{6, K= \|Q\|_{L^2(\R)}}$
in~\eqref{Z1}, 
this divergence of the potential energy seems to create
a potential obstruction
to the construction of the focusing Gibbs
measure in the current setting.
 In \cite{BO97}, Bourgain wrote
``One remarkable point concerning the normalizability problem for Gibbs-measures of NLS in the focusing case is its close relation to blowup phenomena in the classical theory. Roughly speaking, this may be understood as follows. After normalization, the measure would be forced to live essentially on ``blowup data'', which however is incompatible with the invariance properties under the flow.''


In spite of the existence of the minimal mass blowup solution, 
we prove that the focusing critical Gibbs measure
is normalizable at the optimal threshold
  $K=\|Q\|_{L^2(\R)}$.

\begin{theorem}\label{THM:OPT}
Let   $K=\|Q\|_{L^2(\R)}$.
Then, 
the one-dimensional partition function $Z_{6,K}$ in~\eqref{Z1} is finite.
\end{theorem}

In view of the discussion  above, 
Theorem \ref{THM:OPT}
was unexpected and is rather surprising.
Theorem \ref{THM:OPT}  answers an open question posed by 
Lebowitz, Rose, and Speer in \cite{LRS}.
See Section~5 in \cite{LRS}.
Moreover, together with Theorem \ref{THM:LRS}\,(ii), 
Theorem \ref{THM:OPT} shows that 
 the partition function $Z_{6, K}$ is not analytic in the cutoff parameter $K$, 
 thus settling another question posed in \cite[Remark~5.2]{LRS}.
Compare this with 
the  subcritical case ($p< 6$), 
where 
the analyticity result
of the partition function on the parameters (including the inverse temperature, 
which we do not consider here)
was proved by 
Carlen, Fr\"ohlich, and Lebowitz~\cite{CFL} (for slightly different Gibbs measures).

The proof of Theorem \ref{THM:OPT} is presented in 
 Section \ref{SEC:opt} and 
constitutes
the major part of this paper, involving 
ideas and techniques from 
various branches
of mathematics: probability theory, functional inequalities, 
elliptic partial differential equations (PDEs), spectral analysis, etc.
We break the proof into several steps:

\smallskip
\begin{enumerate}
\item
In the first step, we use a profile decomposition
and  establish 
a  stability result for the  GNS inequality \eqref{GNS0};
see Lemma \ref{LEM:u-Q}.
When combined with our proof of Theorem~\ref{THM:LRS},
this stability result shows that if the integration is restricted to 
the complement of 
$U_\eps=\{u \in L^2(\T) : \|u-Q\|_{L^2(\T)}<  \eps\}$ for suitable $\eps \ll 1$, 
then the resulting partition function is finite. (In fact, we must exclude a neighborhood of the orbit of 
the ground state $Q$ under translations, rescalings, and rotations, but we ignore this technicality here.) 
Thus,  the question is reduced to the evaluation of the functional integral
\begin{equation}\label{Hamil}
  \int_{U_\eps} e^{-H(u)}\, \ind_{\{\|u\|_{L^2(\T)}\le K \}}\d u
  \end{equation}

\noi
in the neighborhood $U_\eps$ of (the orbit of) the ground state $Q$, 
where $H(u)$ is as in~\eqref{Gibbs2} with $p = 6$.

\smallskip

\item In the second step, 
we show that when $\eps > 0$ is sufficiently small, 
i.e.~when $U_\eps$ lies in a sufficiently small neighborhood of 
the (approximate) soliton manifold 
$\M = \big\{ 
 e^{i \ta}  Q_{\dl, x_0}^\rho: 
0 < \dl < \dl^* , \, x_0\in  \T,
\text{ and } \ta \in \R  \big\}$, 
where $Q_{\dl, x_0}^\rho = (\tau_{x_0} \rho) Q_{\dl, x_0} $ denotes the dilated and translated ground state
(see \eqref{G3-1} and \eqref{QQ1}) and $\rho$ is a suitable cutoff function 
for working on the torus  $\T\cong [-\frac 12,\frac 12)$ (see \eqref{rhodef}), 
we can endow $U_\eps$ with an orthogonal coordinate system
in terms of the (small) dilation parameter $0< \dl  < \dl^*$, 
the translation parameter $x_0 \in \T$, 
the rotation parameter $\ta \in \R/(2\pi\Z)$, 
and the component $v\in L^2(\T)$ orthogonal  to the soliton manifold $\M$.
See
Propositions \ref{PROP:FOC} and \ref{PROP:FOCN}.

\smallskip

\item 
We then introduce a change-of-variable formula
and reduce the integral \eqref{Hamil} 
to an integral in $\dl$, $x_0$, $\ta$, and $v$.
See Lemma \ref{LEM:surface}.

\smallskip

\item 
In 
Subsection \ref{SUBSEC:X1}, 
we reduce the problem to 
estimating 
a certain Gaussian integral with the integrand
given by 
\begin{align*}
\exp\left(-(1-\eta^2) \langle A w, w\rangle_{H^1(\T)}\right)
\end{align*}

\noi
for some small $\eta > 0$.
Here, 
$A = A(\dl, \eta)$ denotes the operator on $H^1(\T)$ defined 
in~\eqref{eqn: Adef}:
\begin{equation*}
Aw= P^{H^1}_{V'}(1-\dx^2)^{-1}\bigg(\dl^{-2}P^{H^1}_{V'}w-(1+5\eta)(\rho Q_\dl)^4\Big(2\Re (P^{H^1}_{V'}w)
+\frac{1}{2}P^{H^1}_{V'}w\Big)\bigg), 
\end{equation*}

\noi
where $V' \subset H^1(\T)$ is 
as in \eqref{XX8}.
See also 
\eqref{eqn: Vdelta-def}.
We point out that  the operator $A$ is closely related to 
the second variation $\dl^ 2 H$ of the Hamiltonian.
See Lemma~\ref{lem: B-lwr-bd}.
In view of the compactness of the operator
$A$, the issue is further reduced to 
 estimating the eigenvalues of $\frac12\Id+(1-\eta^2)A$.
Subsection \ref{SUBSEC:spec}
is devoted to the spectral analysis of the operator~$A$.

\end{enumerate}

\smallskip

\noi
For readers' convenience, 
we present the summary of the proof 
of Theorem \ref{THM:OPT}
in Subsection~\ref{SUBSEC:pf}.

\begin{remark}\rm
There exists an extensive literature 
on the study of soliton-type behavior to dispersive PDEs on $\R^d$; see, for example, \cite{nsbook}
and the references therein.
In particular, there are existing works on $\R^d$ which are closely related
to Steps (2) and (4) described above.
In Remarks \ref{REM:Leo1} and \ref{REM:Leo2}, 
we provide brief comparison of our work on $\T$
with those on $\R^d$, pointing out similarities and differences.

\end{remark}

We now state 
a dynamical consequence
of Theorems \ref{THM:LRS} and \ref{THM:OPT}.

\begin{corollary}\label{COR:NLS}
Let $p = 6$.
Consider the Gibbs measure $\rho$ with the formal density
\begin{align}
\d\rho & = Z_{6, K}^{-1} 
\, e^{\frac{1}{6}\int_\T |u|^6\,\ds  x}
\, \ind_{\{\|u \|_{L^2(\T)} \leq K\}}\, 
\d \mu 
\label{Gibbs4}
\end{align}

\noi
where $\mu$
is the law of  the Ornstein-Uhlenbeck loop 
in \eqref{bloop3}.
If $K \le \|Q\|_{L^2(\R)}$, 
then
the focusing quintic NLS, 
\eqref{NLS1} with $p = 6$, on $\T$
is almost surely globally well-posed
with respect to the Gibbs measure $\rho$.
Moreover, the Gibbs measure $\rho$ is invariant under the NLS dynamics.

By imposing that 
the Ornstein-Uhlenbeck loop 
in \eqref{bloop3} is real-valued \textup{(}i.e.~$g_{-n} = \cj{g_n}$, $n \in \Z$\textup{)}
or by 
replacing $\mu$ in \eqref{Gibbs4}
with the law $\mu_0$
of  the real-valued mean-zero  Brownian loop in \eqref{bloop} 
 \textup{(}with $g_{-n} = \cj{g_n}$, $n \in \Z\setminus\{0\}$\textup{)}, 
a similar result holds for the focusing quintic generalized KdV,
\eqref{KdV} with $p = 6$.

\end{corollary}

%

When $K <  \|Q\|_{L^2(\R)}$, 
Corollary \ref{COR:NLS}
follows from the deterministic local well-posedness results
for the quintic NLS \cite{BO93}
and the quintic gKdV  \cite{CK} in the spaces
containing the support of the Gibbs measure,
combined with 
Bourgain's invariant measure argument \cite{BO94}.
When $K =   \|Q\|_{L^2(\R)}$, 
the density $e^{\frac{1}{6}\int_\T |u|^6\,\ds  x}
\, \ind_{\{\|u \|_{L^2(\T)} \leq K\}}$ is only in $L^1(\mu)$
and thus  
Bourgain's invariant measure argument is not directly applicable.
In this case, however, 
the desired claim 
follows from 
the corresponding result for $K =   \|Q\|_{L^2(\R)} -\eps$, $\eps > 0$, 
and 
the dominated convergence theorem by taking $\eps \to 0$.

\begin{remark}	\rm
(i)
A  result analogous to Theorem \ref{THM:OPT} presumably holds for the two-dimensional radial Gibbs measure 
on $\D$
studied in Theorem \ref{THM:BB}.
Moreover, we expect the analysis on $\D$ to be slightly simpler
since the problem on $\D$ has fewer symmetries (in particular, no translation invariance).
In order to limit the length of this paper, however, we 
do not pursue this issue here.

\smallskip

\noi 
(ii) 
In \cite{OQ}, Quastel and the first author
constructed the focusing Gibbs measure conditioned 
 at a specified mass, provided that the mass is sufficiently small
 in the critical case ($p = 6$).
This answered another question posed in \cite{LRS}.
See also \cite{CFL, Bre}.
However, the argument in~\cite{OQ}, based on Bourgain's approach, 
is not quantitative. 
Thus,  it would be of interest to further investigate
 this problem to see if
the focusing Gibbs measure conditioned 
 at a specified mass 
  in the critical case ($p = 6$)
 can be indeed constructed
 up to the optimal mass threshold
as in Theorem \ref{THM:OPT}.
We point out that the key ingredients for the proof of Theorem \ref{THM:OPT}
(see the steps (1) - (4) right after the statement of Theorem \ref{THM:OPT})
hold true even in the case of the focusing critical Gibbs measure
(with the critical power $p = 6$) 
restricted to the critical $L^2$-norm  $\|u\|_{L^2(\T)} = \|Q\|_{L^2(\R)}^2$, 
and thus we expect that the focusing critical Gibbs measure 
at a specified mass
is normalizable even at the optimal mass threshold.
We, however, do not pursue this issue further in this paper.

\end{remark}

\begin{remark}\rm
(i) 
In view of the minimal mass blowup solution
to NLS, 
the normalizability of the focusing critical Gibbs measure 
at the optimal threshold
$K=\|Q\|_{L^2(\R)}$
in Theorem~\ref{THM:OPT}
was somehow unexpected.
As an afterthought, 
we may give some reasoning for this phenomenon, 
referring to certain properties of 
the minimal mass blowup solution
on the real line
(which are not known in the periodic setting).

The first result is the 
uniqueness\,/\,rigidity of the minimal mass blowup solution on $\R$
due to Merle \cite{Merle93}, 
which states that 
if an $H^1$-solution  $u$ with $\|u\|_{L^2(\R)} = \|Q\|_{L^2(\R)}$
to the focusing quintic NLS on $\R$ blows up 
in a finite time, it must be the minimal mass blowup solution
up to the symmetries of the equation.\footnote{In a recent preprint \cite{Dodson}, 
Dodson extended this 
uniqueness\,/\,rigidity result of the minimal mass blowup solution to the $L^2(\R)$-setting.
We, however, point out that our understanding of the corresponding problem on the torus $\T$ is rather poor in this direction, 
in particular in a low regularity setting.  For example, 
even local well-posedness  in $L^2(\T)$
of the focusing quintic NLS on $\T$ remains a challenging open problem
after Bourgain's work \cite{BO93}. See also \cite{Kishimoto}.}
See also an extension \cite{LZ}
of this result for rougher $H^s$-solution, $s > 0$,
which holds only on $\R^d$ for $d \geq 4$ under the radial assumption.
While an analogous result is not known in the periodic setting
(and in low dimensions), 
 these results may indicate non-existence of rough
 blowup solutions
at the critical  threshold
$K=\|Q\|_{L^2(\R)}$.
Theorem \ref{THM:OPT} shows that 
this non-existence claim   holds true probabilistically.

Another point is instability of the minimal mass blowup solution.
The minimal mass blowup solution 
 is intrinsically unstable because a mass subcritical perturbation leads to a globally defined solution.
See also \cite{MRS}.
Such instability may be related to the
fact that we do not see (rough perturbations of) the minimal mass blowup solution
probabilistically.

Lastly, we mention the situation for the focusing quintic gKdV on the real line.
While finite time blowup solutions to gKdV ``near'' the ground state
are known to exist \cite{Merle, MM1},
it is also known that there is {\it no}
minimal mass blowup solution to the focusing quintic gKdV on the real line 
\cite{MM2}.
Thus, from the gKdV point of view, 
the normalizability in Theorem~\ref{THM:OPT}
is perhaps naturally expected
(but we point out that analogues of 
the results in 
\cite{Merle, MM1, MM2} are not known on $\T$).

\smallskip

\noi
(ii)
In recent years, there has been  a growing interest in studying 
nonlinear dispersive equations with random initial data; see, for example, 
\cite{BO96, BT1, CO, LM, BT3, BOP1, BOP2, Poc, OP, OOP, BOP3, DNY3, Bring3}.
In particular, 
there are recent works \cite{OOP, Bring, FM} on stability of  finite time blowup solutions
under rough and random perturbations.
The so-called log-log blowup solutions to the mass-critical focusing NLS on $\R^d$ were constructed by 
Perelman \cite{Perel}
and 
Merle and Rapha\"el
\cite{MR1, MR2, MR3, MR4}.
When $d = 1, 2$, these log-log blowup solutions on $\R^d$ are known to be stable
under $H^s$-perturbations for $s > 0$; see \cite{Rap, CR}.
When $d = 2$, 
Fan and Mendelson~\cite{FM}  proved stability of the log-log blowup solutions
to the mass-critical focusing NLS
under random but structured $L^2$-perturbations.
 (Some of) these results on the log-log blowup solutions (at least the deterministic ones)
are expected to hold in the periodic setting due to the local-in-space nature
of the blowup profile.  See, for example,~\cite{PR}
for the construction of the log-log blowup solutions on  
a domain in $\R^2$.
We point out that the log-log blowup solutions mentioned above have mass strictly 
greater than (but close to) the mass of the ground state, 
which is complementary to the regime we study in this paper with regard to   the construction of the focusing critical Gibbs measure.

\end{remark}

\begin{remark}\rm

While the construction of the defocusing Gibbs measures
has been extensively studied
and well understood
due to the strong interest in constructive Euclidean quantum field theory, 
the (non-)normalizability issue of the focusing Gibbs measures, 
going back  to the work of Lebowitz, Rose, and Speer \cite{LRS}
and Brydges and Slade \cite{BS}, 
is not fully explored.
In \cite{BO97}, 
Bourgain wrote 
``It seems worthwhile to investigate this aspect 
[the (non-)normalizability issue of the focusing Gibbs measures]
more as a continuation of 
\big[\cite{LRS}\big] and \big[\cite{BS}\big].''
See   related works
\cite{Rider, BBulut, CFL, OOT, OST, OOT2, TW}
on the non-normalizability (and other issues)
for focusing Gibbs measures.

In a recent series of works 
\cite{OOT, OOT2}, the first and third authors with Okamoto
employed
the variational approach due to Barashkov and Gubinelli \cite{BG}
to  study the following two {\it critical} focusing\footnote{By ``focusing'', 
we also mean the  non-defocusing case, such as the cubic interaction 
appearing in \eqref{FG2}, 
such that the interaction potential 
(for example, $\frac{\s}3 \int_{\T^3} u^3 \d x$ in \eqref{FG2}) is unbounded from above.}  models
on the three-dimensional torus $\T^3$:

\begin{itemize}
\item[(i)]
the focusing Hartree Gibbs measure 
with a Hartree-type quartic interaction, formally written as
\begin{align}
\d\rho(u) = Z^{-1} \exp \bigg(\frac{\s}4 \int_{\T^3} (V*u^2)u^2 \d x\bigg) \d\mu_3(u), 
\label{FG1}
\end{align}

\noi
where the coupling constant $\s > 0$ corresponds to the focusing interaction and 
 $V$ is (the kernel of) 
the Bessel potential 
of order $\be >0$ given by 
\begin{align*}
 V*f = \jb{\nabla}^{-\beta} f = (1-\Dl)^{-\frac{\be}{2}} f.
\end{align*}

\noi
Hereafter, $\mu_3$ denotes the massive Gaussian free field on $\T^3$.
When $\be = 2$, the focusing Hartree model \eqref{FG1} turns out to be critical.

Recall that the  Bessel potential of order $\be$ on $\T^3$
can be written (for some $c>0$) as 
\begin{align} \label{Be2}
V(x) = c |x|^{\be-3} + K(x)
\end{align}
for $0<\be<3$ and $x \in \T^3 \setminus \{ 0 \}$,
where $K$ is a smooth function on $\T^3$.
See Lemma 2.2 in \cite{ORSW}.
Thus,  when $\be = 2$, the potential $V$ essentially corresponds to the Coulomb potential
$V(x) = |x|^{-1}$, which is of particular physical relevance.

\smallskip

\item[(ii)]
the $\Phi^3_3$-measure, formally written as
\begin{align}
\d\rho(u) = Z^{-1} \exp \bigg(\frac{\s}3 \int_{\T^3} u^3 \d x\bigg) \d\mu_3(u), 
\label{FG2}
\end{align}

\noi
where  the coupling constant $\s \in \R\setminus \{0\}$ measures the strength
of the cubic interaction.
Since $u^3$ is not sign definite,
the sign of $\s$ does not play any role
and, in particular, the problem is not  defocusing even if $\s < 0$.
We point out that the $\Phi^3_3$-model makes sense only in 
the real-valued setting.

\end{itemize}

\smallskip

In the three-dimensional setting, 
the massive Gaussian free field $\mu_3$ is supported
on $H^s(\T^3) \setminus H^{-\frac{1}{2}}(\T^3)$ for any $s < -\frac 12$.
Thus, the potentials in \eqref{FG1} and \eqref{FG2}
do not make sense as they are given, and proper renormalizations
need to be introduced.
Furthermore, due to the focusing nature of the problems, 
one needs to 
endow the measures with certain taming.
In 
\cite{OOT, OOT2}, the first and third authors with Okamoto
studied the 
 generalized grand-canonical Gibbs measure
 formulations of the focusing Hartree Gibbs measure in \eqref{FG1}
 and the $\Phi^3_3$-measure in \eqref{FG2}.
 For example, 
 the 
 generalized grand-canonical Gibbs measure
 formulations of the  $\Phi^3_3$-measure in \eqref{FG2}
 is given by 
\begin{align}
\d\rho(u) = Z^{-1}
 \exp \bigg( \frac \s3 \int_{\T^3} :\! u^3\!: \, \d x - A
\bigg|\int_{\T^3} :\! u^2 \! : \, dx\bigg|^3 - \infty \bigg) \d \mu_3(u), 
\label{FG3}
\end{align}

\noi
where $:\! u^k\!:$ denotes the standard Wick renormalization
and the term $- \infty$ denotes another (non-Wick) renormalization.
See  the work by 
Carlen, Fr\"ohlich, and Lebowitz \cite{CFL}
for a discussion of 
the generalized grand-canonical Gibbs measure in the one-dimensional setting.
See also Remark~2.1 in \cite{LRS}.

In \cite{OOT}, the first and third authors with Okamoto
established a phase transition
in the following two respects:
(i.a)  the 
 focusing Hartree Gibbs measure in \eqref{FG1}
is constructible for $\be > 2$, 
while it is not for $\be < 2$
and (i.b)~when $\be = 2$, 
the  focusing Hartree Gibbs measure
is constructible  in the weakly nonlinear regime $0 < \s \ll 1$, 
while it is not in the strongly nonlinear regime $\s \gg 1$.
This shows that the focusing Hartree Gibbs measure is critical when  $\be = 2$.

In terms of scaling, the $\Phi^3_3$-model 
corresponds to the focusing Hartree model in the critical case $\be = 2$.
Indeed, it was shown in \cite{OOT2} that the $\Phi^3_3$-model is also critical, 
exhibiting the following phase transition; 
the $\Phi^3_3$-measure is constructible in the weakly nonlinear regime 
 $0 < |\s| \ll 1$, 
whereas  it is not in the strongly nonlinear regime $|\s| \gg 1$.
While 
the 
focusing Hartree Gibbs measure in \eqref{FG1}
is absolutely continuous with respect the base massive  Gaussian free field $\mu_3$
even in the critical case ($\be = 2$), 
it turned out that 
  the $\Phi^3_3$-measure in~\eqref{FG2} 
  is singular with respect to the base massive Gaussian free field $\mu_3$. 
This singularity of the $\Phi^3_3$-measure
introduced additional difficulties
(as compared to the focusing Hartree Gibbs measures studied in \cite{OOT})
 in both the measure (non-)construction part
and the dynamical part in \cite{OOT2}.
See \cite{OOT2} for a further discussion.

In view of the aforementioned results in \cite{OOT, OOT2}, 
it is of interest to investigate (existence of) 
a threshold value  $\s_* >0$ (depending on the models)
such that the construction of the critical focusing Hartree Gibbs measure with $\be = 2$ 
(and the $\Phi^3_3$-measure, respectively)
holds
for $0 < \s < \s_*$ (and for $0 < |\s| < \s_*$, respectively), 
while 
the non-normalizability of 
the critical focusing Hartree Gibbs measure with $\be = 2$ 
(and the $\Phi^3_3$-measure, respectively)
holds
for $\s > \s_*$ (and for $ |\s| > \s_*$, respectively).
If such a threshold value $\s_*$ could be determined, 
it would also be of interest to study normalizability at the threshold
$\s = \s_*$ in the focusing Hartree case (and $|\s| = \s_*$
in the $\Phi^3_3$-case), analogous to Theorem \ref{THM:OPT}
in the one-dimensional case.
Such a problem, however, requires optimizing all the estimates
in the proofs in \cite{OOT, OOT2}
and is out of reach at this point.

Several comments are in order.
As mentioned above, the $\Phi^3_3$-model can be considered
only in the real-valued setting.
Furthermore, 
the 
 critical focusing Hartree model with $\be = 2$ 
 and the $\Phi^3_3$-model 
 are mass-subcritical (whereas the critical cases studied in Theorems~\ref{THM:LRS}\,(ii), \ref{THM:BB}\,(ii), and \ref{THM:OPT}
 are all mass-critical).
 Hence, the critical nature
 of these models do not seem to have anything to do with 
 finite-time blowup solutions (in particular
 to NLS) unlike Theorems \ref{THM:LRS}, \ref{THM:BB}, and \ref{THM:OPT}
 studied in this paper.
While we mentioned the results 
on the 
 generalized grand-canonical Gibbs measure
 formulations 
 (namely with a taming by the Wick-ordered $L^2$-norm)
 of the  focusing Hartree measures and the $\Phi^3_3$-measure, 
 analogous results hold
 even when we consider the (non-)construction of these measures
 endowed with 
a   Wick-ordered $L^2$-cutoff.
See, for example,  Remark 5.10 in \cite{OOT}.
We point out that even in this latter setting
(namely with a   Wick-ordered $L^2$-cutoff), 
what matters is the size of the coupling constant $\s$
and the size of the Wick-ordered $L^2$-cutoff does not play any role
(unlike Theorems \ref{THM:LRS}, \ref{THM:BB}, and \ref{THM:OPT}).
See also \cite{OST}.

Lastly, let us mention the dynamical aspects of these models.
For both of the models~\eqref{FG1} and \eqref{FG2}, it is possible to study 
the standard (parabolic) stochastic quantization~\cite{PW} (namely the associate
stochastic nonlinear heat equation) 
and the canonical stochastic quantization \cite{RSS}
(namely the associate stochastic damped nonlinear wave equation).
In the parabolic case, well-posedness
follows easily from 
the standard first order expansion as in \cite{McKean, BO96, DPD2};
see \cite[Appendix A]{OOT} and \cite{EJS}.
Due to a weaker smoothing property, 
the well-posedness issue in the hyperbolic setting
becomes more challenging.
By adapting the paracontrolled approach,
 originally introduced in the parabolic setting 
 \cite{GIP}, to the wave setting \cite{GKO2}, 
the first and third authors with Okamoto
\cite{OOT, OOT2} constructed global-in-time dynamics
for the stochastic damped nonlinear wave equations
associated with the focusing Hartree model (for $\be \ge 2$)
and the $\Phi^3_3$-model. 
As for the focusing Hartree model endowed with a Wick-ordered $L^2$-cutoff, 
Bourgain \cite{BO97}
studied the associated
 Hartree NLS on $\T^3$
 and constructed global-in-time dynamics when $\be > 2$.
This result was extended to the critical  
case ($\be = 2$) in \cite{OOT, DNY4}, 
where, in \cite{DNY4}, Deng, Nahmod, and Yue 
proved well-posedness\footnote{The local well-posedness argument in \cite{DNY4} applies
to the range $\be > \be_*$ for some  $\be_* < 1$ sufficiently close to~1.}
of the  associated Hartree NLS on $\T^3$
by using the random averaging operators, 
originally introduced in \cite{DNY2}.

%
%

\end{remark}

\subsection{Notations}
We write $ A \les B $ to denote an estimate of the form $ A \leq CB $
for some $C> 0$.
Similarly, we write  $ A \sim B $ to denote $ A \les B $ and $ B \les A $ and use $ A \ll B $ 
when we have $A \leq c B$ for some small $c > 0$.
We may use subscripts to denote dependence on external parameters; for example,
 $A\les_{p} B$ means $A\le C(p) B$,
 where the constant $C(p)$ depends on a parameter $p$. 

In the following, we deal with complex-valued functions viewed 
as elements in {\it real} Hilbert and Banach spaces.
In particular, with $M = \T$,  $\D$, or $\R$, the inner product on $H^s(M)$
is given by
\begin{align}
 \jb{f, g}_{H^s(M)} = \Re \int_M (1 - \Dl)^s f(x) \cj {g(x)} \, \ds x.
\label{inner1}
\end{align}

\noi
Note that with the inner product \eqref{inner1}, 
the family $\{e^{2\pi in x}\}_{n \in \Z}$ does not form an orthonormal basis of $L^2(\T)$.
Instead, we need to 
use $\{e^{2\pi in x}, i e^{2\pi in x}\}_{n \in \Z}$
as an orthonormal basis of $L^2(\T)$.
A similar comment applies to the case of the unit disc $\D$. 
We  point out that the series representations 
such as \eqref{bloop}
are not affected by whether we use the inner product~\eqref{inner1}
with the real part or that without the real part.
For example, in \eqref{bloop}, 
we have 
\[g_n e^{2\pi i n x} = (\Re g_n) e^{2\pi i n x}  + (\Im g_n)\,  i e^{2\pi i n x}.  \]

\noi
Here,  the right-hand side is more directly associated with  the 
the inner product \eqref{inner1}, while the left-hand side 
 is  associated with  the inner product without the real part.

Given $N \in \NB$, 
we denote by $\pi_N$ the Dirichlet  projection (for functions on $\T$) onto frequencies $\{|n|\leq N\}$:
\begin{align}
 \pi_N f(x) = \sum_{|n| \leq N } \ft f(n) e^{2\pi in x}
\label{EN0}
 \end{align}

\noi
and we set 
\begin{align}
E_N = \pi_N L^2(\T) = \text{span}\big\{e^{2\pi i n x}:  |n| \leq N\big\}.
\label{EN1}
\end{align}

\noi
We also define $\pi_{\ne 0}$ to be the orthogonal projection
onto the mean-zero part of a function: 
\begin{align}
 \pi_{\ne 0} u = \sum_{n \in \Z\setminus\{0\}} \ft u(n) e^{2\pi inx}
\label{EN2}
\end{align}

\noi
and set $\pi_0 = \Id - \pi_{\ne 0}$.


Given $k \in \Z_{\geq 0} := \NB \cup\{0\}$, 
let $P_k$ be  
the  Littlewood-Paley projection onto frequencies of order $2^k$
defined by 
\begin{align}
P_{k} u= \sum_{2^{k-1}<  |n|\leq  2^k} \ft{u}(n)e^{2\pi i nx}.
\label{LP1}
\end{align}

\noi
Similarly, set
\begin{align}
P_{\le k}u & = 
\sum_{j = 0}^k u_j = \sum_{|n|\le 2^k} \ft{u}(n)e^{2\pi i nx},
\label{LP2}\\
P_{\ge k}u & = 
\sum_{j = k}^\infty u_j =
\sum_{|n|> 2^{k-1}} \ft{u}(n)e^{2\pi i nx}.
\label{LP3}
\end{align}

Given measurable sets $A_1, \dots, A_k$, 
we 
 use the following notation:
 \begin{align}
 \E[f(u), A_1, \dots,  A_k]
  = \E\bigg[f(u)\prod_{j = 1}^k \ind_{A_j}\bigg], 
 \label{EXP1}
 \end{align}
 
\noi 
where $\E$ denotes an expectation with respect to 
a probability distribution for $u$ under discussion.

This paper is organized as follows.
In Section \ref{SEC:Bourgain}, we review Bourgain's argument from~\cite{BO94},
which will be used in our proof of the direct implication of Theorem \ref{THM:LRS}\,(ii) in Subsection~\ref{SUBSEC:1d}.
In Section \ref{SEC:GNS}, we  go over 
the Gagliardo-Nirenberg-Sobolev  inequality \eqref{GNS0} on $\R^d$
and discuss its variants on $\T$ and $\D$.
In Section~\ref{SEC:below},  
we then establish the direct implications
of Theorem~\ref{THM:LRS}\,(ii) on the one-dimensional torus $\T$ (Subsection \ref{SUBSEC:1d})
and Theorem~\ref{THM:BB}\,(ii) on the unit disc $\D$ (Subsection~\ref{SUBSEC:2d}).
In Section~\ref{SEC:non}, we prove the non-normalizability claim 
in Theorem~\ref{THM:BB}\,(ii).
Finally, we 
prove normalizability of the focusing critical Gibbs measure
 at the optimal mass threshold
(Theorem \ref{THM:OPT}) in Section~\ref{SEC:opt}.

\section{Review of Bourgain's argument}
\label{SEC:Bourgain}

In this section, we reproduce Bourgain's argument in \cite{BO94} for the proof  of Theorem~\ref{THM:LRS}\,(i).  Part of the argument presented below will be used in Subsection \ref{SUBSEC:1d}.
Let $2 < p \leq 6$ and $u$ denote the 
 mean-zero  Brownian loop $u$ in \eqref{bloop}.
 Rewriting \eqref{Z1} as
\begin{equation*}
Z_{p,K}\le1+\int_0^\infty \ld^{p-1}e^{\frac{1}{p}\ld^p}
\P\big(\|u\|_{L^p(\T)}>\ld, \|u\|_{L^2(\T)}\le K\big)\,\ds \ld,
\end{equation*}

\noi
we see that it suffices to show that there exist $C>0$ and  
 $c>\frac 1p $ such that 
\begin{equation}
\P\big(\|u\|_{L^p(\T)}>\ld, \|u\|_{L^2(\T)}\le K\big)\le C e^{-c\ld^p}
\label{H1}
\end{equation}

\noi
for all sufficiently large $\ld \gg 1$.

Given  $k \in \Z_{\geq 0}$, 
we set 
\[u_k=P_{k} u, \qquad 
u_{\le k}=P_{\le k}u, \qquad \text{and}\qquad 
u_{\ge k}=P_{\ge k}u, \]

\noi
where $P_k$, $P_{\le k}$, and $P_{\ge k}$
are as in \eqref{LP1}, \eqref{LP2}, and \eqref{LP3}.
By subadditivity,  we have for any $k$:
\begin{align}
\begin{split}
\P\big(\|u_{\ge k}\|_{L^p(\T)}>\ld\big) 
& \le \P \bigg(\bigcup_{j = k}^\infty \big\{   \|u_{j}\|_{L^p(\T)}> \ld_j\big\}\bigg)\\
& \le 
\sum_{j=k}^\infty \P( \|u_{j}\|_{L^p(\T)}> \ld_j),
\end{split}
\label{eqn: sacco}
\end{align}

\noi
where $\{\ld_j\}_{j = k}^\infty$ is a sequence of positive numbers such that 
\begin{equation}\label{eqn: grisvold}
\sum_{j=k}^\infty\ld_j = \ld.
\end{equation}

\noi
Then,
by using Sobolev's inequality  in the form of Bernstein's inequality, we have
\begin{equation}\label{eqn: 1stbernstein}
\|u_{j}\|_{L^p(\T)}\le C2^{j(\frac{1}{2}-\frac{1}{p})}\|u_j\|_{L^2(\T)}.
\end{equation}

\noi
Thus, with \eqref{bloop},  the probability 
on the right-hand side of \eqref{eqn: sacco} is bounded by
\begin{align}
\begin{split}
 \P\bigg( \|u_{j}\|_{L^2(\T)}> \frac{\ld_j}{C}2^{j(\frac{1}{p}-\frac{1}{2})} \bigg) 
 &=\P\bigg( \sum_{2^{j-1}\le |n|< 2^j}\frac{|g_n|^2}{n^2}
 > \frac{\ld_j^2}{C_0^2}2^{2j(\frac{1}{p}-\frac{1}{2})}\bigg) \\
&\le \P\bigg( \sum_{2^{j-1}\le |n|< 2^j}|g_n|^2> \frac{\ld_j^2}{4C_0^2}2^{(1+\frac{2}{p})j}\bigg), 
\end{split}
\label{eqn: rebecca}
\end{align}

\noi
where $C_0 =(2\pi)^{-1} C$.

The next lemma follows from a simple calculation involving moment generating functions of 
Gaussian random variables.
See for example  \cite{OQV} for a proof.

\begin{lemma}\label{LEM:Rbound}
Let $\{X_n\}_{n \in \NB}$  be independent standard real-valued Gaussian random variables. 
Then, we have
\begin{equation*} 
\P\bigg(\sum_{n=1}^M X_n^2 \ge R^2\bigg)\le e^{-\frac{R^2}{4}}, 
\end{equation*}

\noi
if $R\ge 3 M^{\frac 12}$.
\end{lemma}

By applying Lemma \ref{LEM:Rbound}, 
we can bound  the probability \eqref{eqn: rebecca}  by 
\begin{equation}\label{eqn: eurus}
\exp\bigg(-\frac{\ld_j^2}{16C_0^2}2^{(1+\frac{2}{p})j}\bigg), 
\end{equation}

\noi
provided
\begin{equation}\label{eqn: tomwaitts}
\ld_j \geq 6  C_0 2^{-\frac{1}{p}j } .
\end{equation}

\noi
By choosing
\[\ld_j = \ld (1-2^{-r})2^{kr}2^{-jr}\]
for $0<r<\frac 1p$, both conditions \eqref{eqn: tomwaitts} and \eqref{eqn: grisvold} are satisfied for all large $k$
(and $j \geq k$). 
For such $k$, the probability in \eqref{eqn: eurus} is then bounded by
\begin{equation*}
\exp\bigg(-\frac{\ld^2(1-2^{-r})^2}{16C_0^2}2^{2kr}2^{(1+\frac{2}{p}-2r)j}\bigg).
\end{equation*}

\noi
Summing over $j\geq k$ in \eqref{eqn: sacco}, we find that
\begin{equation}\label{eqn: vanguard}
\P\big(\|u_{\ge k}\|_{L^p(\T)}>\ld\big)\le C_{r}\exp\bigg(-\frac{\ld^2(1-2^{-r})^2}{16C_0^2}2^{(1+\frac{2}{p})k}\bigg).
\end{equation}

By applying Bernstein's inequality again with the restriction  $\|u\|_{L^2(\T)}\le K$,
we have
\begin{equation}\label{eqn: bertrand}
\|u_{\le k-1}\|_{L^p(\T)}\le C2^{k(\frac{1}{2}-\frac{1}{p})}\|u_{\le k-1}\|_{L^2(\T)}\le C2^{k(\frac{1}{2}-\frac{1}{p})}K .
\end{equation}

\noi
Hence, by  setting
\[k=\log_2 \left(\frac{\ld}{2CK}\right)^{\frac{2p}{p-2}},\]

\noi
it follows from 
\eqref{eqn: bertrand} that
\begin{align}
\|u_{\le k-1}\|_{L^p(\T)}\le \frac{\ld}{2}.
\label{H2}
\end{align}

\noi
Therefore, from  \eqref{eqn: vanguard} and \eqref{H2}, we obtain 
\begin{equation*}
\begin{split}
\P\Big(\| & u\|_{L^p(\T)}  >\ld, \|u\|_{L^2(\T)}\le K\Big)\\
& \le\P\bigg(\|u_{\le k-1}\|_{L^p(\T)}>\frac{\ld}{2}, \|u\|_{L^2(\T)}\le K\bigg)
+ \P\bigg(\|u_{\ge k}\|_{L^p(\T)}>\frac{\ld}{2}\bigg)\\
& \le C_{r}\exp\bigg(-\frac{(1-2^{-r})^2}{64C_0^2\cdot (2CK)^{\frac{2(p+2)}{p-2}}}\ld^{\frac{4p}{p-2}}\bigg)
\end{split}
\end{equation*}

\noi
for all sufficiently large $\ld \gg 1$.
Note that the exponent $\frac{4p}{p-2}$ beats the exponent $p$ in \eqref{H1} 
if (i)  $p<6$ or 
(ii) $p=6$ and $K$ is sufficiently small. Determining the optimal threshold for $K$ would presumably require a delicate optimization of $\ld_j$ in \eqref{eqn: sacco}, an exact Gaussian tail bound to replace the appraisal \eqref{eqn: eurus}, and an optimal inequality to replace the applications of 
Bernstein's inequality in  \eqref{eqn: 1stbernstein} and \eqref{eqn: bertrand} to determine the precise tail behavior of $\|u\|_{L^p(\T)}$ given $\|u\|_{L^2(\T)}\le K$. We did not attempt this calculation. Even if it is possible to carry out, such an approach would likely lead to a less transparent argument than the one we propose in Section \ref{SEC:below}.
Moreover, our argument is easily adapted to the case of the two-dimensional unit disc $\D$.

\section{Sharp Gagliardo-Nirenberg-Sobolev inequality}
\label{SEC:GNS}

The optimizers for the Gagliardo-Nirenberg-Sobolev interpolation inequality
with the optimal constant:
\begin{equation}\label{GNSdisp}
\|u\|_{L^p(\R^d)}^p \le C_{\GNS}(d,p)\|\nb u\|^{\frac{(p-2)d}{2}}_{L^2(\R^d)}\|u\|^{2+\frac{p-2}{2}(2-d)}_{L^2(\R^d)}
\end{equation}

\noi
play an important role in the study of the focusing Gibbs measures.
The following result on the optimal constant $ C_{\GNS}(d,p)$
and optimizers is due to Nagy \cite{nagy} for $d=1$ and Weinstein~\cite{weinstein} for $d\ge 2$.
See also Appendix B in \cite{Tao}.

\begin{proposition}\label{THM:W}
Let $d \geq 1$ and let \textup{(i)} $p > 2$ if $d = 1, 2$ 
and \textup{(ii)} $2 < p < \frac{2d}{d-2}$ if $d \geq 3$.
Consider the functional
\begin{align}
J^{d,p}(u)=\frac{\|\nb u\|^{\frac{(p-2)d}{2}}_{L^2(\R^d)}\|u\|^{2+\frac{p-2}{2}(2-d)}_{L^2(\R^d)}}{\|u\|_{L^p(\R^d)}^{p}}
\label{JJ1}
\end{align}

\noi
on $H^1(\R^d)$. 
Then, the minimum
\begin{align*}
C^{-1}_{\GNS}(d,p):= \inf_{\substack{u\in H^1(\R^d)\\u \ne 0}} J^{d, p}(u)
\end{align*}

\noi
is attained 
at a function  $Q\in H^1(\R^d)$ 
which is a positive, radial, and exponentially decaying  solution to the following semilinear elliptic equation
 on $\R^d$:
\begin{equation}\label{elliptic}
(p-2)d\Delta Q + 4Q^{p-1} -(4+(p-2)(2-d))Q =0
\end{equation}

\noi
with the minimal $L^2$-norm \textup{(}namely, the ground state\textup{)}.
Moreover, 
we have
\begin{equation}\label{GN3}
C_{\GNS}(d,p)= \frac{p}{2}\|Q\|_{L^2(\R^d)}^{2-p}.
\end{equation}

 \end{proposition}

See \cite{frank} for a pleasant exposition, including a proof of the uniqueness of positive solutions to~\eqref{elliptic}, following \cite{kwong}.

\begin{remark}
\label{REM:opt}
 \rm

Recall from Theorem 1.1 in \cite{frank}
that any optimizer $u$ for \eqref{GNSdisp} is 
of the form $u(x) = cQ(b (x - a))$
for some $a \in \R^d$, $b > 0$, and $c \in \C \setminus \{0\}$.
In particular, $u$ is positive 
(up to a multiplicative constant).
\end{remark}

The scale invariance of the minimization problem implies that these inequalities 
also   hold  
 on the finite domains 
$\T$ and $\D$ (essentially) with the same optimal constants.

\begin{lemma}\label{LEM:periodicGN}
\textup{(i)}
Let  $p> 2$.
Then, given any  $m>0$, there is a constant $C=C(m) >0$ such that 
\begin{equation}\label{GNS1}
\|u\|_{L^p(\T)}^p \le \big(C_{\GNS}(1,p)+m\big)\| u'\|_{L^2(\T)}^{\frac{p-2}{2}}\|u\|_{L^2(\T)}^{\frac{p+2}{2}}
+C(m)\|u\|_{L^2(\T)}^p
\end{equation}

\noi
for 
any $u\in H^1(\T)$.
We point out that there exists no $C_0 > 0$ such that the Gagliardo-Nirenberg-Sobolev inequality:
\begin{equation}\label{JJ1a}
\|u\|_{L^p(\T)}^p \le C_0 \| u'\|^{\frac{p-2}{2}}_{L^2(\T)}\|u\|^{\frac{p+2}{2}}_{L^2(\T)}
\end{equation}

\noi
holds
for all functions in $H^1(\T)$.

By restricting our attention to 
mean-zero functions belonging to 
\[\textstyle H^1_0(\T) := \big\{u \in H^1(\T):  \int_\T u \, \d x = 0\big\},\] 

\noi
the 
Gagliardo-Nirenberg-Sobolev inequality \eqref{GNSdisp} holds true on $\T$.
Namely, we have 
\begin{equation}\label{JJ1aa}
\|u\|_{L^p(\T)}^p \le C_{\GNS}(1, p) \| u'\|^{\frac{p-2}{2}}_{L^2(\T)}\|u\|^{\frac{p+2}{2}}_{L^2(\T)}
\end{equation}

\noi
for any function $u \in  H^1_{0}(\T)$.
In fact, 
\eqref{JJ1aa} holds for any $u$ belonging to 
\begin{align}
 H^1_{00}(\T) := \big\{u \in H^1(\T):  
\Re u(x_1) =  \Im u (x_2) = 0 \text{ for some } x_1, x_2 \in \T\big\}.
\label{H00}
\end{align}

\smallskip

\noi
\textup{(ii)} Let $p> 2$.
Then, we have
\begin{equation}\label{GNS2}
\|u\|_{L^p(\D)}^p \le C_{\GNS}(2,p)\|\nb u\|_{L^2(\D)}^{p-2}\|u\|_{L^2(\D)}^2
\end{equation}

\noi
for any $u\in H^1(\D)$ vanishing on $\dd \D$, that is, $u\in H_0^1(\D)$.
\end{lemma}

\begin{proof}
(i) For the proof of \eqref{GNS1},   see Lemma 4.1 in  \cite{LRS}.
As for the failure of the GNS inequality \eqref{JJ1a} for general $u \in H^1(\T)$,  we first note that 
\eqref{JJ1a} does not hold for constant functions.
Moreover, given a function $u \in H^1(\T)$, 
by  simply considering  $u + C$ for large $C\gg 1$
and observing that the left-hand side of
\eqref{JJ1a} on $\T$ grows
faster than the right-hand side as $C \to \infty$, 
we see that  the inequality 
\eqref{JJ1a} on $\T$ does not hold
for (non-constant) functions in general,  unless they have mean zero.
Hence, we restrict our attention to mean-zero   functions.

For 
mean-zero functions on $\T$, 
Sobolev's inequality on $\T$ 
(see \cite{BO})
and  an interpolation yield
the GNS inequality \eqref{JJ1a} with  some constant $C_0 >0$
for any $u \in H^1_0(\T)$.
In fact,   the GNS inequality \eqref{JJ1a} on $\T$
for mean-zero functions holds with 
$C_0 = C_{\GNS}(1,p)$ coming from the GNS inequality \eqref{GNSdisp} on~$\R$.
Suppose that $u \in H^1_0(\T)$ is a real-valued mean-zero function on $\T$.
Then, by the continuity of $u$, there exists a point $x_0 \in \T$ such that $u(x_0) = 0$.
By setting
\[ v(x) = \begin{cases} 
u( x + x_0 - \frac 12), & \quad \text{for } x \in [-\frac 12, \frac 12],\\
0, & \quad \text{for } |x|>\frac12, 
\end{cases}
\]

\noi
we can apply the GNS inequality \eqref{GNSdisp} on $\R$
to $v$ and 
 conclude that the GNS inequality~\eqref{JJ1a} on $\T$ holds for~$u$
 with 
$C_0 = C_{\GNS}(1,p)$.
Now, given complex-valued  $u \in H^1_0(\T)$,
write $u = u_1 + i u_2$, where
$u_1 = \Re u$ and $u_2 = \Im u$.
Since $u$ has mean zero on $\T$, 
its real and imaginary parts also have mean zero.
In particular, there exists $x_j \in \T$ such that $u_j(x_j)=0$, 
$j = 1, 2$.
Hence, by the argument above, 
we see that 
for $u_1$ and $u_2$, the GNS inequality~\eqref{JJ1a}
holds with $C_0 = C_{\GNS}(1,p)$.
We now proceed as in 
Step~2 of the proof of Theorem~A.1 in~\cite{frank}.
By H\"older's inequality (in $j$), we have
\begin{align}
\sum_{j = 1}^2 
\Big(\| u'_j\|^{\frac{p-2}{2}}_{L^2(\T)}\|u_j\|^{\frac{p+2}{2}}_{L^2(\T)}\Big)^\frac{2}{p}
\leq \bigg(\| u'\|^{\frac{p-2}{2}}_{L^2(\T)}\|u\|^{\frac{p+2}{2}}_{L^2(\T)}\bigg)^\frac{2}{p}.
\label{JJ1b}
\end{align}

\noi
Then, by the triangle inequality, 
the GNS inequality \eqref{JJ1a} 
with $C_0 = C_{\GNS}(1,p)$ for $u_j$, $j = 1, 2$,  and \eqref{JJ1b}, we obtain
\begin{align}
\begin{split}
\| u \|_{L^p(\T)}^p 
& = \|u_1^2 + u_2^2\|_{L^\frac{p}{2}(\T)}^\frac{p}{2}
\leq \bigg(\sum_{j = 1}^2 \|u_j\|_{L^p(\T)}^2\bigg)^\frac{p}{2}\\
& \leq C_{\GNS}(1,p)
\| u'\|^{\frac{p-2}{2}}_{L^2(\T)}\|u\|^{\frac{p+2}{2}}_{L^2(\T)}.
\end{split}
\label{JJ1c}
\end{align}

\noi
This proves the GNS inequality \eqref{JJ1aa} 
for $u \in H^1_0(\T)$.
Note that the argument above
shows that  the GNS inequality \eqref{JJ1aa}
on $\T$ 
indeed holds
for any~$u$ belonging to  a larger  class
$H^1_{00}(\T)$ defined in \eqref{H00}.

\smallskip
\noi
(ii) 
Given  a function $u\in H^1(\D)$ vanishing on $\dd \D$,  we can extend $u$ on $\D$ 
to $\bar{u}\in H^1(\R^2)$ 
by setting $\bar u \equiv 0$ on $\R^2\setminus\D$.
Then,  applying \eqref{GNSdisp} on $\R^2$, 
we obtain 
the sharp GNS inequality~\eqref{GNS2} on~$\D$ 
\end{proof}

We conclude this section by stating non-existence
of optimizers
for the  Gagliardo-Nirenberg-Sobolev inequality
\eqref{JJ1aa} on $\T$ among mean-zero functions
and for the  Gagliardo-Nirenberg-Sobolev inequality \eqref{GNS2}
on $\D$ among functions in $H^1(\D)$, vanishing on $\dd \D$.

\begin{lemma}
\label{LEM:opt}

\textup{(i)}
There exists no optimizer in $H^1_0(\T)$ for  the Gagliardo-Nirenberg-Sobolev inequality~\eqref{JJ1aa}
on $\T$.

\smallskip
\noi
\textup{(ii)} 
There exists no optimizer in $H^1(\D)$, vanishing on $\dd \D$,  for  the Gagliardo-Nirenberg-Sobolev inequality~\eqref{GNS2}
on $\D$.


\end{lemma}

\begin{proof}
(i)
Define the functional $J_\T^{1, p}(u)$ by 
\begin{align*}
J_\T^{1,p}(u)
=\frac{\| u'\|^\frac{p-2}{2}_{L^2(\T)}\|u\|^\frac{p+2}{2}_{L^2(\T)}}{\|u\|_{L^p(\T)}^{p}}
\end{align*}

\noi
as in \eqref{JJ1}
and 
 consider the minimization problem
 over $ H^1_{00}(\T)$:
\begin{align}
 \inf_{\substack{u\in H^1_{00}(\T)\\u \ne 0}} J_\T^{1, p}(u).
\label{JJ2}
\end{align}

\noi
It follows from \eqref{JJ1c}  that this infimum is bounded below by 
$ C_{\GNS}(1,p)^{-1}>0$.
Suppose that there exists a mean-zero  optimizer $u_* 
\in H^1_{0}(\T) \subset  H^1_{00}(\T)$
for \eqref{JJ2}.
Then, by adapting the argument in  Step~2 of the proof of Theorem A.1 in \cite{frank}
(see \eqref{JJ1b} and \eqref{JJ1c} above)
to the case of the one-dimensional torus $\T$, 
we see that either (i)~one of $\Re u_*$ or $\Im u_*$
is identically equal to 0, 
or (ii)~both $\Re u_*$ and $\Im u_*$
are optimizers for \eqref{JJ2}
and $|\Re u_*|= \ld|\Im u_*|$ for some $\ld > 0$.
In either case, 
one of $\Re u_*$ or $\Im u_*$ is an  optimizer for \eqref{JJ2}.
Without loss of generality, 
suppose that 
$\Re u_*$ is a (non-zero) optimizer for \eqref{JJ2}.
Note that 
$\Re u_* \in H^1_{00}(\T)$
and that its positive and negative parts
also belong to $H^1_{00}(\T)$
(but not to $H^1_{0}(\T)$).
Then, by the argument in Step~2 of the proof of Theorem~A.1 in~\cite{frank}
once again, 
we conclude that 
$\Re u_*$ is either non-negative or non-positive.
This, however, is a contradiction 
to the fact that $u$ (and hence $\Re u_*$) has mean zero on $\T$.
This argument shows that there
 is no mean-zero optimizer
for the GNS inequality  \eqref{JJ1aa} on $\T$.

\smallskip

\noi
(ii)
Suppose that there exists
an optimizer $u$ in $H^1(\D)$, vanishing on $\dd \D$.
We can   extend $u$ on $\D$
to a function in $ H^1(\R^2)$ 
by setting $ u \equiv 0$ on $\R^2\setminus\D$, 
which would be a (non-zero) non-negative
 optimizer for \eqref{GNSdisp} on $\R^2$ with compact support, 
which is a contradiction (see Theorem~1.1 in \cite{frank}).
\end{proof}

\section{Integrability below the threshold}
\label{SEC:below}

\subsection{On the one-dimensional torus}\label{SUBSEC:1d}
We first  present the proof of 
 the direct implication in Theorem~\ref{THM:LRS}\,(ii). 
 Namely, we show that the partition function $Z_{6, K}$
 in \eqref{Z1} is finite, 
 provided that $K < \|Q\|_{L^2(\R)}$.
  See \cite[Theorem 2.2 (b)]{LRS} for the converse.
All the norms are taken over the one-dimensional torus $\T$, 
unless otherwise stated.

In the following, 
$\E$ denotes an expectation with respect to 
the mean-zero  Brownian loop
$u$ in \eqref{bloop}
(in particular $u$ has mean-zero on $\T$).
Given  $\ld>0$,  we use the notation \eqref{EXP1} and write
\begin{equation*}
\begin{split}
Z_{p, K} & = \E\Big[e^{\frac{1}{p}\int_\T |u(x)|^p\,\ds  x}, 
\|u\|_{L^2}\le K\Big]\\
&= \E\Big[e^{\frac{1}{p}\int_\T |u(x)|^p\,\ds  x}, \|u_{\ge 0}\|_{L^p}\le \ld, \|u\|_{L^2}\le K\Big]\\
&\hphantom{X}+ \E\Big[e^{\frac{1}{p}\int_\T |u(x)|^p\,\ds  x}, \|u_{\ge 0}\|_{L^p}> \ld,  \|u\|_{L^2}\le K\Big].
\end{split}
\end{equation*}

\noi
Here,  we used the fact that $P_{\geq 0} = \text{Id}$ on mean-zero functions.

For $k \geq 1$, define  $E_k$ by 
\begin{align}
\begin{split}
E_k&=\big\{\|u_{\geq 0} \|_{L^p}>\ld, \ldots, \|u_{\ge k-1}\|_{L^p}>\ld, \|u_{\ge k}\|_{L^p}\le \ld \big\}\\
&\subset \big\{\|u_{\ge k-1}\|_{L^p}>\ld, \|u_{\ge k}\|_{L^p}\le \ld \big\}.
\end{split}
\label{H2ab}
\end{align}

\noi
Note that the sets $E_k$'s are disjoint
and that  $\sum_{k=1}^N \mathbf{1}_{E_k}$ increases to $ \ind_{\{\|u\|_{L^p}>\ld\}}$
almost surely as $N \to \infty$ since $u$ in \eqref{bloop}
belongs almost surely to $L^p(\T)$ for any finite $p$.
Hence, by the  monotone convergence theorem, we obtain
\begin{equation}\label{eqn: livius}
\begin{split}
Z_{p, K} 
&= \E\Big[e^{\frac{1}{p}\int_\T |u(x)|^p\,\ds  x}, \|u\|_{L^p}\le\ld, \,\|u\|_{L^2}\le K\Big]\\
&
\hphantom{X}
+\sum_{k=1}^\infty \E\Big[e^{\frac{1}{p}\int_\T |u(x)|^p\,\ds  x}, E_k, \,\|u\|_{L^2}\le K\Big].
\end{split}
\end{equation}

The first term on the right-hand side of \eqref{eqn: livius} is clearly finite
for any finite $\ld, K > 0$.
Hence,  in view of \eqref{H2ab}, it suffices to show that
\begin{equation}\label{eqn: andre}
\E\Big[e^{\frac{1}{p}\int_\T |u(x)|^p\,\ds  x}, \|u_{\ge k-1}\|_{L^p}>\ld, \|u_{\ge k}\|_{L^p}\le \ld, \,\|u\|_{L^2}\le K\Big]
\end{equation}

\noi
is summable in $k \in \NB$.

Given an integer $p$, we have
\begin{align*}
|u|^p = |u_{\le k-1}+u_{\ge k}|^p&\le \sum_{\l=0}^p \binom{p}{\l} |u_{\le k-1}|^{p-\l}|u_{\ge k}|^\l.
\end{align*}

\noi
Integrating and applying  H\"older's inequality
followed by Young's inequality, 
we have,  for any $u$ satisfying $\|u_{\ge k}\|_{L^p}\le \ld$, 
\begin{align}\label{eqn: ezekiel}
\begin{split}
\int_\T|u|^p\, \ds  x 
& \le \int_\T|u_{\le k-1}|^p\,\ds  x+\sum_{\l=1}^p \binom{p}{\l}  \|u_{\le k-1}\|_{L^p}^{p-\l} \ld^{\l}\\
& \le \int_\T|u_{\le k-1}|^p\,\ds  x+\sum_{\l=1}^p \binom{p}{\l}  
\bigg( \frac{p-\l}{p}\eps \|u_{\le k-1}\|_{L^p}^p+ \frac{\l}{p} \eps^{- \frac{p-\l}{\l}}\ld^p\bigg)\\
& \le 
\big(1+ (2^p-1)\eps\big)
\int_\T|u_{\le k-1}|^p\,\ds  x+  C_p( \eps)\ld^p
\end{split}
\end{align}

\noi
for some small $\eps > 0$ (to be chosen later), 
where in the last step we used
\[\sum_{\l = 1}^p\frac{p-\l}{p}\binom p\l
\le \sum_{\l = 1}^p\binom p\l
\leq 2^p -1.\]

\noi
Hence,  by letting
\begin{equation}
\dl =\dl(p, \eps) = (2^p-1)\eps, 
\label{H2a}
\end{equation}

\noi
the quantity in \eqref{eqn: andre} is bounded by
\begin{equation}\label{eqn: titus}
e^{C_p(  \eps) \ld^p }
 \E\Big[e^{\frac{(1+\dl)}{p}\int_\T |u_{\le k-1}(x)|^p\, \ds  x}, \|u_{\ge k-1}\|_{L^p}>\ld ,\|u\|_{L^2}\le K\Big].
\end{equation}

Now, let $\ld = 1$ and $p = 6$. 
By Lemma \ref{LEM:periodicGN}\,(i)
 for some small $m>0$ (to be chosen later)
 with $\|u\|_{L^2}\leq K$, there is a constant $C(m)>0$ such that 
 \eqref{eqn: titus} is  now bounded by
\begin{align}
& e^{C_6( \eps)   + C(m)K^6}\E\Big[e^{\frac{(C_{\GNS}(1,6)+m)K^4(1+\dl)}{6}\int_\T |u_{\le k-1}'(x)|^2\,\ds  x},\|u_{\ge k-1}\|_{L^6}>1 \Big] \notag 
\intertext{By H\"older's inequality, }
& \leq e^{C_6(\eps)   + C(m)K^6}\bigg\{\E\Big[e^{\frac{(C_{\GNS}(1,6)+m)K^4(1+\eta)(1+\dl)}{6}\int_\T |u_{\le k-1}'(x)|^2\,\ds  x}
\Big]\bigg\}^{\frac{1}{1+\eta}}\notag \\
& \hphantom{X}
\times \Big\{\P\big(\|u_{\ge k-1}\|_{L^6}> 1\big)\Big\}^{\frac{\eta}{1+\eta}}.
\label{eqn: haddock}
\end{align}

%

By  \eqref{eqn: vanguard}, we have
\begin{equation}
\Big\{\P\big(\|u_{\ge k-1}\|_{L^6}> 1\big)\Big\}^{\frac{\eta}{1+\eta}}
\le C^{\frac{\eta}{1+\eta}}\exp\left(-C\frac{\eta}{1+\eta}2^{\frac 43k}\right).
\label{H3}
\end{equation}

\noi
Recall that
\begin{align}
\E[e^{tX^2}] = (1-2t)^{-\frac{1}{2}}
\label{H3a}
\end{align}
 
\noi
for $X \sim \N_\R(0, 1)$
and 
$t < \frac 12$.
Then, using \eqref{H3a} 
and \eqref{GN3} in Proposition~\ref{THM:W}, 
 the expectation in~\eqref{eqn: haddock} is bounded by
\begin{align}
\begin{split}
\E & \Big[  e^{\frac{(C_{\GNS}(1,6)+m)K^4(1+\eta)(1+\dl)}{6}\int_\T |u_{\le k-1}'(x)|^2\,\ds  x}\Big] \\
&= \prod_{1\le |n| \le 2^{k-1}} \E\Big[ e^{ \frac{(C_{\GNS}(1,6)+m)K^4(1+\eta)(1+\dl)}{6}|g_n|^2}\Big]\\
 &\le \bigg(1-2\frac{(C_{\GNS}(1,6)+m) K^4 (1+\eta)(1+\dl)}{6}\bigg)^{-2^{k}}\\
&= \bigg( 1 - \frac{K^4}{\|Q\|_{L^2(\R)}^4}\Big(1+\frac{m}{3}\|Q\|_{L^2(\R)}^4\Big)
(1+\eta)(1+\dl)\Bigg)^{-2^{k}}.
\end{split}
\label{H4}
\end{align}

\noi
Note that, under $K < \|Q\|_{L^2(\R)}$ and \eqref{H2a}, 
 we can choose $m, \eps, \eta > 0$ sufficiently small such that
\[\frac{K^4}{\|Q\|_{L^2(\R)}^4}\Big(1+\frac{m}{3}\|Q\|_{L^2(\R)}^4\Big)(1+\eta)(1+\dl)< c<1,\]

\noi
guaranteeing the application of \eqref{H3a} in the computation above.

Finally, summing \eqref{eqn: andre} over $k \in \NB$
with \eqref{eqn: titus}, \eqref{eqn: haddock}, 
\eqref{H3}, and \eqref{H4}, we have
\begin{align}
\begin{split}
\sum_{k = 1}^\infty
&  \E\Big[e^{\frac{1}{6}\int_\T |u(x)|^6\,\ds  x}, \|u_{\ge k-1}\|_{L^6}>\ld, \|u_{\ge k}\|_{L^6}\le \ld, \,\|u\|_{L^2}\le K\Big]\\
& \leq e^{C(\eps, m, K) }
C^{\frac{\eta}{1+\eta}}
\sum_{k = 1}^\infty
\exp\bigg(-C\frac{\eta}{1+\eta}2^{\frac{4}{3}k}\bigg)
\exp\bigg(\frac{2^k}{1+\eta}\log\frac{1}{1-c}\bigg)\\
& < \infty.
\end{split}
\label{H4a}
\end{align}

\noi
Therefore, we conclude that 
 the partition function $Z_{6, K}$ is finite for any $K<\|Q\|_{L^2(\R)}$.
This completes the proof of Theorem \ref{THM:LRS}.

\begin{remark}\rm
\label{REM:gauss}

As mentioned in Remark \ref{REM:loop}, 
Theorem \ref{THM:LRS}\,(ii)  holds
for the Ornstein-Uhlenbeck loop 
in  \eqref{bloop3}.
We first note that
  \eqref{eqn: vanguard}
  also holds
for the Ornstein-Uhlenbeck loop  $u$
since \eqref{eqn: rebecca}  holds
even if we replace $n^2$ by $(2\pi)^{-2}+  n^2$.
Hence, from \eqref{bloop3} and~\eqref{H3}, we have 
\begin{align*}
\eqref{eqn: haddock}
& \leq e^{C_6(\eps)   + C(m)K^6}\bigg\{\E_\mu \Big[e^{\frac{(C_{\GNS}(1,6)+m)K^4(1+\eta)(1+\dl)}{6}
\int_\T |\jb{\dx}u_{\le k-1}(x)|^2\,\ds  x}
\Big]\bigg\}^{\frac{1}{1+\eta}}\notag \\
& \hphantom{X}
\times \Big\{\P\big(\|u_{\ge k-1}\|_{L^6}> 1\big)\Big\}^{\frac{\eta}{1+\eta}}\\
& 
\leq 
C^{\frac{\eta}{1+\eta}}\exp\left(-C\frac{\eta}{1+\eta}2^{\frac 43k}\right)
\prod_{0\le |n| \le 2^{k-1}} \E\Big[ e^{ \frac{(C_{\GNS}(1,6)+m)K^4(1+\eta)(1+\dl)}{6}|g_n|^2}\Big], 
\end{align*}

\noi
where $\E_\mu$ denotes an expectation
with respect to the law of  the Ornstein-Uhlenbeck loop $u$
and $\jb{\dx} = \sqrt{1 - \dx^2}$.
The rest of the argument follows as above, 
thus establishing 
Theorem~\ref{THM:LRS}\,(ii)
for the Ornstein-Uhlenbeck loop $u$ in \eqref{bloop3}.

\end{remark}

\subsection{On the two-dimensional disc}
\label{SUBSEC:2d}

Next, we 
prove 
the normalizability of the Gibbs measure on $\D$ stated in Theorem \ref{THM:BB}\,(ii).
Namely, we show that the partition function $\wt Z_{4, K}$
 in \eqref{Z3} is finite, 
 provided that $K < \|Q\|_{L^2(\R^2)}$.
 The proof is based on a computation analogous to that in Subsection \ref{SUBSEC:1d}.
 As we see below, however, 
 we need to proceed with more care, partially due to the eigenfunction estimate,
 which makes the computation barely work on $\D$; 
 compare \eqref{H4a} and \eqref{H4b}.

We first recall  the following simple corollary of  Fernique's theorem \cite{fernique}.
See also Theorem~2.7 in \cite{DZ}.\footnote
{In the context of Theorem 2.7 in \cite{DZ}, 
we set $x = \frac{X}{\E[\|X\|_B]}$.
Then, by Markov's inequality
and choosing $r\gg 1$, we have
\[ \log \bigg(\frac{\mu(\|x\|_B > r)}{\mu(\|x\|_B \leq r)} \bigg)
= \log \bigg(\frac{\mu\big(\|X\|_B > r\E[\|X\|_B]\big)}{1- \mu\big(\|X\|_B > r\E[\|X\|_B]\big)} \bigg)
\leq \log \frac{1}{r-1} \leq -2\]

\noi
without using any fine property of $X$.
Then,  \eqref{exp1} and \eqref{exp2} follow in view of 
Remark 2.8 in \cite{DZ}.
}

\begin{lemma}\label{LEM:Fer}
There exists a  constant $c>0$ such that if $X$ is a mean-zero Gaussian process   with values 
in a separable Banach space $B$ with $\E\big[\|X\|_{B}\big]<\infty$, then
\begin{align}
\int e^{ c \frac{\|X\|_B^2}{(\E[\|X\|_B])^2}}\,\ds \P <\infty.
\label{exp1}
\end{align}

\noi
In particular, we have
\begin{equation}
\P\Big(\|X\|_B \ge t\E\big[\|X\|_B\big]\Big)\le e^{-ct^2}
\label{exp2}
\end{equation}

\noi
for any $t>1$.
\end{lemma}

Recall from \cite[Lemmas 2.1 and 2.2]{Tzv1} 
the asymptotic formula for the eigenvalue $z_n$:
\begin{align}
z_n = \pi\bigg(n - \frac 14\bigg) + O\bigg(\frac{1}{n}\bigg)
\label{H5}
\end{align}

\noi
and 
 the eigenfunction estimate: 
\begin{align}
\|e_n \|_{L^4(\D)} \les \big[\log(2 + n)\big]^\frac{1}{4}.
\label{H6}
\end{align}

\noi
As in Section \ref{SEC:Bourgain}, 
we define the spectral projection of $v$ by 
\[v_k:= \sum_{2^{k-1}<  n\leq  2^k} \ft{v}(n)e_n\]

\noi
for $k \in \Z_{\geq 0}$, where $\ft v(n) = \int_\D v e_n \ds x$.
We also define $v_{\leq k}$ and $v_{\geq k}$ in an analogous manner.

In the following, $\E$ denotes an expectation with respect to 
the random Fourier series $v$ in \eqref{bloop2}
and all the norms are taken over the unit disc $\D \subset \R^2$
unless otherwise stated.
By Minkowski's integral inequality with~\eqref{H5} and~\eqref{H6}, 
we have
\begin{align}
\begin{split}
\E\big[\|v_j \|_{L^4}\big]
& \les 
\bigg\|\bigg( \sum_{ 2^{j-1}<  n \le 2^j} \frac{1}{{z_n}^2}e_n^2\bigg)^\frac{1}{2}\bigg\|_{L^4}\\
& \le 
\bigg( \sum_{ 2^{j-1}<  n \le 2^j} \frac{1}{{z_n}^2}\|e_n\|_{L^4}^2\bigg)^\frac{1}{2}
\les  \jb{j}^\frac {1}{4}2^{-\frac{1}{2}j}
\end{split}
\label{H7}
\end{align}

\noi
for any  $j \in \Z_{\geq 0}$.
Then, applying Lemma \ref{LEM:Fer} and \eqref{H7}
 with suitable $\eps_j \sim \jb{j}^{-2}$ such that $\sum_{j\ge k}\eps_j \leq  1$, we have
\begin{equation}
\label{eqn: jessica} 
\begin{split}
\P\big( \|v_{\ge k}\|_{L^4}\ge \ld\big)
 &\le \sum_{j\ge k}\P\big(\|v_j\|_{L^4} \ge \eps_j \ld\big)\\
 &=  \sum_{j\ge k}\P
 \bigg(\|v_j\|_{L^4} \ge \frac{\eps_j \ld}{\E\big[\|v_j \|_{L^4}\big]}\E\big[\|v_j \|_{L^4}\big]\bigg)\\
&\le \sum_{j\ge k} e^{-c \eps_j^2 \ld^2 \jb{j}^{-\frac 12} 2^{j}}\\
&\le Ce^{-c' \ld^2 \jb{k}^{-\frac 92}  2^{k}}
\end{split}
\end{equation}

\noi
for some constant $C>0$, uniformly in $k \in \Z_{\geq 0}$ and $\ld \geq 1$.
Then, it follows from  the Borel-Cantelli lemma
that 
\begin{align}
\limsup_{k \to \infty} \frac{\|v_{\geq k}\|_{L^4}}{\jb{k}^3} = 0
\label{H7a}
\end{align}

\noi
with probability 1.

Define a set $F_k$ by 
\begin{align*}
F_k&=\big\{\|v_{\geq j} \|_{L^4}> \jb{j}^3, 
\ j = 0, 1, \dots, k-1, \quad \text{and}\quad 
\|v_{\ge k}\|_{L^p}\le \jb{k}^3  \big\}.
\end{align*}

\noi
By definition,  $F_k$'s are disjoint
and, from \eqref{H7a}, we have
\[ \P\bigg( \bigcup_{k = 1}^\infty F_k\bigg) = 1.\]

\noi
This implies 
that   $\sum_{k=1}^N \mathbf{1}_{F_k}$ increases to $ \ind_{\{\|v\|_{L^4}>1\}}$
almost surely as $N \to \infty$.

Starting from $\wt Z_{4,K}$ in \eqref{Z3}, we reproduce the computations in 
\eqref{eqn: livius}, \eqref{eqn: ezekiel}, 
and \eqref{eqn: titus} with $p=4$ and $\ld = 1$
by replacing 
$u$ in \eqref{bloop},  
the sets $E_k$,  
and  the integrals over $[0,1]$
with $v$ in~\eqref{bloop2}, 
the sets $F_k$, and  
integrals over $\D$, respectively.
 We find
\begin{align}\label{eqn: baptiste}
\begin{split}
\wt Z_{4,K}
& \leq 
 \E\Big[e^{\frac{1}{4}\int_\D |v(x)|^4\,\ds  x}, \|v\|_{L^4}\le 1, \,\|v\|_{L^2}\le K\Big]\\
& \hphantom{X}
+ \sum_{k=1}^\infty 
\E\Big[e^{\frac{1}{4}\int_{\D} |v(x)|^4\, \ds  x}, \|v_{\ge k-1}\|_{L^4}>\jb{k-1}^3 ,\\
& \hphantom{XXXXXXXXXXXil}
\|v_{\ge k}\|_{L^4}\leq \jb{k}^3, 
\|v\|_{L^2}\le K\Big]\\
& \leq 
C 
+ \sum_{k=1}^\infty e^{C_4(\eps) \jb{k}^{12} }
\E\Big[e^{\frac{(1+\dl)}{4}\int_{\D} |v_{\le k-1}(x)|^4\, \ds  x}, \\
& \hphantom{XXXXXXXXXXXil}
\|v_{\ge k-1}\|_{L^4}>\jb{k-1}^3 ,\|v\|_{L^2}\le K\Big],
\end{split}
\end{align}
where
$\dl=\dl(4)=48\eps$ is as in \eqref{H2a} with $p = 4$.
As before, it remains to show that the series in~\eqref{eqn: baptiste} is convergent.

From Lemma \ref{LEM:periodicGN}\,(ii), we have
\begin{equation*}
\|v\|^4_{L^4(\D)}\le C_{\GNS}(2,4)\|v\|_{L^2(\D)}^2\|\nb v\|_{L^2(\D)}^2
\end{equation*}

\noi
for any $v\in H^1_0(\D) := \{ u \in H^1(\D): u \equiv 0 \text{ on }  \dd \D\}$
with 
\[C_{\GNS}(2,4)=2\|Q\|^{-2}_{L^2(\R^2)},\]

\noi
where 
 $Q$ is the optimizer for the Gagliardo-Nirenberg-Sobolev inequality \eqref{GNSdisp} on $\R^2$.
Then, by recalling
\begin{align*}
\int_{\D} |\nb v_{\le k-1}(x)|^2\,\ds  x &= -\int_{\D} v_{\le k-1}(x) \cj{\Delta v_{\le k-1}(x)}\,\ds  x\\
&= -\sum_{n\le 2^{k-1}} \frac{|g_n|^2 }{z_n^2} \int_0^1 e_n(r)\Delta_r e_n(r)r\,\ds  r\\
&=\sum_{n\le 2^{k-1}} |g_n|^2
\end{align*}

\noi
and applying 
H\"older's inequality, 
the expectation in the summands on the right-hand side of \eqref{eqn: baptiste} is bounded by
\begin{equation}\label{eqn: cassandra}
\begin{split}
&\E\Big[e^{\frac{C_{\GNS}(2,4)K^2(1+\dl)}{4}\int_{\D}| \nb v_{\le k-1}(x)|^2\,\ds  x},
\|v_{\ge k-1}\|_{L^4}>\jb{k-1}^3 \Big]\\
& \hphantom{X}
=  \E\Big[e^{\frac{C_{\GNS}(2,4)K^2(1+\dl)}{4}\sum_{n\le 2^{k-1}}|g_n|^2} , \|v_{\ge k-1}\|_{L^4}> 
\jb{k-1}^3\Big]\\
& \hphantom{X}
 \leq \bigg\{\E\Big[e^{\frac{C_{\GNS}(2,4)K^2(1+\eta)(1+\dl)}{4}\sum_{n\le 2^{k-1}}|g_n|^2}\Big]\bigg\}^{\frac{1}{1+\eta}}\Big\{\P\big( \|v_{\ge k-1}\|_{L^4}> \jb{k-1}^3 \big)\Big\}^{\frac{\eta}{1+\eta}}.
\end{split}
\end{equation}

\noi
The first factor on the right-hand side of \eqref{eqn: cassandra}
can be computed exactly as in \eqref{H4}, using~\eqref{H3a}.
Namely,   provided that
\begin{align}
\label{eqn: mathaeus}
\begin{split}
C(K,\dl,\eta) & :  =C_{\GNS}(2,4) \frac{K^2}{2}(1+\eta)(1+\dl) \\
&\hphantom{:}
 =\frac{K^2}{\|Q\|_{L^2(\R^2)}^2}(1+\eta)(1+\dl)
<1, 
\end{split}
\end{align}

\noi
it is finite and equals
\begin{align}
\bigg(1-\frac{K^2}{\|Q\|_{L^2(\R^2)}^2}(1+\eta)(1+\dl)\bigg)^{-2^{k-1}}.
\label{H8}
\end{align}

Given  $K<\|Q\|_{L^2(\R^2)}$,  we choose $\dl$ and $\eta$ sufficiently small 
such that  \eqref{eqn: mathaeus} holds. 
Then, 
from 
\eqref{eqn: cassandra}
with 
 \eqref{eqn: jessica} and \eqref{H8}, 
 we conclude  that 
$\wt Z_{4, K}$ in \eqref{eqn: baptiste}
is bounded by  
 \begin{align}
 \begin{split}
\wt Z_{4,K}
& \leq 
C + \sum_{k=1}^\infty e^{C_4(\eps) \jb{k}^{12} }\\
& \hphantom{XXXX}
\times 
\exp\bigg(\frac{2^{k-1}}{1+\eta} \log \frac{1}{1-C(K,\dl,\eta)}\bigg)\exp\bigg(-\frac{c''\eta}{1+\eta}
\jb{k}^\frac{3}{2}2^k\bigg)\\
& < \infty.
\end{split}
\label{H4b}
\end{align}

\noi
This proves the normalizability of the Gibbs measure on $\D$ claimed in Theorem \ref{THM:BB}\,(ii).

\section{Non-integrability above the threshold on the disc}\label{SEC:non}

In this section, we  discuss the non-normalizability 
of the Gibbs measure on $\D$ stated
in Theorem \ref{THM:BB} (ii).
Namely, when 
 $K> \|Q\|_{L^2(\R^2)}$, 
 we prove  
  \begin{align}
\wt Z_{4,K}= \E\Big[e^{\frac{1}{4}\int_{\D}|v|^4\,\ds  x}\mathbf{1}_{\{\|v\|_{L^2(\D)}\le K\}}\Big]
= \infty.
\label{BB1}
\end{align}

\noi
Here, $Q$ is the ground state on $\R^2$ as in Proposition~\ref{THM:W}.
Since 
the proof of \eqref{BB1} is essentially identical to  that 
for the non-normalizability of the Gibbs measure on the torus $\T$
(see \cite[Theorem~2.2~(b)]{LRS}), 
we keep our presentation brief.

Let $K_0 = \|Q\|_{L^2(\R^2)}$
and fix $K > K_0$.
Choose  $\al > 1$ such that 
\begin{align}
\| \al Q\|_{L^2(\R^2)} < K.
\label{BB1a}
\end{align}

\noi
Then, by setting
\begin{align}
H_{\R^2}(v) := \frac 12\int_{\R^2} |\nb v|^2 \ds x - \frac 14\int_{\R^2} |v|^4 \ds x, 
\label{BB2}
\end{align}

\noi
it follows from Proposition~\ref{THM:W}
(in particular 
\eqref{GN3})
that 
 $H_{\R^2}(Q) = 0$.
As a result, we have
\begin{align*}
H_{\R^2} (\al Q) < 0. 
\end{align*}

Given $\rho > 0$, define the $L^2$-invariant scaling operator $D_\ld$ on $\R^2$
by setting
\[ D_\ld (f) = \ld^{-1} f(\ld^{-1}x).\]

\noi
Then, we have
\begin{align}
\begin{split}
\| D_{\rho^{-1}} (f) \|_{L^2(\R^2)}^2 & =  \|f\|_{L^2(\R^2)}^2, \\
\| D_{\rho^{-1}} (f) \|_{L^4(\R^2)}^4 & = \rho^2 \|f\|_{L^4(\R^2)}^4, \\
\| \nb D_{\rho^{-1}} (f) \|_{L^2(\R^2)}^2 & = \rho^2 \|\nb f\|_{L^2(\R^2)}^2,\\
\| \dd_r^2 D_{\rho^{-1}} (f) \|_{L^1(\R^2)} & = \rho \|\dd_r^2 f\|_{L^1(\R^2)}, 
\end{split}
\label{BB4}
\end{align}

\noi
where we assume $f$ to be radial for the last identity
and $\dd_r$ denotes the directional derivative in the radial direction.

For each $\rho \gg 1$, 
we define  $Q_{\al, \rho} \in C^\infty_{\text{rad}}(\D)$
by setting
\[ Q_{\al, \rho}(x) =  \al D_{\rho^{-1}}(Q)(x) = \al \rho Q(\rho x)\]

\noi
for $|x| \leq 1 - \frac {1}{\rho}$
and  $Q_{\al, \rho}  \equiv 0$ on $\dd \D$.
Then, thanks to the exponential decay of the ground state $Q$ on $\R^2$, 
it follows from \eqref{BB1a}, \eqref{BB2}, and \eqref{BB4} that 
\begin{align}
\begin{split}
H_{\D} (Q_{\al, \rho}) 
 : = \frac 12\int_{\D} |  \nb Q_{\al, \rho}|^2 & \ds x   - \frac 14\int_{\D} |Q_{\al, \rho}|^4 \ds x
\le - A_1 \rho^2 < 0, \\
\| Q_{\al, \rho} \|_{L^2(\D)}^2 & < K^2 - \eta_0, \\
\| Q_{\al, \rho} \|_{L^4(\D)}^4 & \le A_2 \rho^2, \\
\| \dd_r^2 Q_{\al, \rho} \|_{L^1(\D)} & \leq A_3 \rho 
\end{split}
\label{BB4a}
\end{align}

\noi
for some $A_1, A_2, A_3 > 0$ and sufficiently small $\eta_0 > 0$, 
uniformly in large $\rho \gg1 $.

We need  the following simple calculus lemma
which follows from a straightforward computation.  
See Lemma 4.2 in \cite{LRS}.

\begin{lemma}\label{LEM:BB1}
Given $\eta > 0$, there exists small $\eps > 0$ such that 
for any $v_1, v_2 \in C_0(\D)$ satisfying
$\| v_1 -  v_2\|_{L^\infty(\D)} < \eps $, 
we have 
\[ \Big| \| v_1\|_{L^p(\D)}^p
- \| v_2\|_{L^p(\D)}^p\Big| < \eta \big( \| v_1\|_{L^p(\D)}^p + 1\big)\]

\noi
for $p = 2, 4$.
Here, 
$C_0(\D)$ denotes the collection of continuous complex-valued functions on $\D$, 
vanishing on the boundary $\dd \D$.

\end{lemma}

Given  $\eps > 0$, 
let $B_\eps (Q_{\al, \rho}) \subset C_0(\D)$
be the ball of radius $\eps$ centered at $Q_{\al, \rho}$:
\[B_\eps (Q_{\al, \rho}) = \big\{ v  \in C_0(\D): \| v - Q_{\al, \rho}\|_{L^\infty(\D)} < \eps \big\}.\]

\noi
Let $ \eta \leq \frac{\eta_0}{K^2 + 1-\eta_0}$.
Then, from   Lemma \ref{LEM:BB1} and \eqref{BB4a}, 
there exists small $\eps > 0$ such that 
  \begin{align}
\begin{split}
\wt Z_{4,K}
& = \int_{C_{\text{rad}, 0}(\D)}e^{\frac{1}{4}\int_{\D}|v(x)|^4\,\ds  x}\mathbf{1}_{\{\|v\|_{L^2(\D)}\le K\}}
 \P(\ds v)\\
& \geq 
\int_{B_\eps(Q_{\al, \rho})}
e^{\frac{1}{4}\int_{\D}|v(x)|^4\,\ds  x}\, \P(\ds v)\\
& \geq 
\exp\bigg( \frac{1-\eta}{4}\int_{\D}|Q_{\al, \rho}(x)|^4\,\ds  x - \frac{\eta}{4}\bigg)
\P \big(B_\eps(Q_{\al, \rho})\big), 
\end{split}
\label{BB5}
\end{align}

\noi
where 
$C_{\text{rad}, 0}(\D) =   \big\{ v\in C_0(\D): v,  \,\text{radial}\big\}$
and 
$\P$ denotes the Gaussian probability measure with respect to the random series \eqref{bloop2}.

Since $Q_{\al, \rho} \in C^\infty_{\text{rad}, 0}(\D) \subset H^1_{\text{rad, 0}}(\D)
=   \big\{ v\in H^1_0(\D):  v,\,  \text{radial}\big\}$, 
we can apply the Cameron-Martin theorem
(Theorem 2.8 in \cite{DP})
and obtain 
\begin{align}
\begin{split}
\P \big(B_\eps(Q_{\al, \rho})\big)
& = 
e^{-\frac 12 \int_{\D} |\nb Q_{\al, \rho}(x)|^2 \ds x}
\int_{B_\eps(0)} 
  e^{\Re \int_\D v \Dl  Q_{\al, \rho} \ds x}\, \P(\ds v)\\
& \geq   
e^{-\frac 12 \int_{\D} |\nb Q_{\al, \rho}(x)|^2 \ds x}
  e^{-\eps \| \dd_r^2  Q_{\al, \rho}\|_{L^1(\D)}} \P(B_\eps(0))\\
& \geq   
e^{-\frac 12 \int_{\D} |\nb Q_{\al, \rho}(x)|^2 \ds x}
  e^{-\eps A_3 \rho } \P\big(B_\eps(0)\big), 
\end{split}
\label{BB6}
\end{align}

\noi
where we used \eqref{BB4a} in the last step.

Hence, by further imposing
$\eta \leq \frac{2 A_1}{A_2}$, 
we conclude from  \eqref{BB5} and \eqref{BB6}
with $ \P\big(B_\eps(0)\big) = C_\eps >  0$ and \eqref{BB4a} that 
  \begin{align*}
\wt Z_{4,K}
& \geq C_\eps 
\exp\bigg( 
- H_\D(Q_{\al, \rho})
- \frac{\eta}{4}\|Q_{\al, \rho}\|_{L^4(\D)}^4 
-\eps  A_3\rho
- \frac{\eta}{4}\bigg)\\
& \geq C_\eps 
\exp\bigg( 
 A_1 \rho^2 
- \frac{\eta}{4} A_2 \rho^2
-\eps  A_3\rho
- \frac{\eta}{4}\bigg)\\
& \too \infty, 
\end{align*}

\noi
as $\rho\to \infty$.
This proves the non-normalizability of the Gibbs measure on $\D \subset \R^2$,  
when  $K> \|Q\|_{L^2(\R^2)}$.

\section{Integrability at the optimal mass threshold}\label{SEC:opt}

\numberwithin{equation}{subsection}

In this section, we present the proof of Theorem \ref{THM:OPT}.
Namely, we show that   $Z_{6,K}<\infty$ when
\[K=\|Q\|_{L^2(\R)}.\]

\noi
We fix 
$d = 1$,  $p = 6$, and  
$C_{\GNS} = C_{\GNS}(1, 6)$.
In this section, we prove the normalizability
when $u$ is 
 the Ornstein-Uhlenbeck loop  in \eqref{bloop3}:
 \begin{equation}
\label{Q0}
u(x)=\sum_{n \in \Z}\frac{g_n(\o)}{\jb{n}}e^{2\pi inx}, 
\end{equation}

\noi
where $\jb{n} = (1 + 4\pi^2 |n|^2)^\frac{1}{2}$.
 Namely,  
expectations are taken with respect 
to the law $\mu$ of 
 the Ornstein-Uhlenbeck loop  in~\eqref{Q0}.
In this non-homogeneous setting, the problem is not scaling invariant 
and some extra care is needed.
 See for example the proofs of 
Lemma \ref{LEM:u-Q} and 
Proposition~\ref{PROP:FOC}.
In Remark \ref{REM:mean}, 
we indicate the necessary modifications
for handling the case
of  the mean-zero  Brownian loop in 
\eqref{bloop}.

\subsection{Rescaled and translated ground state}

In the following,  we compare a function~$u$ on the circle $\T$ 
to translations and rescalings of the ground state $Q$.
For this purpose, we introduce the $L^2$-invariant scaling operator
 $D_\ld$ on $\R$ by
\begin{align}
D_\ld f(x)=\ld^{-\frac 12}f(\ld^{-1}x).
\label{G2}
\end{align}

\noi
Then, given $\dl > 0$, we set
\begin{align*}
Q_{\dl}(x)= D_\dl Q (x) = \dl^{-\frac 12}Q(\dl^{-1} x) .
\end{align*}

\noi
While the ground state $Q$ is defined on $\R$, 
we now interpret it as a function of $x\in [-\frac 12,\frac 12)$
which  we identify with  the torus $\T\cong [-\frac 12,\frac 12)$. 
Then, 
it follows from \eqref{elliptic} that
\begin{align}\dx^2 Q_\dl+Q^5_\dl - 2\dl^{-2}Q_\dl = 0
\label{elliptic2}
\end{align}

\noi
for $x\in (-\frac 12,\frac 12)$.
As a consequence, we have
\begin{equation}\label{G3}
\jb{\dx^2 Q_\dl, v}_{L^2(\T)} 
= -\jb{Q^5_\dl, v}_{L^2(\T)}+2\dl^{-2}\jb{Q_\dl,v}_{L^2(\T)}.
\end{equation}

Given $x_0 \in \T$, we also introduce 
the translation operator $\tau_{x_0}$ defined by 
\begin{align}
\tau_{x_0}f(x) = f(x - x_0),
\label{trans1}
\end{align}

\noi
where $x- x_0$ is interpreted mod 1, 
taking values in $[-\frac 12, \frac 12)$.
The rescaled and translated  version of the ground state:
\begin{align}
Q_{\dl, x_0 }(x)
: =
\tau_{x_0} D_\dl  Q(x) = \dl^{-\frac 12}Q(\dl^{-1}(x- x_0))
\label{G3-1}
\end{align}

\noi
plays an important role in our analysis.
 By this definition, we have  $Q_\dl=Q_{\dl,0}$.
In the following, we use $\dd_\dl $
and $\dd_{x_0}$ to denote differentiations with respect
to the scaling and translation parameters, respectively.
When there is no confusion, we also denote by $D_\ld$
and $\tau_{x_0}$ the scaling and translation operators
for  functions on the real line.
 
When restricted to the torus, the rescaled and translated version of the ground state belongs to $H^1(\T)$, but not to $H^2(\T)$ (nor to $H^k(\T)$ for any higher $k$). 
This is due to the fact that 
$Q'_\dl$ on $[-\frac 12, \frac 12)$, truncated at $x = \pm \frac 12$, does not respect the periodic boundary condition.
%
%
In order to overcome this problem, we introduce an even cutoff function 
$\rho \in C^\infty(\R; [0, 1])$ such that
\begin{equation} \label{rhodef}
 \supp(\rho) \subset \big[-\tfrac14, \tfrac 14\big]
 \qquad \text{and}\qquad  \rho \equiv 1 \text{ on } \big[-\tfrac18, \tfrac 18\big],
\end{equation}

\noi
and consider $\rho Q_\dl$, which is smooth at $x = \pm \frac 12$.
In the following, we will view 
 \[ Q^\rho_{\dl} := \rho Q_\dl  = \rho D_\dl Q\]

\noi
as a function on the torus.
By setting
\begin{align}
 Q_{\dl,x_0}^\rho := \tau_{x_0} (\rho Q_\dl)
= \rho(x - x_0) \dl^{-\frac 12}Q(\dl^{-1}(x- x_0)), 
\label{QQ1}
\end{align}

\noi
we  will also view $ Q^\rho_{\dl,x_0}$ as a function on the torus.   
 
Note that since $Q$ decays exponentially as $|x| \to \infty$ on $\R$, we have
\begin{equation}\label{corr1}
\begin{split}
\| Q_{\dl, x_0}^\rho\|_{L^2(\T)}^2 &=\dl^{-1}\int_{\R} \rho^2(x) Q^2(\dl^{-1}x)\, \ds  x  \\
&=\int_{\R}\rho^2(\dl x)  Q^2(x)\, \ds  x\\
&= \|Q\|_{L^2(\R)}^2+O(\exp(-c\dl^{-1}))
\end{split}
\end{equation}

\noi
as $\dl \to 0$.

\begin{remark}\label{REM:ortho}\rm

From \eqref{G3-1}, we have
\begin{align}
\dd_\dl Q_{\dl}(x) = -\tfrac 12 \dl^{-\frac 32} Q(\dl^{-1}x)
- \dl^{-\frac 32} (\dl^{-1} x) Q'(\dl^{-1}x), 
\label{Q1}
\end{align}

\noi
which is an even function. Similarly, $ \dd_\dl Q_\dl^\rho = \dd_\dl (\rho Q_\dl)$ is an even function. 
On the other hand, we have
\begin{align}
\dd_{x_0} Q_{\dl}(x) 
= \dd_{x_0} Q_{\dl, x_0}(x)|_{x_0 = 0}
=- \dl^{-\frac 32} Q'(\dl^{-1}x),
\label{Q2}
\end{align}

\noi
which is an odd function.
Similarly, 
\begin{align}
\dd_{x_0} Q_\dl^\rho = \dd_{x_0}Q_{\dl, x_0}^\rho|_{x_0 = 0}
= \dd_{x_0} (\tau_{x_0} (\rho Q_\dl))|_{x_0 = 0} 
\label{QQ0}
\end{align}
 is also 
an odd function.
By writing
\[ 
\jb{\dd_\dl Q_{\dl}^\rho, 
\dd_{x_0}  Q_{\dl}^\rho}_{H^{k}(\T)}
= 
\jb{ \dd_\dl Q_{\dl}^\rho, 
(1 - \dx^2)^{k} \dd_{x_0}  Q_{\dl}^\rho}_{L^2(\T)}, \]

\noi
we see that 
$\dd_\dl Q_{\dl}^\rho$ and 
$\dd_{x_0}  Q_{\dl}^\rho$ are orthogonal in $H^{k}(\T)$, 
$k \in \Z_{\ge 0}$, 
since $(1 - \dx^2) \dd_{x_0}  Q_{\dl}^\rho$ is an odd function.

Similarly, 
for given  $x_0 \in \T$, 
we have 
\begin{align}
\begin{split}
\dd_\dl Q_{\dl, x_0}(x) & = -\tfrac 12 \dl^{-\frac 32} Q(\dl^{-1}(x - x_0))
- \dl^{-\frac 32} (\dl^{-1} (x - x_0) )Q'(\dl^{-1}(x - x_0)), \\
\dd_{x_0} Q_{\dl, x_0}(x) & = - \dl^{-\frac 32} Q'(\dl^{-1}(x - x_0)).
\end{split}
\label{Q3}
\end{align}

\noi
By  parity considerations of these functions centered at $x = x_0$, 
we also conclude that $\dd_\dl Q_{\dl, x_0}^\rho$ and 
$\dd_{x_0}  Q_{\dl, x_0}^\rho$ are
orthogonal on $H^k(\T)$, 
$k \in \Z_{\ge 0}$.

Given $ \dl > 0$ and $x_0 \in \R$, 
let 
$Q_{\dl, x_0 }
= \tau_{x_0} D_\dl  Q$ as a function on the real line.
Then, 
a similar consideration shows orthogonality of 
$\dd_\dl Q_{\dl, x_0}$
and 
$\dd_{x_0} Q_{\dl, x_0}$
in $\dot H^k(\R)$, $k \in \Z_{\geq 0}$. 

\end{remark}

\begin{remark} \label{REM:LeoM}
\rm The set 
$$\M = \big\{ e^{i \ta}  Q^\rho_{\dl,x_0} 
= e^{i \ta}  \tau_{x_0}(\rho Q_{\dl} )
: \dl \in \R_+, \,  x_0 \in \T, 
\text{ and } \ta \in \R \big\} \subset H^1(\T)$$ 
is a smooth manifold of dimension 3, embedded in $H^1(\T)$. The presence of the cutoff function~$\rho$ is fundamental for this to be true. 
Indeed, if we instead consider the set 
$$\M' = \{ e^{i \ta} Q_{\dl,x_0} 
: \dl \in \R_+, \,  x_0 \in \T, 
\text{ and } \ta \in \R \} \subset H^1(\T),$$

\noi
then it turns out that $\M'$
is  \emph{not}  a Lipschitz submanifold of $H^1(\T)$ (at least with the parametrization induced by 
$(\dl, x_0, \ta)$). Indeed, if $\dl = 1$, $0 < x_0 \ll 1$, and $\ta = 0$, recalling that $Q'$ is odd, we have  
\begin{align*}
\| Q_{1,0} - Q_{1,x_0} \|_{{\dot H}^1(\T)}^2 
&= \int_{-\frac 12 + x_0}^{\frac 12} |Q'(x) - Q'(x-x_0)|^2 dx \\
& \hphantom{XX}+ \int_{-\frac 12}^{-\frac12 + x_0} |Q'(x) - Q'(x-x_0 + 1)|^2 dx \\
&\ges \int_{-\frac 12}^{-\frac12 + x_0}  \big|Q'\big(-\tfrac 12\big)\big|^2 dx \\
&\gtrsim |x_0|,
\end{align*}
and thus $\| Q_{1,0} - Q_{1,x_0} \|_{\dot H^1(\T)} \gtrsim |x_0|^\frac 12$,
which shows that the dependence of $Q_{\dl, x_0}$ in $x_0$ is not Lipschitz.
 Similar considerations hold for $e^{i\ta} Q_{\dl,x_0}$ for every value of $\ta, \dl, x_0$. In view of Lemmas \ref{LEM:surface} and \ref{LEM:meas} below, it is crucial that the set $\M$ is a $C^2$-manifold. 
\end{remark}

\subsection{Stability of the optimizers of  the GNS inequality}
\label{SUBSEC:stab}
We begin by establishing stability of the 
optimizers  of the Gagliardo-Nirenberg-Sobolev inequality \eqref{GNSdisp};
 if $u$ is ``far'' from all rescalings and translations of 
 the ground state $Q$  in the $L^2$-sense, 
 then the GNS inequality \eqref{GNSdisp} is ``far'' from being sharp in the sense of \eqref{G1} below.

Given $\g > 0$, 
define the set $S_\g$ 
by setting 
\begin{equation}
\begin{split}
S_\g = \big\{ & u  
\in  L^2(\T): \,   \|u\|_{L^2(\T)}\le \|Q\|_{L^2(\R)}, \\
&  \|P_{\leq k} \pi_{\ne 0 }u\|_{L^6(\T)}^6
  \le (C_{\GNS}-\g)\|P_{\le k} u \|_{\dot{H}^1(\T)}^2\|Q\|^4_{L^2(\R)}
\text{  for all $k\ge 1$}\big\}.
\end{split}
\label{G1}
\end{equation}

\noi
Here,  
$P_{\leq k} u$ denotes the Dirichlet projector onto the frequencies $\{|n| \leq 2^k \}$ defined 
in \eqref{LP2}
and 
$\pi_{\ne 0 }$ denotes the projection onto the mean-zero part defined in \eqref{EN2}.

Fix $\g > 0$.
Young's inequality and \eqref{G1} yield
\begin{align}
\begin{split}
  \|P_{\leq k} u\|_{L^6(\T)}^6
& \le \frac{C_{\GNS}-\frac \g2}{C_{\GNS}-\g}
\|P_{\leq k} \pi_{\ne 0} u\|_{L^6(\T)}^6
+ C_\g |\pi_0 u |^6\\
& 
  \le \big(C_{\GNS}-\tfrac \g2\big)\|P_{\le k} u \|_{\dot{H}^1(\T)}^2\|Q\|^4_{L^2(\R)}
+   C_\g \|Q\|^6_{L^2(\R)}
\end{split}
\label{G1x}
\end{align}

\noi
for any $u \in S_\g$, 
where $\pi_0 = \Id - \pi_{\ne 0}$.
Then, 
   by repeating the argument in Subsection~\ref{SUBSEC:1d}
   with \eqref{G1x}, 
we can show that
\begin{equation}
  \E\Big[e^{\frac{1}{6}\int_\T|u(x)|^6\,\ds  x}, 
  S_\g\Big]\le C(\g)<\infty. 
  \label{E-upper}
\end{equation}

The main goal of this subsection is to prove 
 the following ``stability'' result.
 In view of~\eqref{E-upper}, 
this lemma allows us to restrict our attention to 
 a small neighborhood of the orbit of the ground state $Q$
 in the subsequent subsections.

\begin{lemma}\label{LEM:u-Q}
Let $\rho$ be as in \eqref{rhodef}.
Given any $\eps > 0$ and $\dl^* > 0$, 
there exists $\g(\eps, \dl^*) > 0$ such that the following holds;
suppose that 
a  
function  $u \in L^2(\T)$ with 
 $\|u\|_{L^2(\T)} \le \|Q\|_{L^2(\R)}$
 satisfies
  \begin{align}
 \|u- e^{i \ta} Q_{\dl, x_0}^\rho\|_{L^2(\T)}\ge \eps
 \label{GX1}
  \end{align}
 
 \noi
 for all    $0<\dl < \dl^*$, $x_0\in  \T$, and $\ta \in \R$, 
 where  $Q_{\dl,x_0}^\rho
 = \tau_{x_0} (\rho Q_\dl)$
is   as in \eqref{QQ1}.
   Then,  we have  $u\in S_{\g(\eps,\dl^*)}$.
\end{lemma}

Recall that 
there exists no optimizer for 
the GNS inequality   \eqref{JJ1aa} on the torus $\T$; see Lemmas \ref{LEM:periodicGN} and \ref{LEM:opt}.
Lemma \ref{LEM:u-Q} states that 
``almost optimizers'' of 
the GNS inequality~\eqref{JJ1aa} on  $\T$
exist only in a small $L^2(\T)$-neighborhood of
$e^{i \ta} Q_{\dl, x_0}^\rho$
for $0 < \dl < \dl^*$, where $\dl^* > 0$ is any  given  small number.
This is not surprising since
(i)~the ground state $Q$ (up to symmetries)  is the unique optimizer
for the GNS inequality~\eqref{GNSdisp} on the real line $\R$ (see Remark \ref{REM:opt})
and 
(ii)~we see from the definition 
\eqref{G3-1}  that, 
as $\dl \to 0$, 
$Q_{\dl}^\rho = \rho Q_\dl $ on $ \T \cong [-\tfrac 12, \tfrac 12) $ 
becomes a 
more and more accurate approximation of (a dilated copy of) the ground state $Q$ on $\R$
(and thus becomes a more and more accurate
 almost optimizer for 
the GNS inequality   \eqref{JJ1aa} on the torus $\T\cong [-\tfrac 12, \tfrac 12) $).

\begin{proof}
We first make preliminary computations
which allow us to reduce the problem to the mean-zero case.
Suppose that 
$u \in H^1(\T)$ satisfies
 \begin{align}
 \|u\|_{L^2(\T)} \le \|Q\|_{L^2(\R)}
 \label{GX1a}
 \end{align}
but
$u \notin S_\g$ for some $\g > 0$.
Then, from the GNS inequality \eqref{JJ1aa} on~$\T$
for mean-zero functions (see Lemma \ref{LEM:periodicGN})
and the definition \eqref{G1} of $S_\g$, there exists $k \in \NB$ such that 
\begin{align}
\begin{split}
 (C_{\GNS}-\g)\|P_{\le k} u \|_{\dot{H}^1(\T)}^2\|Q\|^4_{L^2(\R)}
& <  \|P_{\leq k} \pi_{\ne 0 }u\|_{L^6(\T)}^6\\
&  \le C_{\GNS}\|P_{\le k} u \|_{\dot{H}^1(\T)}^2\| \pi_{\ne 0 }u\|_{L^2(\T)}^4.
\end{split}
\label{GX1b}
\end{align}

\noi
Thus, from \eqref{GX1a} and \eqref{GX1b}, we obtain
\begin{align*}
\begin{split}
| \pi_0 u |^2 & = \| \pi_0 u \|_{L^2(\T)}^2 
=  \|  u \|_{L^2(\T)}^2  - \| \pi_{\ne 0} u \|_{L^2(\T)}^2 \\
& < \bigg(1 - \sqrt{\frac {C_{\GNS}-\g}{C_{\GNS}}}\bigg)\|Q\|_{L^2(\R)}^2
= O(\g) \|Q\|_{L^2(\R)}^2
\end{split}
\end{align*}

\noi
as $\g \to 0$.
Hence, if we have \eqref{GX1} for some $\eps > 0$, 
then there exists $\g_0 = \g_0(\eps) > 0$ such that 
 \begin{align}
 \|\pi_{\ne 0} u- e^{i \ta} Q_{\dl, x_0}^\rho\|_{L^2(\T)}\ge \frac{\eps}{2}.
 \label{GX3}
  \end{align}

\noi
for any   $u\not \in S_{\g}$
with  $0 < \g < \g_0$.

We prove the lemma by contradiction.
Suppose that there is no such  $\g(\eps,\dl^*)$. 
 Then,  there exist
$\eps > 0$, $\dl^* > 0$, 
$\{u_n\}_{n \in \NB} \subset  L^2 (\T)$ with 
\begin{align}
\|u_n \|_{L^2(\T)} \leq \|Q\|_{L^2(\R)},
\label{G3a}
\end{align}

\noi 
and
$\g_n\to0$  
such that 
  \begin{align}
  \|u_n-e^{i\ta} Q_{\dl,x_0}^\rho\|_{L^2(\T)}\ge \eps
\label{G3b}
  \end{align} 
 \noi
 for any  $0<\dl < \dl^*$, $x_0 \in  \T$, and $\ta \in \R$
 but 
 $u_n\notin S_{\g_n}$. 
By the discussion above, in particular from \eqref{GX3}, 
there exists $N_0(\eps) \in \NB$ such that 
  \begin{align*}
  \|\pi_{\ne 0 } u_n- e^{i\ta} Q_{\dl,x_0}^\rho\|_{L^2(\T)}\ge \frac \eps2
  \end{align*} 
 \noi
 for any  $0<\dl < \dl^*$,  $x_0 \in  \T$, $\ta \in \R$, 
 and  $n \geq N_0$
 but  $\pi_{\ne 0} u_n\notin S_{\g_n}$. 
Therefore, without loss of generality, 
we may assume that $u_n$, satisfying \eqref{G3a}
and \eqref{G3b}, has mean zero (i.e.~$u_n = \pi_{\ne 0} u_n$)
and derive a contradiction.

 By definition \eqref{G1} of $S_\g$, for each $n \in \NB$, 
 there exists $k_n \in \NB$ such that 
\begin{align}
\|P_{\le k_n}u_n\|^6_{L^6(\T)}
&>(C_{\GNS}-\g_n)\|P_{\le k_n} u_n\|^2_{\dot{H}^1(\R)}\|Q\|_{L^2(\R)}^4,
\label{G4}
\end{align}

\noi
since $u_n\notin S_{\g_n}$. 
Then, from \eqref{GNSdisp} and \eqref{G4} with \eqref{G3a} and $\g_n \to 0$, 
we see that 
 \begin{align}
 \|P_{\leq k_n} u_n \|_{L^2(\T)} \too \|Q\|_{L^2(\R)}
 \label{G4b}
 \end{align}

\noi
 as $n \to \infty$.
In view of the upper bound \eqref{G3a}, 
 we also  have  $\|u_n\|_{L^2(\T)} \to \|Q\|_{L^2(\R)}$.
 Hence, by the Pythagorean theorem, 
 we obtain
\begin{align}
\|u_n-P_{\le k_n} u_n \|_{L^2(\T)} \too0.
\label{G4-1}
\end{align}

\noi
Furthermore, we claim that 
\begin{align}
\|P_{\leq k_n} u_n \|_{\dot H^1(\T)} \too \infty.
\label{G4b1}
\end{align}

\noi
Otherwise, we would have 
$\|P_{\leq k_n} u_n \|_{\dot H^1(\T)} \leq C <  \infty$ for all $n \in \NB$
and thus there exists a subsequence, still denoted by $\{P_{\leq k_n} u_n\}_{n \in \NB}$, 
converging weakly to some $u$ in $\dot H^1(\T)$.
Then, from the compact embedding\footnote{Recall that we work with mean-zero functions on $\T$.} 
of $\dot H^1(\T)$ into $L^2(\T)\cap L^6(\T)$, 
we see that $P_{\leq k_n} u_n$ converges strongly to $u$ in $L^2(\T)\cap L^6(\T)$.
Hence, from 
\eqref{G3a} and then \eqref{G4}, we obtain
\begin{align*}
C_{\GNS}\|u \|^2_{\dot{H}^1(\T)}\|u \|_{L^2(\T)}^4
& \leq \liminf_{n \to \infty} (C_{\GNS} - \g_n)  \|P_{\leq k_n} u_n \|^2_{\dot{H}^1(\T)}\|Q \|_{L^2(\R)}^4 \notag \\
& \leq \liminf_{n \to \infty} \|P_{\le k_n}u_n\|^6_{L^6(\T)}
=  \|u \|^6_{L^6(\T)}.
\end{align*}

\noi
This would imply that $u$ is a mean-zero optimizer
of the Gagliardo-Nirenberg-Sobolev  inequality \eqref {GNSdisp}
on the torus $\T$, which is a contradiction to Lemma \ref{LEM:opt}. 
Therefore,~\eqref{G4b1} must hold.

By continuity, there exists a point $x_n \in \T$ 
such that 
\begin{equation}\label{eqn: zero}
|P_{\le k_n}u_n (x_n)| \leq \| P_{\le k_n}u_n \|_{L^2(\T)}.
\end{equation}

\noi
With $\be_n = \|P_{\leq k_n} u_n \|_{\dot H^1(\T)}^{-1}\to 0$,  
define $v_n:\R\to \mathbb{C}$ by
\begin{align}
v_n(x)=\begin{cases}
P_{\le k_n}u_n (x+x_n-\frac12), &\quad  \text{for } x\in [-\frac12,\frac12], \\
0, & \quad \text{for } |x|>\frac12 + 
\be_n, 
\end{cases}
\label{G4-2}
\end{align}

\noi
and by  linear interpolation
for  $\frac{1}{2} < |x|
\le \frac{1}{2} + 
\be_n$, 
where the addition here is understood mod 1.
Then, 
from   \eqref{eqn: zero}, we have $|v_n(\pm \frac12)| \leq \| P_{\le k_n}u_n \|_{L^2(\T)}$.
Moreover, from  \eqref{G4-2} with \eqref{G3a} and  \eqref{G4b1}, 
we have
  $v_n\in H^1(\R)$ with 
\begin{align}
\|v_n\|_{L^2(\R)}\leq \sqrt{1 + 2\be_n}\, \|P_{\le k_n} u_n \|_{L^2(\T)}
\leq \sqrt{1 + 2 \be_n} \, \|Q\|_{L^2(\R)}
\label{G4c}
\end{align}

\noi
for any $n \in \NB$, 
and 
\begin{align}
\|v_n\|_{\dot H^1(\R)}^2 
= 
\|P_{\le k_n} u_n \|_{\dot H^1(\T)}^2
+ 2\be_n^{-1} |P_{\le k_n}u_n (x_n)|^2 \too \infty,
\label{G4c1}
\end{align}

\noi
as $n \to \infty$.
Hence, from \eqref{G4}, \eqref{G4c}, and \eqref{G4c1}
with \eqref{eqn: zero} and \eqref{G3a}, we have
\begin{align}
 \|v_n\|_{L^6(\R)}^6>(C_{\GNS}-\g_n)
 \al_n  
 \|v_n\|_{\dot{H}^1(\T)}^2\|v_n\|_{L^2(\R)}^4, 
\label{G5x}
 \end{align}

\noi
where 
\begin{align*}
\al_n  = 
 (1+2\be_n)^{-2}
(1 - 2\be_n \|Q\|_{L^2(\R)}^2).
\end{align*}

\noi
Since $\be_n = \|P_{\leq k_n} u_n \|_{\dot H^1(\T)}^{-1}\to 0$,  
we have $\al_n \to 1$ as $n \to \infty$.
Hence, we conclude from~\eqref{G5x}
that there exists $\wt \g_n \to 0 $ such that 
\begin{align}
 \|v_n\|_{L^6(\R)}^6>(C_{\GNS}-\wt \g_n)
 \|v_n\|_{\dot{H}^1(\T)}^2\|v_n\|_{L^2(\R)}^4.
\label{G5}
 \end{align}

With the scaling operator $D_\ld$ as in \eqref{G2}, 
let 
\begin{align}w_n=D_{\|v_n\|_{\dot{H}^1(\R)}}v_n. 
\label{G5-1}
\end{align}

\noi
Then, from \eqref{G4b}, \eqref{G4c}, and \eqref{G5}, 
we have 
\begin{align}
& \|w_n\|_{\dot{H}^1(\R)}=1, \label{G5a}\\
& \|w_n\|_{L^2(\R)}=\|v_n\|_{L^2(\R)} \leq 
\sqrt{1 + 2 \be_n} \,
\|{Q}\|_{L^2(\mathbb R)}
 \text{ and $\|w_n\|_{L^2(\R)} \to  \|{Q}\|_{L^2(\mathbb R)}$}, 
\label{G5b}\\
& \|w_n\|^6_{L^6(\R)}>(C_{\GNS}-\wt\g_n)\|w_n\|_{L^2(\R)}^4.
\label{G5c}
\end{align}

\noi
Since $\{w_n\}_{n \in \NB}$
is a bounded sequence in $H^1(\R)$, 
we can invoke
 the profile decomposition \cite[Proposition 3.1]{HK}
for the (subcritical) Sobolev embedding: $H^1(\R)\hookrightarrow L^6(\R)$.
See also Theorem 4.6 in \cite{KV}.
There exist $J^* \in \Z_{\geq 0} \cup\{\infty\}$, 
a sequence $\{\phi^j \}_{j = 1}^{J^*}$ of non-trivial $H^1(\R)$-functions, 
and a sequence $\{x^j_n\}_{j  = 1}^{J^*}$ for each $n \in \NB$
such that 
up to a subsequence, still denoted by $\{w_n \}_{n \in \NB}$, we have
 \begin{align}
 w_n(x)=\sum_{j=1}^J \phi^j(x-x_n^j)+r^J_n(x)
 \label{G5d}
 \end{align}
 
\noi
for each finite $0 \leq J \leq J^*$, 
where the remainder term $r^J_n$
satisfies
\begin{align}
 \lim_{J \to \infty} \limsup_{n \to \infty} \| r^J_n\|_{L^6(\R)} = 0.
 \label{G5e}
\end{align}

\noi
Here,  $\lim_{J \to \infty} f(J) := f(J^*)$ if $J^* < \infty$.
Moreover, for any finite $0 \leq J \leq J^*$, we have
\begin{align}
\|w_n\|_{L^2(\R)}^2&  =  \sum_{j = 1}^J \|\phi^j \|_{L^2(\R)}^2 + \|r_n^J\|_{L^2(\R)}^2
+ o(1),\label{G6a}\\
\|w_n\|_{\dot H^1(\R)}^2 & =  \sum_{j = 1}^J \|\phi^j \|_{\dot H^1(\R)}^2 + \|r_n^J\|_{\dot H^1(\R)}^2
+ o(1),
\label{G6aa}
\end{align}

\noi
as $n \to \infty$, and 
\begin{align}
\limsup_{n \to \infty}
\|w_n\|_{L^6(\R)}^6  = &  \sum_{j = 1}^{J^*} \|\phi^j \|_{L^6(\R)}^6 .
\label{G6b}
\end{align}

\noi
From \eqref{G5a} and \eqref{G6aa} and taking $n \to \infty$, 
we obtain
\begin{align}
\sup_{j = 1, \dots, J^*} \|\phi^j \|_{\dot H^1(\R)}^2 
\leq \sum_{j = 1}^{J^*} \|\phi^j \|_{\dot H^1(\R)}^2 \leq 1,
\label{G7}
\end{align}

\noi
where an equality holds at the first inequality  if and only if $J^* = 0 $ or $1$.

From \eqref{GNSdisp} and \eqref{G5c}
with $\wt \g_n \to 0$, we have
\begin{align*}
C_{\GNS}\lim_{n\to \infty} \|w_n\|^4_{L^2(\R)}
&=\limsup_{n\to \infty} \|w_n\|^6_{L^6(\R)}. 
\end{align*}

\noi
By \eqref{G6b} followed by \eqref{GNSdisp}, \eqref{G7}, 
and \eqref{G6a} with $\l^2\subset \l^4$,
\begin{align}
C_{\GNS}\lim_{n\to \infty} \|w_n\|^4_{L^2(\R)} 
&=  \sum_{j=1}^{J^*} \|\phi^j\|^6_{L^6(\R)} \notag \\
&\le C_{\GNS}\sum_{j=1}^{J^*} \|\phi^j\|^2_{\dot{H}^1(\R)}\|\phi^j\|_{L^2(\R)}^4 
\label{G8}\\
&\le C_{\GNS}\sum_{j=1}^{J^*} \|\phi^j\|^4_{L^2(\R)}\nonumber\\
&\le C_{\GNS}\lim_{n\to \infty} \|w_n\|^4_{L^2(\R)}.\nonumber
\end{align}

\noi
Here, an equality holds if and only if $J^* = 0$ or $1$.
If $J^* = 0$, 
then it follows from \eqref{G6b}
that $w_n$ tends to 0 in $L^6(\R)$
as $n \to \infty$.
Then, from  \eqref{G5c}, 
we see that $w_n$ tends to 0 in $L^2(\R)$.
This is a contradiction to \eqref{G5b}.
Hence, we must have $J^* = 1$.
In this case,~\eqref{G8} (with $J^* = 1$)
holds with equalities and thus 
we see that $\phi^1$ is 
an optimizer for the Gagliardo-Nirenberg-Sobolev inequality on $\R$.

Hence, 
we conclude from Remark \ref{REM:opt} that 
 there exist  $\s\neq 0$, $\ld>0$,  and $x_0\in \R$ such that 
\[\phi^1=\s  \tau_{x_0} D_\ld Q.\]

\noi
From the last step in \eqref{G8} (with $J^* = 1$ and an equality) with \eqref{G5b}, 
we have 
\[\|\phi^1 \|_{L^2(\R)} = \lim_{n\to \infty} \|w_n\|_{L^2(\R)} = \|Q\|_{L^2(\R)},\]

\noi 
which implies  that $\s=e^{i\ta}$ for some $\ta \in \R$. 
From~\eqref{G5d}
and~\eqref{G5e}
with $J = J^* = 1$,  
we see that $ \tau_{- x_n^1} w_n$ converges weakly 
to $ \phi^1$ in $L^2(\R)$
(which follows from  the (weak) convergence 
of $ \tau_{- x_n^1} w_n$ 
to $ \phi^1$ 
in $L^6(\R)$), 
while~\eqref{G5b} implies convergence of the $L^2$-norms.
Hence,  we obtain strong convergence in $L^2(\R)$:
\[\| w_n- e^{i\ta}\tau_{  x_0  +  x_n^1} D_\ld Q\|_{L^2(\R)}\too 0.\]

\noi
Hence,  from \eqref{G5-1}, we have
\begin{align}
\| v_n-e^{i\ta} \tau_{  y_n }D_{\ld_n}Q\|_{L^2(\R)}\too 0,
\label{G9}
\end{align}

\noi
where 
$y_n =  \|v_n \|_{\dot H^1(\R)}^{-1}(x_0 + x_n^1)$ and 
$\ld_n =  \|v_n \|_{\dot H^1(\R)}^{-1}\ld \to 0$ in view of \eqref{G4c1}.

Recall that  $v_n$ is supported on $[-\frac12 - \be_n, \frac 12+ \be_n]$
with the $L^2$-norm on 
$[-\frac12, \frac 12]^c$
bounded by $\sqrt{ 2 \be_n} \, \|Q\|_{L^2(\R)}$
(which tends to 0 as $n \to \infty$).
Thus, 
it follows from \eqref{G9} that 
\begin{align}
\begin{split}
J_n :  \! & = \int_{|y| \geq \frac{1}{2\ld_n}}
|Q( y - \ld_n^{-1} y_n) |^2 dy \\
& = 
 \| \tau_{y_n} D_{\ld_n}Q\|_{L^2([-\frac 12, \frac 12]^c)}
 \too 0, 
 \end{split}
 \label{G10}
\end{align}

\noi
as $n \to \infty$.
In the following, 
we denote by $Q_{\ld_n, y_n}$ 
the dilated and translated ground state $Q$,  viewed as 
 a periodic function on $\T$, 
and by $ \tau_{y_n} D_{\ld_n} Q$ 
the dilated and translated ground state viewed as 
 a function on the real line.
By possibly choosing a subsequence, we assume that 
$y_n \ge 0$ without loss of generality.
Write  
\[\T \cong [-\tfrac 12, \tfrac 12) 
= [-\tfrac 12+ y_n , \tfrac 12) \cup [-\tfrac 12, - \tfrac 12 + y_n)  
=: I_{1, n} \cup I_{2, n}. \]

\noi
Note that while  $Q_{\ld_n, y_n}$ 
and $ \tau_{y_n} D_{\ld_n} Q$ 
coincide on $I_{1, n}$, they 
do not coincide on $I_{2, n}$.
Thanks to the exponential decay of the ground state $Q$, we have
\begin{align}
\begin{split}
\|  \tau_{y_n} D_{\ld_n}Q \|_{L^2(I_{2, n})} 
& \leq 
 \|D_{\ld_n}Q\|_{L^2((-\infty, - \frac 12])}
 = O(\exp(-c \ld_n^{-1}))
\too 0.
\end{split}
\label{G12}
\end{align}

\noi
On the other hand, on $I_{2, n}$, 
we have
$Q_{\ld_n, y_n}(x) = \ld_n^{-\frac 12}
Q(\ld_n^{-\frac 12}(x+1 - y_n))$.
Thus, from a change of variables and \eqref{G10}, we have
\begin{align}
\begin{split}
 \| Q_{\ld_n, y_n} \|_{L^2(I_{2, n})}
& \leq \|\tau_{y_n} D_{\ld_n}Q\|_{L^2([\frac 12, \infty))} 
\too 0.
\end{split}
\label{G12a}
\end{align}

\noi
Therefore, from \eqref{G9},  \eqref{G12} and \eqref{G12a}, we obtain
\begin{align}
\begin{split}
\| v_n & - e^{i \ta}  Q_{\ld_n, y_n} \|_{L^2(\T)}\\
& \leq
\| v_n -e^{i \ta}  \tau_{y_n} D_{\ld_n} Q \|_{L^2([-\frac 12, \frac 12])}
+
 \| Q_{\ld_n, y_n} \|_{L^2(I_{2, n})}
+ \|  \tau_{y_n} 
D_{\ld_n}
Q \|_{L^2(I_{2, n})} \\
& \too 0, 
\end{split}
\label{G13}
\end{align}

\noi
as $n \to \infty$. Moreover, since $\ld_n \to 0$, 
it follows from \eqref{rhodef} that 
\begin{align} \label{G13b}
\begin{split}
\|  e^{i \ta}  Q_{\ld_n, y_n}^\rho - e^{i \ta}  Q_{\ld_n, y_n} \|_{L^2(\T)} 
&= \|  \rho  Q_{\ld_n, 0} -   Q_{\ld_n, 0} \|_{L^2(\T)} \\
&\le \| Q_{\ld_n, 0} \|_{L^2([-\frac 18, \frac18]^c)} \\
&= \| Q \|_{L^2([-\frac {\ld_n}8, \frac {\ld_n}8]^c)} \too 0.
\end{split}
\end{align}
Finally, by combining \eqref{G4-1}, 
\eqref{G4-2}, \eqref{G13}, and \eqref{G13b}, 
we obtain 
a contradiction
to~\eqref{G3b}.
This completes the proof of Lemma \ref{LEM:u-Q}.
\end{proof}

\subsection{Orthogonal coordinate system in a neighborhood of the soliton manifold}

In view of 
 Lemma \ref{LEM:u-Q} and \eqref{E-upper}, 
 in order to  prove  $Z_{6,K=\|Q\|_{L^2(\R)}}<\infty$, it  suffices to show that
\begin{equation}
\E\Big[e^{\frac{1}{6}\int_\T|u(x)|^6\,\ds  x}
\ind_{\{\|u\|_{L^2(\T)} \leq  \|Q\|_{L^2(\R)}\}}
,  U_\eps(\dl^*) \Big]<\infty
\label{eqn: to-bound}
\end{equation}

\noi
for some small $\eps, \dl^*> 0$, where $ U_{\eps}(\dl^*)$ is defined by 
\begin{equation}
\begin{split}
U_{\eps}(\dl^*)=
\big\{u\in L^2(\T): 
\ &  \|u-e^{i \ta} Q_{\dl, x_0}^\rho\|_{L^2(\T)}<\eps \\
&  \text{for some }  0 < \dl < \dl^* , \, x_0\in  \T,
\text{ and } \ta \in \R  \big\}
\end{split}
\label{Ueps}
\end{equation}

\noi
with  $Q_{\dl,x_0}^\rho
 = \tau_{x_0} (\rho Q_\dl)$
   as in \eqref{QQ1}.
Namely, the domain of integration $U_\eps(\dl^*) $ is an $\eps$-neighborhood of the (approximate) soliton manifold
$\M = \M (\dl^*)$:
\[ \M = \big\{ 
e^{i \ta} Q_{\dl, x_0}^\rho = e^{i \ta} \tau_{x_0} (\rho Q_\dl)
: 
0 < \dl < \dl^* , \, x_0\in  \T,
\text{ and } \ta \in \R  \big\}.\]

\noi
Here, we say ``approximate'' due to the presence of the cutoff function $\rho$.
See Remark \ref{REM:LeoM}.
Recall that the expectation in \eqref{eqn: to-bound}
is taken with respect to 
 the Ornstein-Uhlenbeck loop  in~\eqref{Q0} whose law is given by the Gaussian measure 
$\mu$ in~\eqref{Gibbs3a} with the Cameron-Martin space $H^1(\T)$.
Our main goal in this subsection is to endow $U_\eps(\dl^*)$
with an ``orthogonal'' coordinate system, 
where the orthogonality is measured in terms of $H^1(\T)$.
This then allows us to introduce a change of variables
for the integration in 
\eqref{eqn: to-bound}; see Subsection~\ref{SUBSEC:change}.

Given  $0 < \dl  \ll 1$, $x_0 \in \T$, and $0\le \theta <  2\pi$, we define 
$V_{\dl, x_0, \theta}
= V_{\dl, x_0, \theta}(\T)$ by 
\begin{align}
\begin{split}
V_{\dl,x_0,\theta}=
\big\{u\in L^2(\T): \, 
& \jb{u, (1-\dx^2)e^{i\theta}\dd_\dl  Q_{ \dl, x_0}^\rho}_{L^2(\T)}=0, \\
& \jb{u, (1-\dx^2)e^{i\theta}\partial_{x_0} Q_{\dl,x_0}^\rho}_{L^2(\T)}=0,\\
& {\jb{u, (1-\dx^2)ie^{i\theta} Q_{\dl, x_0}^\rho}}_{L^2(\T)} = 0 \big\}.
\label{eqn: Vdelta-def}
\end{split}
\end{align}

\noi
We point out that, due to the insufficient regularity of $u \in L^2(\T)$, 
instead of the $H^1$-inner product, 
we  measure the orthogonality 
in \eqref{eqn: Vdelta-def}
with respect to the $L^2$-inner product with a weight 
$(1- \dx^2)$ so that the inner products 
in \eqref{eqn: Vdelta-def} are well defined for $u \in L^2(\T)$.
Here, we emphasize that the inner product on $L^2(\T)$ defined in \eqref{inner1}
is real-valued.
It is easy to check that  $(\tau_{x_0}\rho)\dd_\dl Q_{ \dl, x_0}$ and $\dd_{x_0} (\tau_{x_0}\rho)Q_{ \dl, x_0}$
are orthogonal in $L^2(\T)$ and $H^1(\T)$
(see Remark~\ref{REM:ortho}), and similarly that they are orthogonal to $i(\tau_{x_0}\rho)Q_{ \dl, x_0}$ in $H^1(\T)$ 
(viewed  as a real vector space with the inner product in \eqref{inner1}). 
Hence,  
the space $V_{\dl, x_0,\theta}$ denotes a \emph{real} subspace of codimension 3
in $L^2(\T)$, 
 orthogonal (with the weight $(1 - \dx^2)$) 
 to the tangent vectors
$e^{i\theta}\dd_\dl Q_{ \dl, x_0}^\rho$,   $e^{i\theta}\dd_{x_0} Q_{ \dl, x_0}^\rho
= e^{i\theta}\dd_{x_0}\big( \tau_{x_0} (\rho Q_\dl)\big)$, 
and $\dd_\ta(e^{i\theta} Q_{\dl, x_0}^\rho ) = ie^{i\theta} Q_{\dl, x_0}^\rho$
of the soliton manifold $\M$
(with $0 < \dl \ll 1$).
The following proposition shows that  a small neighborhood of $\M$
can be endowed with an orthogonal coordinate system in terms of 
 $e^{i\theta}\dd_\dl Q_{ \dl, x_0}^\rho$, 
 $e^{i\theta}\dd_{x_0} Q_{ \dl, x_0}^\rho$, 
$ie^{i\theta} Q_{\dl, x_0}^\rho$,
and $V_{\dl, x_0,\theta}$.

\begin{proposition}\label{PROP:FOC}
Given small $\eps_1 > 0$, 
there exist $\eps = \eps(\eps_1) > 0$ and 
$ \dl^*=  \dl^*(\eps_1)  > 0$
 such that 
\begin{align*}
\begin{split}
 U_\eps(\dl^*) \subset \big\{u\in L^2(\T):\
& \|u-
e^{i\ta} Q_{\dl,x_0}^\rho\|_{L^2(\T)}<\eps_1, 
\ u-e^{i \ta} Q_{\dl, x_0}^\rho\in V_{\dl, x_0,\theta} (\T) \\
&  \text{for some  }  0 < \dl <  \dl^* , \, x_0\in \T, \text{ and }\ta \in \R  
\big\}.
\end{split}
\end{align*}

\end{proposition}

\begin{remark} \label{REM:Leo1}
\rm
The series of works by Nakanishi and Schlag \cite{ns1,ns2,ns3,ns4} and Krieger, Nakanishi, and Schlag \cite{kns1,kns2,kns3,kns4,kns5} use a coordinate system similar to the one given by Proposition \ref{PROP:FOC},\footnote{Without the multiplication by (a translate) of the cutoff function $\rho$.}
  but centered around a \emph{single} soliton. 
See, for example,  Section 2.5 in~\cite{kns2}. Thanks to the symmetries of the problem
on $\R^d$, it is easy to extend the coordinate system to a tubular neighborhood of the soliton manifold in these works. 
In the setting of Proposition~\ref{PROP:FOC} on the torus $\T$, however, we lack dilation symmetry, which makes it impossible to use such a soft argument to conclude
the desired result
(namely, endow $U_\eps(\dl^*)$
with an orthogonal coordinate system).
\end{remark}

\begin{proof}

We first show that the claimed result
holds true in the case of the real line (without the extra factor $\tau_{x_0} \rho$).
Given $\g_0 \in \R$,  $0 < \dl  \ll 1$, $x_0 \in \T$, and $0\le \theta <  2\pi$, 
we define $V^{\g_0}_{\dl, x_0, \theta} = V^{\g_0}_{\dl, x_0, \theta}(\R)$ by 
\begin{equation}
\begin{split}
V_{\dl,x_0,\theta}^{\g_0}
 = \big\{u\in L^2(\R): \, 
& \jb{u, (\g_0^2-\dx^2)e^{i\theta}\dd_\dl  Q_{ \dl, x_0}}_{L^2(\R)}=0, \\
& \jb{u, (\g_0^2-\dx^2)e^{i\theta}\partial_{x_0}  Q_{\dl,x_0}}_{L^2(\R)}=0,\\
& {\jb{u, (\g_0^2-\dx^2)ie^{i\theta} Q_{\delta, x_0}}}_{L^2(\R)} = 0\big\}, 
\label{eqn: Vdelta0-def} 
\end{split}
\end{equation}

\noi
where 
 the inner product on $L^2(\R)$ is real-valued as defined in \eqref{inner1}.

Consider the map $H: \big(\R \times L^2(\R)\big) \times 
\big(\R_+\times \R\times
(\R/2\pi\Z)\times V^0_{1,0,0}\big) \to L^2(\R)$ defined by
\begin{equation}\label{F1}
H(\g_0,u,\delta,x_0,\theta,w)=u - e^{i\theta}Q_{\dl,x_0} - P_{V_{\dl,x_0,\ta}^{\g_0}}w, 
\end{equation}

\noi
where $P_{V^{\g_0}_{\dl,x_0,\ta}}$ is the projection onto $V^{\g_0}_{\dl,x_0,\ta}(\R)$ in $L^2(\R)$.
It is easy to see that $H$ is Fr\'echet differentiable and $H(0,Q_{1,0},1,0,0,0) = 0$.
Moreover, the Fr\'echet derivative of $H$ in the 
$\big(\R_+ \times \R \times (\R/2\pi\Z) \times V^0_{1,0,0}\big)$-variable 
at $(0,Q_{1,0},1,0,0,0)$ is
\[(\ds_2 H)(0,Q_{1,0},1,0,0,0)= 
- \ds  (e^{i\ta}Q_{\dl,x_0})|_{(\dl, x_0, \ta) = (1,  0,  0)} -  \Id_{V^{0}_{1,0,0}},\]

\noi
where
\[\big(\ds  (e^{i\ta} Q_{\dl,x_0})\big)(\al, \be , \g)
= \al e^{i\ta}\partial_{\dl}Q_{\dl,x_0}
+ \be e^{i\ta}\partial_{x_0}Q_{\dl,x_0}
+ \g  ie^{i\ta}Q_{\dl,x_0}
\]

\noi
for 
 $(\al, \be, \g) \in \R \times \R \times \R$. 


The image of $Z = \ds  (e^{i\ta}Q_{\dl,x_0})|_{(\dl, x_0, \ta) = (1, 0, 0)}$ 
has dimension 3, while 
the image of $\Id_{V^{0}_{1,0,0}}$,
namely the subspace $V^{0}_{1,0,0}$,  has codimension 3.
Moreover,  if $u$ lies in the intersection
of the  image of $Z$
and the image of $\Id_{V^{0}_{1,0,0}}$, 
then the definition 
\eqref{eqn: Vdelta0-def} of $V^{0}_{1,0,0}$
and 
the orthogonality of
$e^{i\theta}\dd_\dl Q_{ \dl, x_0}$, 
$e^{i\theta}\dd_{x_0} Q_{ \dl, x_0}$, 
and $ie^{i\theta} Q_{\dl, x_0}$ in $\dot H^1(\R)$
(see Remark \ref{REM:ortho})
allow us to conclude that $u = 0$. 
Hence, $\ds_2 H$ is invertible at $(0,Q_{1,0},1,0,0,0)$. 
By the implicit function theorem, (\cite[Theorem 26.27]{driver}), 
there exists a
neighborhood in $\R \times L^2(\R)$ of $(0, Q_{1,0})$
of  the form:
\begin{align}
W = \big\{|\g| < \g_*\big\} \times \big\{\norm{u-Q_{1,0}}_{L^2 (\R)} < \eps\big\},
\label{J0}
\end{align}

\noi
and a $C^1$-function $b = (\dl, x_0, \ta, w)
: W \to \R_+\times \R\times
(\R/2\pi\Z)\times V^0_{1,0,0}$
%
such that 
\begin{align}
b(0, Q_{1, 0}) = (1, 0, 0, 0) \qquad \text{and}\qquad  
H\big(\g_0, u, b(\g_0, u)\big) = 0.
\label{F2}
\end{align}

\noi
Namely, from \eqref{F1} and \eqref{F2}, we have
\begin{equation}\label{Seps-form}
u = e^{i\ta}Q_{\dl,x_0}+v
\end{equation}

\noi
for some $(\dl, x_0,\ta, v ) \in \R_+\times \R\times (\R/2\pi\Z)
\times  V_{\delta,x_0,\theta}^{\g_0}$.
Moreover, by continuity of $b$ with \eqref{F2}, 
if we choose $\g_*= \g_*(\eps_0)$
and $ \eps=\eps(\eps_0)$ 
sufficiently small, 
then we can guarantee
\begin{align}
\|v\|_{L^2(\R)} < \eps_0.
\label{F3}
\end{align}

Now,  suppose that 
\begin{align}
\|u-e^{i \ta_0}Q_{\dl_0,x_0}\|_{L^2(\R)} <  \eps
\label{J1}
\end{align}

\noi
for some $( \dl_0, x_0, \ta_0) \in \R_+ \times \R\times ( \R/2\pi\Z)$
with $0 < \delta_0 < \g_*$.
Note that we have
\begin{align}
 Q_{\dl_0, x_0} =  \tau_{x_0} D_{\dl_0} Q_{1,0}
 \label{J1a}
\end{align}

\noi
where $D_{\dl_0}$ is the scaling operator on the real line defined in \eqref{G2}
and $\tau_{x_0}$ denotes the translation operator
for functions on the real line.
By setting  
\begin{align}
T = e^{i\ta_0}  \tau_{x_0}  D_{\dl_0},
\label{J1b}
\end{align}

\noi 
we can rewrite~\eqref{J1} as
\begin{align*}
\|T^{-1} u-Q_{1,0}\|_{L^2(\R)} <  \eps
\end{align*}

\noi
since $T$ is an isometry on $L^2(\R)$.
Then, it follows from \eqref{J0} that 
 $(\delta_0, T^{-1} u) \in W$.
 Hence, from 
\eqref{Seps-form}, we have
\begin{equation}
T^{-1} u=e^{i\ta}Q_{\dl,x}+v
\label{J3}
\end{equation}

\noi
for some $(\dl, x,\ta,v)
\in \R_+\times \R\times (\R/2\pi\Z) \times V^{\dl_0}_{\dl,x,\ta}$ near $(1,0,0,0)$. 
Therefore, we obtain 
\begin{equation}
u=e^{i \ta_1} Q_{\dl_1,x_1}+ Tv
\label{J4}
\end{equation}

\noi
for some $(\dl_1, x_1,\ta_1 ) \in \R_+ \times \R\times (\R/2\pi\Z)$
such that 
\begin{align}
e^{i\ta_1}Q_{\dl_1,x_1} = Te^{i\ta}Q_{\dl, x}.
\label{J5}
\end{align}

\noi
From \eqref{J1b} and \eqref{J5}
with \eqref{F3}, 
we can easily check that  $Tv \in V^1_{\dl_1, x_1,\ta_1}$
and  $\|T v\|_{L^2(\R)} < \eps_0$
for $v \in V^{\delta_0}_{\dl,x,\ta}$
in  \eqref{J3}. 
Hence, we  conclude that $u$ in \eqref{J4} has the desired form
in this real line case.

\medskip

We now prove the  claim in the case of the torus $\T$.
In this case, a scaling argument such as \eqref{J1a} no longer works
and we need to proceed with care.
By a translation and a rotation, we may assume that $u\in L^2(\T)$ satisfies  
\begin{align}
\|u- \rho Q_{\dl_0}\|_{L^2(\T)}<  \eps,
\label{QQ2}
\end{align}

\noi
where $Q_{\dl_0} =  Q_{\dl_0,0}$.  
Extending $u$ by 0 outside $[-\frac 12,\frac12]$, 
we obtain a function in $L^2(\R)$, which we compare to the rescaled  soliton on the real line.
From \eqref{rhodef} and \eqref{QQ2}, we have 
\begin{align*}
\begin{split}
\|u-Q_{\dl_0}\|_{L^2(\R)}^2
&= \int_\R |u(x) - \rho(x) Q_{\dl_0}(x) - (1-\rho(x)) Q_\dl(x)|^2 \, \ds x\\
&= \int_{-\frac 12}^{\frac 12}|u(x)-\rho(x) Q_{\dl_0}(x)|^2\,\ds  x\\
&\hphantom{X}
+2\Re\Big( \int_{\R}(u(x) - \rho(x) Q_{\dl_0}(x))(1-\rho(x)) Q_\dl(x)\,\ds  x\Big)\\
&\hphantom{X}
+ \int_{\R} |(1-\rho(x))Q_\dl(x)|^2 \,\ds x\\
&\le \int_{-\frac 12}^{\frac 12}|u(x)-\rho(x) Q_{\dl_0}(x)|^2\,\ds  x \\
&\hphantom{X}
+ 2\bigg(\int_{-\frac 12}^{\frac 12}|u(x)-\rho(x) Q_{\dl_0}(x)|^2\,\ds  x\bigg)^\frac 12 \bigg(\int_{\R} |(1-\rho(x))Q_\dl(x)|^2 \,\ds x\bigg)^\frac 12 \\
&\hphantom{X}
+ \int_{\R} |(1-\rho(x))Q_\dl(x)|^2 \,\ds x\\
 &\le \eps +O(\exp(-c\dl_0^{-1}))
\end{split}
\end{align*}

\noi
and thus, for  sufficiently small $\dl_0 = \dl_0(\eps) >0$, we  have 
\begin{align*}
\|u-Q_{\dl_0}\|_{L^2(\R)} <  2\eps.
\end{align*}

\noi
Hence,  from the discussion above on the real line case, 
 we have 
\begin{align}
u 
=  e^{i\ta_1}Q_{\dl_1,x_1}+v
\label{J5a}
\end{align}

\noi
for some $(\dl_1,x_1,\ta_1)$ near $(\dl_0,0,0)$ and 
$v\in V^1_{\dl_1,x_1,\ta_1}(\R)$
with  $\|v\|_{L^2(\R)} < \eps_0 \ll 1$. 
Since $u=0$ and $Q_{\dl_1,x_1}=O_{L^2(\R)}(\exp(-c\dl_1^{-1}))$ outside $[-\frac 12,\frac 12]$, we have
\begin{align}
\|v\|_{L^2([-\frac 12,\frac 12]^c)}=O(\exp(-c\dl_1^{-1})).
\label{J5b}
\end{align}
Moreover, from \eqref{rhodef} and \eqref{QQ1}, we have 
\begin{align}
\|e^{i\ta_1}Q_{\dl_1,x_1} -  e^{i\ta_1}Q_{\dl_1,x_1}^\rho\|_{L^2(\R)} 
= O(\exp(-c\dl_1^{-1}))
\label{J5c}
\end{align}

Define  the map $F = F_{\dl_1, x_1,\ta_1} :\R_+\times \R\times
(\R/2\pi\Z) \times V_{\dl_1,x_1,\ta_1}(\T)\to L^2(\T)$ by
\begin{equation*}
F(\dl, x,\ta,v)=e^{i\ta}Q_{\dl,x}^\rho+P_{V_{\dl,x,\ta}}v,
\end{equation*}

\noi
where 
$Q_{\dl,x}^\rho = (\tau_x \rho) Q_{\dl,x}$ and 
$P_{V_{\dl,x,\ta}}$ is the projection onto $V_{\dl,x,\ta}(\T)$
in $L^2(\T)$. 
Then, from~\eqref{J5a}, we claim that 
\begin{equation}\label{eqn: img-dist}
  \|u-F(\dl_1,x_1,\ta_1,v_1)\|_{L^2(\T)}= O(\exp(-c\dl_1^{-1})),
\end{equation}
\noi
where $v_1 = P_{V_{\dl_1,x_1,\ta_1}}(v|_{\T})$.
Note that we have
\begin{align}
\|v_1\|_{L^2(\T)} \les \eps_0 \ll 1. 
\label{J6a}
\end{align}

\noi
Let us verify \eqref{eqn: img-dist}.
For $v \in V_{\dl_1, x_1,\ta_1}^1(\R)$ in \eqref{J5a}, we have
\begin{align}
\begin{split}
& \jb{v, (1-\dx^2)e^{i \ta_1}\dd_{\dl} Q_{ \dl_1, x_1}}_{L^2(\R)}
= \jb{v, (1-\dx^2)e^{i \ta_1} \partial_{x_0} Q_{\dl_1,x_1}}_{L^2(\R)}\\
& \hphantom{X}
= \jb{v, (1-\dx^2)ie^{i\theta_1} Q_{\dl_1, x_1}}_{L^2(\R)} =0.
\end{split}
\label{J7}
\end{align}

\noi
Recalling from Remark \ref{REM:ortho}  
the orthogonality of 
$e^{i\ta_1}\dd_{\dl} Q_{ \dl_1, x_1}^\rho$, 
$e^{i\ta_1}\partial_{x_0}  Q_{\dl_1,x_1}^\rho$,
and $ie^{i\ta_1} Q_{ \dl_1, x_1}^\rho$ in $H^2(\T)$, we have 
\begin{align}
\begin{split}
v_1 
& = P_{V_{\dl_1,x_1,\ta_1}}(v|_{\T})\\
& = v  
-\frac{\jb{v, (1-\dx^2)e^{i\ta_1}\dd_{\dl} Q_{ \dl_1, x_1}^\rho}_{L^2(\T)}}{
\|(1-\dx^2)\partial_{\dl} Q_{\dl_1, x_1}^\rho\|^2_{L^2(\T)}}(1-\dx^2)e^{i\ta_1}\dd_{\dl} Q_{ \dl_1, x_1}^\rho\\
& \hphantom{XX}- \frac{\jb{v, (1-\dx^2)e^{i\ta_1}\dd_{x_0}Q_{ \dl_1, x_1}^\rho}_{L^2(\T)}}{
\|(1-\dx^2)\partial_{x_0} Q_{\dl_1, x_1}^\rho\|^2_{L^2(\T)}}(1-\dx^2)e^{i\ta_1}
\dd_{x_0} Q_{ \dl_1, x_1}^\rho\\
&\hphantom{XX}- \frac{\jb{v, (1-\dx^2)ie^{i\ta_1}Q_{ \dl_1, x_1}^\rho}_{L^2(\T)}}
{\|(1-\dx^2)Q_{ \dl_1, x_1}^\rho\|_{L^2(\T)}^2} (1-\dx^2)i e^{i\ta_1} Q_{ \dl_1, x_1}^\rho.
\end{split}
\label{J7a}
\end{align}

\noi
Then, from the exponential decay of the ground state, 
\eqref{J5b}, \eqref{J5c}, and \eqref{J7}, 
we obtain 
\begin{align}
\|v - v_1\|_{L^2(\T)} = O(\exp(-c\dl_1^{-1}))
\label{J8}
\end{align}

\noi
for $0 < \dl_1 \ll 1$.
Here, we used the polynomial  bounds (in $\dl_1^{-1}$)
on  $\|(1-\dx^2)\dd_{\dl} Q_{\dl_1, x_1}^\rho\|^2_{L^2(\T)}$, 
 $\|(1-\dx^2)\partial_{x_0}Q_{\dl_1, x_1}^\rho\|^2_{L^2(\T)}$, 
 and $\|(1-\dx^2)Q_{\dl_1,x_1}^\rho\|_{L^2}^2$.
See 
\eqref{JA3}, \eqref{JA5}, and \eqref{JA8}  below.
Hence, from \eqref{J5a}, \eqref{J8}, and 
$F(\dl_1,x_1,\ta_1,v_1) = e^{i\ta_1}Q_{\dl_1, x_1} + v_1$, 
we obtain~\eqref{eqn: img-dist}.

By another translation and rotation, we may assume that $x_1 = 0$ and $\ta_1=0$ in \eqref{eqn: img-dist}.
Namely, we have 
\begin{equation}
  \| u-F(\dl_1,0 ,0,v_1 )\|_{L^2(\T)}= O(\exp(-c\dl_1^{-1})).
\label{J9}
\end{equation}

\noi
Hence, 
to finish the proof, 
we  show that  \eqref{J9} 
guarantees that 
$ u $ lies in the image of $F = F_{\dl_1, 0, 0}$.
For this purpose, we recall the following version of 
the inverse function theorem; see \cite[Theorem~26.29]{driver}. 
See also \cite[Lemma~26.28]{driver}.
In the following, $\|\cdot \|$ denotes the operator norm.

\begin{lemma}\label{LEM: ift}
Given Banach spaces $X$ and $Y$, 
let $f:U\to Y$ be a $C^1$-map from an open subset $U\subset X$ 
to $Y$.
Suppose $x_0\in U$  such that
$\ds f(x_0)$ is invertible.
If there exists 
 $R>0$ such that 
  $\cj{B^X(x_0,R)}\subset U$  and
\begin{equation}\label{eqn: alpha-def}
\sup_{x\in B^X(x_0,R)}\|(\ds  f(x_0))^{-1}\ds  f(x)-\Id_X\|=\kk<1,
\end{equation}

\noi
then $f$ is invertible \textup{(}with a $C^1$-inverse\textup{)} on $B^X(x_0,R)$. Moreover, letting $y_0=f(x_0)\in Y$, we have
  \begin{equation*}
B^Y(y_0,r)\subset f(B^X(x_0,(1-\kk)^{-1}\|\ds  f(x_0)^{-1}\|r))
\end{equation*}

\noi
for any $r<(1-\kk)R/\|\ds  f(x_0)^{-1}\|$.
Here, $B^Z(z_0,R)$ denotes the ball of radius $R$ in $Z = X$ or $Y$ centered at $z_0$.

  \end{lemma}

  Our goal is thus to estimate the quantity $\kk$ in \eqref{eqn: alpha-def}, with $f$ replaced by 
  $F = F_{\dl_1, 0,0}$, to conclude from \eqref{J9} that 
  $ u$ 
  is in the image of $F$. 
 We begin by computing $\big(\ds  F(\dl_1,0,0,v_1)\big)^{-1}$ and its norm.

\begin{lemma}\label{LEM:dF0}
There is a constant $C>0$ such that 
\begin{equation*}
\big\|\big(\ds  F_{\dl_1, 0, 0}(\dl_1,0,0, v_1)\big)^{-1}\big\|\le C
\end{equation*}

\noi
for all sufficiently small $\dl_1 > 0$
and $\|v_1\|_{L^2(\T)} \ll 1$. 

\end{lemma}

We assume  Lemma \ref{LEM:dF0} for now
and proceed with the proof of Proposition \ref{PROP:FOC}.
The proofs of this lemma 
and Lemma \ref{LEM:dF} below
will be presented 
at the end of this subsection.

Let $(\dl, x_0,\ta_0, w) \in \R_+\times \R \times (\R/2\pi\Z) \times V_{\dl_1, 0,0}(\T)$
 such that $|\dl-\dl_1|\ll 1$, $|x_0|\ll 1$, $|\ta_0|\ll1$,   and  $\|w\|_{L^2(\T)}\ll 1$. 
 In the following, 
 we compute
\begin{align*}
\big(\ds  F(\dl_1,0,0,v_1)\big)^{-1}\ds  F(\dl,x_0, \ta_0,w)
: 
 &\,  \R\times\R\times \R\times V_{\dl_1,0,0}(\T)\\
 &\to \R\times\R\times \R\times V_{\dl_1,0,0}(\T)
\end{align*}

\noi
and compare it with the identity operator. 
Given  a vector  $\textbf{v}=(\al, \be,\ga, v)\in \R\times\R\times\R\times V_{\dl_1,0,0}(\T)$, 
we have
\begin{align}
\begin{split}
\wt u :\! 
& = \big(\ds  F(\dl,x_0,\ta_0,w)\big) 
\textbf{v}\\
&= \al  e^{i\ta_0}\partial_\dl Q_{\dl,x_0}^\rho
+\be e^{i\ta_0}\partial_{x_0}Q_{\dl,x_0}^\rho
 +\ga ie^{i\ta_0}Q_{\dl,x_0}^\rho\\
& \hphantom{X}+\big(\ds  P_{V_{\dl, x_0,\ta_0}}(\al, \be,\ga)\big)w
+ P_{V_{\dl, x_0,\ta_0}}v, 
\end{split}
\label{J9a}
\end{align}

\noi
where $\big(\ds  P_{V_{\dl, x_0,\ta_0}}(\al, \be,\ga)\big)w$ is given by 
\begin{align}
\big(\ds  P_{V_{\dl, x_0,\theta_0}}(\al, \be,\ga)\big)w = 
\al \dd_\dl P_{V_{\dl, x_0,\theta_0}}w + \be \dd_{x_0} P_{V_{\dl, x_0,\ta_0}}w
+\gamma \partial_{\ta_0}P_{V_{\dl, x_0,\ta_0}} w.
\label{J10}
\end{align}

Suppose that 
for $\wt{\textbf{v}} = (\wt \al,\wt \beta,\wt\ga, \wt v)$, we have
\begin{equation}
  \big(\ds  F(\dl,x_0,\ta_0,w)\big)\textbf{v}=\wt u =\big(\ds  F(\dl_1,0,0,v_1)\big)\wt{\textbf{v}}.
\label{J10a}
\end{equation}

\noi
Namely, we have
\begin{align}
\begin{split}
\wt{\textbf{v}}
& =  (\wt \al,\wt \beta,\wt\ga, \wt v)\\
& =\big(\ds  F(\dl_1,0,0,v_1)\big)^{-1}\wt u=\big(\ds  F(\dl_1,0,0,v_1)\big)^{-1}
\big(\ds  F(\dl,x_0,\ta_0,w)\big)\textbf{v}.
\end{split}
\label{J11}
\end{align}

\noi
Then, in view of the hypothesis in Lemma \ref{LEM: ift}, our goal
is to estimate $\textbf{v} - \wt{\textbf{v}}$.

\begin{lemma}\label{LEM:dF}
There exist  $0 < \dl_1, \eps_0 \ll 1$
such that 
given any $v_1 \in V_{\dl_1, 0, 0}(\T)$ with $\|v_1\|_{L^2(\T)} \les \eps_0 \ll 1$
and 
given any 
 $(\dl, x_0,\ta_0, w) \in \R_+\times \R \times (\R/2\pi\Z)\times  V_{\dl_1,0, 0}(\T)$
with  $|\dl-\dl_1| \ll \dl_1$, $|x_0|\ll 1$, $|\ta_0| \ll 1$,  and  $\|w - v_1\|_{L^2(\T)}\ll 1$, 
 we have
  \begin{align}
 \wt \al &=
 \al + O\Big( A_{\dl_1, v_1, \dl, x_0, \ta_0, w}( \al, \be, \ga, v) \Big) , 
\label{eqn: alpha-alphatilde}\\
  \wt    \be&=
  \be + O\Big(
  A_{\dl_1, v_1, \dl, x_0, \ta_0, w}( \al, \be, \ga, v) \Big) , 
\label{eqn: beta-betatilde}\\
  \wt \ga &=
  \ga + O\Big(\dl_1^{-1}
  A_{\dl_1,v_1,  \dl, x_0, \ta_0, w}( \al, \be, \ga, v) \Big) ,
\label{GG1}\\
 \wt     v  & = v + O_{L^2(\T)}\Big(\dl_1^{-1}
 A_{\dl_1, v_1, \dl, x_0, \ta_0, w}( \al, \be, \ga, v) \Big) \label{eqn: v-w0}
  \end{align}

\noi
for any $(\al, \be, \ga, v) \in \R\times\R \times \R \times V_{\dl_1,0, 0}(\T)$,
where $(\wt \al, \wt \be ,\wt\ga, \wt v)$ is defined by \eqref{J11}
and 
$ A_{\dl_1, v_1, \dl, x_0, \ta_0, w}( \al, \be, \ga, v)$ is defined by 
\begin{align*}
\begin{split}
A_{\dl_1,v_1,  \dl, x_0, \ta_0, w}( \al, \be, \ga, v)
& = \Big(\dl_1^{-1}  (|\dl-\dl_1| + |x_0|) + |\ta_0| +  \exp(- c \dl_1^{-1})+ \|w- v_1\|_{L^2(\T)}\Big)\\
& \hphantom{X}
\times 
  \Big( |\al| +  |\be|+ \dl_1 |\ga| +  \dl_1 \|v\|_{L^2(\T)} + \|v_1\|_{L^2(\T)}\Big) .
\end{split}
\end{align*}

\end{lemma}

Now, given small $\dl_1 > 0$, 
let us choose 
$|\dl-\dl_1|  + |x_0| + |\ta_0|+ \|w - v_1\|_{L^2(\T)} \les \dl_1^3$.
Then,  Lemmas~\ref{LEM:dF0} and ~\ref{LEM:dF}
allow us to  apply Lemma \ref{LEM: ift}
with $R \sim  \dl_1^3$ and $\kk \sim \dl_1$
and 
 conclude that the image of $F= F_{\dl_1, 0,0}$ contains a ball of radius $r \sim\dl_1^{3}$ around $F(\dl_1,0,0,v_1)$. 
 Recalling~\eqref{J9}, we see that $ u$ 
 indeed lies in the image of $F$.
Lastly, we need to choose $\dl_1 = \dl_1 (\eps_1)>0$ sufficiently small such that  
$R = c  \dl_1^3 \leq \eps_1$.
This concludes the proof of  Proposition~\ref{PROP:FOC}.
\end{proof}

We conclude   this subsection by presenting the proofs
 of Lemmas~\ref{LEM:dF0} and~\ref{LEM:dF}.
In the following, $\jb{\,\cdot, \cdot\,}$ denotes the inner product in $L^2(\T)$.

\begin{proof}[Proof of Lemma \ref{LEM:dF0}]

Let $(\al , \beta,\ga, v)  \in \R\times \R\times \R \times  V_{\dl_1,0,0}(\T)$.
Then, with $F = F_{\dl_1, 0,0}$, we have 
\begin{align}
\begin{split}
\wt v :\! & =\big(\ds  F(\dl_1,0,0,v_1)\big)(\al,\beta,\ga,v)\\
& =\al \dd_{\dl}Q_{\dl_1}^\rho+\beta 
\dd_{x_0} Q_{\dl_1}^\rho+\ga i Q_{\dl_1}^\rho\\
& \hphantom{X}
+ \big(\ds  P_{V_{\dl,x_0,\ta_0}}(\al, \be,\ga)\big)v_1|_{(\dl, x_0, \ta_0)
= (\dl_1, 0, 0)}
+v\in L^2(\T),
\end{split}
\label{JA0}
\end{align}

\noi
where 
$\dd_{x_0} Q_{\dl_1}^\rho = 
\partial_{x_0} \big((\tau_{x_0}\rho) Q_{\dl_1, x_0})\big)|_{x_0 = 0}$
is as in \eqref{QQ0}
and 
the fourth term on the right-hand side is as in \eqref{J10}.

\smallskip

\noi
$\bullet$
{\bf Case 1:}
We first consider the case  when $v_1 = 0$.
\\
\indent 
Then, 
by the orthogonality of 
$\dd_{\dl}Q_{\dl_1}^\rho$, $\partial_{x_0} Q_{\dl_1}^\rho$,  and $i Q_{\dl_1}^\rho$
in $H^1(\T)$ (see Remark \ref{REM:ortho})
and the definition~\eqref{eqn: Vdelta-def} of $V_{\dl_1,0, 0}$, we have
\begin{align}
\al&= \frac{\langle \wt v, (1-\dx^2)\partial_\dl Q_{\dl_1}^\rho \rangle }
{\langle  \partial_\dl  Q_{\dl_1}^\rho,(1-\dx^2)\partial_\dl Q_{\dl_1}^\rho\rangle}, \label{eqn: alpha-coord}\\
\beta&= \frac{\langle \wt v,(1-\dx^2)\partial_{x_0} Q_{\dl_1}^\rho\rangle}{\langle \partial_{x_0}
Q_{\dl_1}^\rho,(1-\dx^2)\partial_{x_0} Q_{\dl_1}^\rho\rangle}, \label{eqn: beta-coord}\\
\ga&=  \frac{\jb{\wt v,(1-\dx^2)i Q_{\dl_1}^\rho}}{\jb{ Q_{\dl_1}^\rho,(1-\dx^2) Q_{\dl_1}^\rho}},
 \label{eqn: gamma-coord}\\
v&= \wt v-\al  \partial_\dl Q_{\dl_1}^\rho -\beta \partial_{x_0}  Q_{\dl_1}^\rho - \ga i Q_{\dl_1}^\rho. \label{eqn: v-coord}
\end{align}

From \eqref{eqn: alpha-coord}, we obtain
\begin{align}
|\al|\le \|\wt v\|_{L^2(\T)}\frac{\|(1-\dx^2) \partial_\dl Q_{\dl_1}^\rho\|_{L^2(\T)}}
{|\langle  \partial_{\dl}Q_{\dl_1}^\rho,(1-\dx^2) \partial_\dl Q_{\dl_1}^\rho\rangle|}.
\label{JA1}
\end{align}

\noi
By a direct computation with \eqref{corr1}
and \eqref{Q1}, we have 
\begin{equation}\label{JA2}
\begin{split}
\langle  \partial_{\dl} & Q_{\dl_1}^\rho, (1-\dx^2) \partial_{\dl}Q_{\dl_1}^\rho\rangle \\
& = \langle \rho  \partial_{\dl}  Q_{\dl_1}, (1-\dx^2) (\rho \partial_{\dl}Q_{\dl_1})\rangle \\
& =\dl_1^{-4}\int_{-\frac{1}{2\dl_1}}^{\frac{1}{2\dl_1}}
\rho(\dl_1 x) A_1(x)\\
& \hphantom{X}
\times  \Big(\dl_1^2 \big(\rho(\dl_1 x) - \rho''(\dl_1 x)\big) A_1(x)  - 2 \dl_1 \rho'(\dl_1 x)A_2(x) - A_3(x) \Big) \ds x\\
& =\dl_1^{-4}\int_{\R}
A_1(x) \big(\dl_1^2 A_1(x) -A_3 (x)\big) \ds x
+ O(\exp(-c\dl^{-1}_1))\\
& \sim \dl_1^{-4}
\end{split}
\end{equation}

\noi
for $0 < \dl_1 \ll 1$, 
where 
\begin{align}
\begin{split}
A_1(x)   : & \! = - \frac 12 Q(x) - x Q'(x),\qquad  A_2(x) := - \frac 32 Q'(x) - x Q''(x),\\
& \text{and} \qquad  A_3 (x) : = -\frac 52  Q''(x) - xQ'''(x).
\end{split}
\label{JA2a}
\end{align}

\noi
By a similar computation, we have
\begin{align}
\begin{split}
  \|(1-\dx^2) \partial_{\dl}Q_{\dl_1}^\rho\|^2_{L^2(\T)}
&   =\dl^{-6}_1 \int_\R 
  \big(\dl_1^2 A_1(x) - A_3(x)\big)^2
\ds  x+O\big(\exp(-c\dl^{-1}_1)\big)\\
& \sim \dl_1^{-6}
\end{split}
\label{JA3}
\end{align}

\noi
for $0 < \dl_1 \ll 1$.
Hence, from \eqref{JA1}, \eqref{JA2}, and \eqref{JA3}
with \eqref{JA0}, we obtain
\begin{equation}\label{JA3a}
|\al|\lesssim \dl_1\big\|\big(\ds  F(\dl_1,0,0,0)\big)(\al, \be, \ga, v)\big\|_{L^2(\T)}
\end{equation}

\noi
for $0 < \dl_1 \ll 1$.

Proceeding analogously, we have 
\begin{align}
 \langle \partial_{x_0} Q_{\dl_1}^\rho, (1-\dx^2)\partial_{x_0} Q_{\dl_1}^\rho\rangle &\sim \dl_1^{-4},  \label{JA4}\\
\|(1-\dx^2)\partial_{x_0}  Q_{\dl_1}^\rho\|^2_{L^2(\T)}&\sim \dl_1^{-6}  \label{JA5} 
\end{align}

\noi
for $0 < \dl_1 \ll 1$.
Hence, from  \eqref{eqn: beta-coord}, 
\eqref{JA4}, and \eqref{JA5}, 
we obtain
\begin{align}
|\beta|&\lesssim \dl_1 \big\|\big(\ds  F(\dl_1,0,0,0)\big)(\al, \be,\ga, v)\big\|_{L^2(\T)} \label{JA6}
\end{align}

\noi
for $0 < \dl_1 \ll 1$.
Similarly, we have
\begin{align}
\jb{ Q_{\dl_1}^\rho,(1-\dx^2) Q_{\dl_1}^\rho} & \sim \delta_1^{-2},  \label{JA7}\\
\|(1-\dx^2) Q_{\dl_1}^\rho\|_{L^2(\T)}^2 & \sim \delta_1^{-4} \label{JA8}
\end{align}
for $0 < \delta_1 \ll 1$. Hence, 
from \eqref{eqn: gamma-coord}, 
\eqref{JA7}, and \eqref{JA8}, 
we obtain
\begin{align}
|\gamma|&\lesssim \big\|\big(\ds  F(\dl_1,0,0,0)\big) (\al, \be, \ga, v)\big\|_{L^2(\T)}. \label{JA9}
\end{align}

Lastly, from \eqref{eqn: v-coord}, \eqref{JA3a}, \eqref{JA6}, and \eqref{JA9}
with $\| \dd_\dl Q_{\dl_1}^\rho\|_{L^2(\T)}\sim \|\dd_{x_0}  Q_{\dl_1}^\rho\|_{L^2(\T)}\sim \dl_1^{-1}$ 
and $\| Q_{\dl_1}^\rho\|_{L^2(\T)} \sim 1$, we obtain
\begin{align}
\|v\|_{L^2(\T)}
\lesssim \big\|\big(\ds  F(\dl_1,0,0,0)\big)(\al, \be,\ga, v)\big\|_{L^2(\T)}
\label{JA10}
\end{align}

\noi
\noi
for $0 < \dl_1 \ll 1$.
Therefore, we conclude from \eqref{JA3a}, \eqref{JA6}, \eqref{JA9},  and \eqref{JA10}
that  the inverse $\big(\ds  F(\dl_1,0,0, 0)\big)^{-1}$ of $\ds  F$ at $(\dl_1,0,0, 0)$ has a norm bounded by a constant, uniformly in 
sufficiently small $ \dl_1 >0$.

\medskip

\noi
$\bullet$
{\bf Case 2:}
Next, we  consider the case  when $v_1 \ne 0$
with $\|v_1\|_{L^2(\T)} \ll 1$. 
\\
\indent
In this case, from \eqref{JA0}, we have
\begin{align}
\al&= \frac{\langle \wt v, (1-\dx^2)\partial_\dl Q_{\dl_1}^\rho\rangle }
{\langle  \partial_\dl Q_{\dl_1}^\rho,(1-\dx^2) \partial_\dl Q_{\dl_1}^\rho\rangle}
- \frac{\langle V_1, (1-\dx^2) \partial_\dl Q_{\dl_1}^\rho\rangle }
{\langle  \partial_\dl Q_{\dl_1}^\rho,(1-\dx^2) \partial_\dl Q_{\dl_1}^\rho\rangle}, \label{JA11}\\
\beta&= \frac{\langle \wt v,(1-\dx^2)\partial_{x_0}  Q_{\dl_1}^\rho\rangle}
{\langle \partial_{x_0} Q_{\dl_1}^\rho,(1-\dx^2)\partial_{x_0}  Q_{\dl_1}^\rho\rangle}
- \frac{\langle V_1,(1-\dx^2)\partial_{x_0}  Q_{\dl_1}^\rho\rangle}
{\langle \partial_{x_0} Q_{\dl_1}^\rho,(1-\dx^2)\partial_{x_0}  Q_{\dl_1}^\rho\rangle}, \label{JA12}\\
\ga&=  \frac{\jb{\wt v,(1-\dx^2)i Q_{\dl_1}^\rho}}{\jb{ Q_{\dl_1}^\rho,(1-\dx^2) Q_{\dl_1}^\rho}}
- \frac{\jb{V_1(1-\dx^2)i Q_{\dl_1}^\rho}}{\jb{ Q_{\dl_1}^\rho,(1-\dx^2) Q_{\dl_1}^\rho}}, 
 \label{JA13}\\
v&= \wt v-\al  \partial_\dl Q_{\dl_1}^\rho-\beta \partial_{x_0}  Q_{\dl_1}^\rho - \ga i  Q_{\dl_1}^\rho
- V_1,  \label{JA14}
\end{align}

\noi
where $V_1 = V_1(\al, \be, \g)$ is given by 
\begin{align}
V_1: = \ds  P_{V_{\dl_1,0,0}}(\al, \be,\ga)\big)v_1
= 
\big(\ds  P_{V_{\dl,x_0,\ta_0}}(\al, \be,\ga)\big)v_1|_{(\dl, x_0, \ta_0)
= (\dl_1, 0, 0)}.
\label{JA14a}
\end{align}

From~\eqref{J10} and \eqref{J7a}, we have 
\begin{align}
(\ds &  P_{V_{\dl, x_0,\ta_0}}  (\al,\be,\ga))v_1 \notag\\
& =-\al \partial_\dl \Bigg(\frac{\langle v_1, (1-\dx^2)e^{i\ta_0}\partial_\dl Q_{\dl,x_0}^\rho\rangle}{\|(1-\dx^2) \dd_\dl Q_{\dl,x_0}^\rho\|_{L^2(\T)}^2}(1-\dx^2)e^{i\ta_0}\partial_\dl Q_{\dl,x_0}^\rho \Bigg)
\notag \\
&\hphantom{X}
-\be \partial_{x_0}\Bigg(\frac{\langle v_1, (1-\dx^2)e^{i\ta_0}\partial_\dl Q_{\dl,x_0}^\rho\rangle}{\|(1-\dx^2) \dd_\dl Q_{\dl,x_0}^\rho\|_{L^2(\T)}^2}(1-\dx^2)e^{i\ta_0}\partial_\dl Q_{\dl,x_0}^\rho \Bigg)
\notag \\
&\hphantom{X}
-\ga \partial_{\ta_0}\Bigg(\frac{\langle v_1, (1-\dx^2)e^{i\ta_0}\partial_\dl Q_{\dl,x_0}^\rho\rangle}{\|(1-\dx^2) \dd_\dl Q_{\dl,x_0}^\rho\|_{L^2(\T)}^2}(1-\dx^2)e^{i\ta_0}\partial_\dl Q_{\dl,x_0}^\rho \Bigg)
\notag \\
&\hphantom{X}
-\al \partial_\dl\Bigg(\frac{\langle v_1, (1-\dx^2)e^{i\ta_0}\partial_{x_0} Q_{\dl,x_0}^\rho\rangle}{\|(1-\dx^2)\partial_{x_0}Q_{\dl,x_0}^\rho\|_{L^2(\T)}^2}(1-\dx^2)e^{i\ta_0}
\partial_{x_0} Q_{\dl,x_0}^\rho \Bigg)
\notag\\
&\hphantom{X}
-\be \partial_{x_0}\Bigg(\frac{\langle v_1, (1-\dx^2)e^{i\ta_0}\partial_{x_0} Q_{\dl,x_0}^\rho\rangle}{\|(1-\dx^2)\partial_{x_0}Q_{\dl,x_0}^\rho\|_{L^2(\T)}^2}(1-\dx^2)e^{i\ta_0}\partial_{x_0}
Q_{\dl,x_0}^\rho \Bigg)
\notag\\
&\hphantom{X}
-\ga \partial_{\ta_0}\Bigg(\frac{\langle v_1, (1-\dx^2)e^{i\ta_0}\partial_{x_0} Q_{\dl,x_0}^\rho\rangle}{\|(1-\dx^2)\partial_{x_0}Q_{\dl,x_0}^\rho\|_{L^2(\T)}^2}(1-\dx^2)e^{i\ta_0}\partial_{x_0} Q_{\dl,x_0}^\rho \Bigg)
\notag\\
&\hphantom{X}
-\al \partial_\dl\Bigg(\frac{\langle v_1, (1-\dx^2)i e^{i\ta_0} Q_{\dl,x_0}^\rho\rangle}
{\|(1-\dx^2) Q_{\dl,x_0}^\rho\|_{L^2(\T)}^2}(1-\dx^2)i e^{i\ta_0}Q_{\dl,x_0}^\rho \Bigg)
\notag\\
&\hphantom{X}
-\be \partial_{x_0}\Bigg(\frac{\langle v_1, (1-\dx^2)i e^{i\ta_0} Q_{\dl,x_0}^\rho\rangle}
{\|(1-\dx^2) Q_{\dl,x_0}^\rho\|_{L^2(\T)}^2}(1-\dx^2)i e^{i\ta_0}Q_{\dl,x_0}^\rho \Bigg)
\notag\\
&\hphantom{X}
-\ga \partial_{\ta_0}\Bigg(\frac{\langle v_1, (1-\dx^2)i e^{i\ta_0} Q_{\dl,x_0}^\rho\rangle}
{\|(1-\dx^2) Q_{\dl,x_0}^\rho\|_{L^2(\T)}^2}(1-\dx^2)i e^{i\ta_0}Q_{\dl,x_0}^\rho \Bigg).
\label{JB1}
\end{align}

\noi
The derivatives appearing above are given by 
\begin{align}\label{JB2}
\begin{split}
\partial_{\kk_1}  \Bigg( &  \frac{\langle v_1, (1-\dx^2)\partial_{\kk_2} (e^{i\ta_0}Q_{\dl,x_0}^\rho)\rangle}
{\|(1-\dx^2)\partial_{\kk_2} (e^{i\ta_0}Q_{\dl,x_0}^\rho)\|_{L^2(\T)}^2}
(1-\dx^2)\partial_{\kk_2}(e^{i\ta_0}Q_{\dl,x_0}^\rho) \Bigg)\\
& =\frac{\langle v_1, (1-\dx^2)\partial_{\kk_2} (e^{i\ta_0}Q_{\dl,x_0}^\rho)\rangle}
{\|(1-\dx^2)\partial_{\kk_2} (e^{i\ta_0}Q_{\dl,x_0}^\rho)\|_{L^2(\T)}^2}
(1-\dx^2)\partial^2_{\kk_1, \kk_2}(e^{i\ta_0}Q_{\dl,x_0}^\rho)\\
& \hphantom{X}
+\frac{\langle v_1, (1-\dx^2)\partial^2_{\kk_1, \kk_2}(e^{i\ta_0}Q_{\dl,x_0}^\rho)\rangle}
{\|(1-\dx^2)\partial_{\kk_2} (e^{i\ta_0}Q_{\dl,x_0}^\rho)\|_{L^2(\T)}^2}
(1-\dx^2)\partial_{\kk_2}(e^{i\ta_0}Q_{\dl,x_0}^\rho) \\
& \hphantom{X}- 2 \langle (1-\dx^2)\partial_{\kk_2} (e^{i\ta_0}Q_{\dl,x_0}^\rho) ,
 (1-\dx^2)\partial^2_{\kk_1, \kk_2}(e^{i\ta_0}Q_{\dl,x_0}^\rho)\rangle \\
& \hphantom{X}
\times 
\frac{\langle v_1, (1-\dx^2)\partial_{\kk_2} (e^{i\ta_0}Q_{\dl,x_0}^\rho)\rangle}
{\|(1-\dx^2)\partial_{\kk_2}( e^{i\ta_0}Q_{\dl,x_0}^\rho)\|_{L^2(\T)}^4}
(1-\dx^2)\partial_{\kk_2}(e^{i\ta_0}Q_{\dl,x_0}^\rho)
\end{split}
\end{align}

\noi
for $\kk_1, \kk_2 \in \{\dl, x_0,\ta_0\}$.
By estimating \eqref{JB2}
with \eqref{JA2}, \eqref{JA3}, \eqref{JA4}, \eqref{JA5}, \eqref{JA7}, \eqref{JA8},
 and similar computations  for $(1-\dx^2)\partial^2_{\kk_1, \kk_2} (e^{i\ta_0}Q_{\dl,x_0}^\rho)$, 
it is a simple task to show 
\begin{equation}\label{JB2a}
  \big\|\big(\ds  P_{V_{\dl, x_0,\theta_0}}(\al,\be,\ga)\big) v_1\big\|_{L^2(\T)}
  \lesssim \big( \dl^{-1}(|\al|+|\be|)+|\ga|\big)\|v_1\|_{L^2(\T)}, 
\end{equation}
uniformly in $x_0\in \T$ and  $\ta_0 \in \R/2\pi\Z$.
In particular, with \eqref{JA14a}, we have 
\begin{equation}\label{JA15}
 \|V_1\|_{L^2(\T)}
  \lesssim \big( \dl_1^{-1}(|\al|+|\be|)+|\ga|\big)\|v_1\|_{L^2(\T)}.
\end{equation}

%
%

From \eqref{JA11}, \eqref{JA2}, and \eqref{JA3}, we obtain
\begin{equation}\label{JA16}
|\al|\lesssim \dl_1\big\|\big(\ds  F(\dl_1,0,0,v_1)\big)(\al, \be, \ga, v)\big\|_{L^2(\T)}
+ \dl_1 \|V_1 \|_{L^2(\T)}.
\end{equation}

\noi
Similarly, from \eqref{JA12}, we have 
\begin{equation}\label{JA17}
|\be|\lesssim \dl_1\big\|\big(\ds  F(\dl_1,0,0,v_1)\big)(\al, \be, \ga, v)\big\|_{L^2(\T)}
+ \dl_1 \|V_1 \|_{L^2(\T)}.
\end{equation}

\noi
Compare these with \eqref{JA3a} and \eqref{JA6}.
From \eqref{JA16} and \eqref{JA17} with \eqref{JA15}, 
we then obtain
\begin{equation}\label{JA18}
|\al| +|\be| \lesssim \dl_1\big\|\big(\ds  F(\dl_1,0,0,v_1)\big)(\al, \be, \ga, v)\big\|_{L^2(\T)}
+ \dl_1 |\g| \|v_1\|_{L^2(\T)}, 
\end{equation}

\noi
provided that $\|v_1\|_{L^2(\T)} \ll 1$.
By a similar  computation
with \eqref{JA13},
 we have
\begin{equation}\label{JA19}
|\g|\lesssim \big\|\big(\ds  F(\dl_1,0,0,v_1)\big)(\al, \be, \ga, v)\big\|_{L^2(\T)}
+  \|V_1 \|_{L^2(\T)}.
\end{equation}

\noi
Then, from \eqref{JA19}, \eqref{JA15}, and \eqref{JA18}, we obtain 
\begin{equation}\label{JA20}
|\g|\lesssim \big\|\big(\ds  F(\dl_1,0,0,v_1)\big)(\al, \be, \ga, v)\big\|_{L^2(\T)}.
\end{equation}

\noi
Hence, from \eqref{JA18} and \eqref{JA20}, we obtain
\begin{equation}\label{JA21}
|\al| +|\be| \lesssim \dl_1\big\|\big(\ds  F(\dl_1,0,0,v_1)\big)(\al, \be, \ga, v)\big\|_{L^2(\T)}. 
\end{equation}

Proceeding as in Case 1 with 
\eqref{JA14}, \eqref{JA20}, and \eqref{JA21}, 
we obtain 
\begin{align}
\|v\|_{L^2(\T)}
\lesssim \big\|\big(\ds  F(\dl_1,0,0, v_1)\big)(\al, \be,\ga, v)\big\|_{L^2(\T)}
+   \|V_1 \|_{L^2(\T)}.
\label{JA22}
\end{align}

\noi
Then, using \eqref{JA15}, \eqref{JA20}, and \eqref{JA21}, 
we conclude from \eqref{JA22} that 
\begin{align*}
\|v\|_{L^2(\T)}
\lesssim \big\|\big(\ds  F(\dl_1,0,0,v_1)\big)(\al, \be,\ga, v)\big\|_{L^2(\T)}.
\end{align*}

This completes the proof of Lemma \ref{LEM:dF0}.
\end{proof}

Next, we present the proof of  Lemma \ref{LEM:dF}.

\begin{proof}[Proof of Lemma \ref{LEM:dF}]

From \eqref{J10a}, we have
\begin{align}
\begin{split}
\wt u 
&  =\big(\ds  F(\dl_1,0,0,v_1)\big)(\wt \al,\wt \beta,\wt\ga, \wt v)\\
& =\wt \al  \dd_{\dl}Q_{\dl_1}^\rho+\wt \beta \partial_{x_0} Q_{\dl_1}^\rho+ 
\wt\ga i Q_{\dl_1}^\rho + 
 V_1 + \wt v, 
\end{split}
\label{JB2b}
\end{align}

\noi
where $V_1 = V_1(\wt \al, \wt\be, \wt \g)
= \big( \ds  P_{V_{\dl_1,0,0}}(\wt \al,\wt  \be,\wt \ga)\big)v_1
$ is as in \eqref{JA14a} with $\al, \be, \g$
replaced by $\wt \al, \wt \be, \wt \g$.
Then, using  \eqref{J9a}, 
we can write 
 the first component $\wt \al$ of $\wt{\mathbf{v}}$ as 
\begin{equation}
  \begin{split}
\wt     \al&=
\frac{\langle \wt u ,(1-\dx^2) \dd_{\dl} Q_{\dl_1}^\rho\rangle}
{ \langle \partial_{\dl}  Q_{\dl_1}^\rho,(1-\dx^2) \partial_{\dl}Q_{\dl_1}^\rho\rangle}
- \frac{\langle V_1 ,(1-\dx^2) \dd_{\dl} Q_{\dl_1}^\rho\rangle}
{ \langle \partial_{\dl}  Q_{\dl_1}^\rho,(1-\dx^2) \partial_{\dl}Q_{\dl_1}^\rho\rangle}\\
\\
&=
\frac{1}{ \|  \partial_{\dl} Q_{\dl_1}^\rho\|_{H^1(\T)}^2}
\langle 
\al   e^{i\ta_0}\partial_\dl Q_{\dl,x_0}^\rho
+\be e^{i\ta_0}\partial_{x_0}Q_{\dl,x_0}^\rho + \ga ie^{i\ta_0}Q_{\dl, x_0}^\rho \\
& \hphantom{XXXXXXXX}
+ P_{V_{\dl, x_0,\ta_0}}v ,
(1-\dx^2) \partial_{\dl} Q_{\dl_1}^\rho\rangle\\
& \hphantom{X}
+ \frac{\langle W - V_1 ,(1-\dx^2) \dd_{\dl} Q_{\dl_1}^\rho\rangle}
{ \langle  \partial_{\dl}  Q_{\dl_1}^\rho,(1-\dx^2) \partial_{\dl}Q_{\dl_1}^\rho\rangle}, 
    \end{split}
\label{JB3}
\end{equation}

\noi
where
$W = W(\al, \be, \g) = \big(\ds  P_{V_{\dl, x_0,\ta_0}}(\al, \be,\ga)\big)w$.

The main term of the numerator in the last expression is 
$ \al  \langle   e^{i\ta_0}\partial_\dl Q_{\dl,x_0}^\rho, (1-\dx^2) \partial_\dl Q_{\dl_1}^\rho\rangle$. 
Since $\rho$ and $ Q$ are real-valued, 
we see that 
\begin{align*}
\begin{split}
  \langle  e^{i\ta_0}\partial_\dl Q_{\dl,x_0}^\rho, (1-\dx^2) \partial_\dl Q_{\dl_1}^\rho\rangle 
& = \cos \ta_0 \cdot \langle   \partial_\dl Q_{\dl,x_0}^\rho, (1-\dx^2) \partial_\dl Q_{\dl_1}^\rho\rangle\\
& = \big(1 + O(\ta_0^2)\big)  \langle   \partial_\dl Q_{\dl,x_0}^\rho, (1-\dx^2) \partial_\dl Q_{\dl_1}^\rho\rangle
\end{split}
\end{align*}

\noi
for $|\ta_0|\ll 1$.
Recall from \eqref{Q3} and \eqref{JA2a} that
\begin{align}
\partial_\dl Q_{\dl,x_0}&=\dl^{-\frac 32}A_1(\dl^{-1}(x-x_0))
\qquad \text{and}\qquad 
\partial_\dl Q_{\dl_1} = \dl_1^{-\frac 32}A_1(\dl_1^{-1}x).
\label{JB3aa}
\end{align}

\noi
Without loss of generality, assume that $0 < x_0 \ll 1$.
By 
the mean value theorem, we have
\begin{align}
\dl^{-\frac 32 }\big(A_1(\dl^{-1} (x - x_0)) -A_1 (\dl^{-1}x)\big)
\sim \dl^{-\frac 52} A_1' (\dl^{-1}x)|x_0|
\label{JB3a}
\end{align}

\noi
for $x \in [- \frac 38, \frac 38]$.
Then, by \eqref{JB3aa}, \eqref{JB3a},  
and repeating a computation analogous to~\eqref{JA2}
with~\eqref{QQ1} 
and the mean value theorem applied to $\tau_{x_0}\rho-\rho$, 
we have 
\begin{align}
&\langle \partial_{\dl} Q_{\dl,x_0}^\rho, (1-\dx^2)\partial_\dl Q_{\dl_1}^\rho\rangle \nonumber \\
& =\langle  \partial_{\dl}Q_{\dl}^\rho, (1-\dx^2)\partial_\dl Q_{\dl_1}^\rho\rangle 
+\langle \rho \dl^{-\frac 32 }(A_1(\dl^{-1} (\cdot - x_0)) -A_1 (\dl^{-1}\cdot)), 
(1-\dx^2)\partial_\dl Q_{\dl_1}^\rho\rangle \notag \\
& \quad 
+ \langle (\tau_{x_0}\rho-\rho) \partial_{\dl} Q_{\dl,x_0}, 
(1-\dx^2)\partial_\dl Q_{\dl_1}^\rho\rangle \nonumber \\
& =\langle \partial_\dl Q_{\dl_1}^\rho,(1-\dx^2)\partial_\dl Q_{\dl_1}^\rho\rangle 
+ O\big(\dl_1^{-5}(|x_0|+|\dl-\dl_1|)
\big)
+ O(\exp(- c \dl_1^{-1})). \label{eqn: ddeltadelta0-overlap}
\end{align}

\noi
Hence, from \eqref{JA2} and 
\eqref{eqn: ddeltadelta0-overlap}, we obtain 
\begin{equation}\label{eqn: gamma-1}
\frac{\langle  e^{i\ta_0} \partial_{\dl}Q_{\dl,x_0}^\rho, (1-\dx^2) \partial_\dl Q_{\dl_1}^\rho\rangle}
{ \langle  \partial_{\dl} Q_{\dl_1}^\rho,(1-\dx^2) \partial_{\dl}Q_{\dl_1}^\rho\rangle}
=1+O\big(\ta_0^2 + \dl_1^{-1}(|x_0|+ |\dl-\dl_1|)
+ \exp(- c \dl_1^{-1})\big).
\end{equation}

A similar computation 
with the orthogonality of  
$ \partial_\dl Q_{\dl_1}^\rho$ and $\partial_{x_0}  Q_{\dl_1}^\rho$ in $H^1(\T)$ 
(where 
$\dd_{x_0} Q_{\dl_1}^\rho$ is as in \eqref{QQ0})
gives 
\begin{align}
\langle \partial_{x_0}Q_{\dl,x_0}^\rho, (1-\dx^2) \partial_\dl Q_{\dl_1}^\rho\rangle
&=O\big(\dl_1^{-5}(|x_0|+|\dl-\dl_1|)
+ \exp(- c \dl_1^{-1})\big), \nonumber
\end{align}

\noi
and together with  \eqref{JA2}, we obtain
\begin{align}
\frac{\langle e^{i\ta_0} \partial_{x_0} Q_{\dl,x_0}^\rho, (1-\dx^2) \partial_\dl Q_{\dl_1}^\rho\rangle}
{ \langle \partial_{\dl} Q_{\dl_1}^\rho,(1-\dx^2) \partial_{\dl}Q_{\dl_1}^\rho\rangle}
&=O\big(\dl_1^{-1}(|x_0|+ |\dl-\dl_1|)
+ \exp(- c \dl_1^{-1})\big).
\label{eqn: gamma-2}
\end{align}

\noi
By Cauchy-Schwarz inequality
with  \eqref{JA2}, \eqref{JA3},  and $\|Q_{\delta,x_0}^\rho\|_{L^2(\T)}^2 \sim 1$, 
we have
\begin{equation}
\begin{split}
\bigg|\frac{\langle i  e^{i\ta_0} Q_{\dl,x_0}^\rho, (1-\dx^2) \partial_\dl Q_{\dl_1}^\rho\rangle}
{ \langle  \partial_{\dl} Q_{\dl_1}^\rho,(1-\dx^2) \partial_{\dl}Q_{\dl_1}^\rho\rangle} \bigg|
& \les  \dl_1 |\sin(\ta_0)| 
= O( \dl_1  |\ta_0| ).
\end{split}
\label{eqn: gamma-3}
\end{equation}

From \eqref{Q3} and \eqref{JA2a}, we have
\begin{align*}
\dx^2 \partial_\dl Q_{\dl,x_0}&=\dl^{-\frac 72}A_3(\dl^{-1}(x-x_0))
\qquad \text{and}\qquad 
\dx^2  \partial_\dl Q_{\dl_1} = \dl_1^{-\frac 72}A_3(\dl_1^{-1}x).
\end{align*}

\noi
Then, by \eqref{eqn: Vdelta-def} and repeating a similar computation
as above with the mean value theorem, we have
\begin{align}
\begin{split}
 \bigg| & \frac{\langle P_{V_{\dl, x_0,\theta_0}}v, 
  (1-\dx^2)\partial_\dl Q_{\dl_1}^\rho\rangle}{ \langle  \partial_{\dl} Q_{\dl_1}^\rho,(1-\dx^2) \partial_{\dl}Q_{\dl_1}^\rho\rangle}\bigg|
\\
& \hphantom{X}
 = \bigg|\frac{\langle P_{V_{\dl, x_0,\theta_0}}v, 
(1-\dx^2)(  \partial_\dl Q_{\dl_1}^\rho - e^{i\ta_0} \partial_\dl Q_{\dl, x_0}^\rho\rangle}
{ \langle  \partial_{\dl} Q_{\dl_1}^\rho,(1-\dx^2) \partial_{\dl}Q_{\dl_1}^\rho\rangle}\bigg|
\\
& \hphantom{X} 
\les \big( \dl_1 |\ta_0| + |x_0|+ |\dl-\dl_1|
+ \exp(- c \dl_1^{-1})\big)\|v\|_{L^2(\T)}.
\end{split}
\label{eqn: PVw-overlap}
\end{align}

We now consider the last term in \eqref{JB3}.
We write $W - V_1$ as 
\begin{align}
\begin{split}
W - V_1
& = \big(\ds  P_{V_{\dl, x_0,\ta_0}}(\al, \be,\ga)\big)(w - v_1)\\
& \hphantom{X}
+ \big(\ds  P_{V_{\dl, x_0,\ta_0}}(\al, \be,\ga)
- \ds  P_{V_{\dl_1,0,0}}( \al,  \be, \ga)\big) v_1\\
& \hphantom{X}
+ \big(
 \ds  P_{V_{\dl_1,0,0}}(\al, \be,\ga)
- \ds  P_{V_{\dl_1,0,0}}(\wt \al,\wt  \be,\wt \ga)\big) v_1\\
& =: \1 + \II + \III.
\end{split}
\label{JBY1}
\end{align}

\noi
By Cauchy-Schwarz inequality
with \eqref{JB2a},   \eqref{JA2},  and \eqref{JA3}, we have
\begin{align} 
\begin{split}
\bigg|\frac{\langle \1,  (1-\dx^2) \partial_\dl Q_{\dl_1}^\rho\rangle}
{ \langle  \partial_{\dl} Q_{\dl_1}^\rho,(1-\dx^2) \partial_{\dl}Q_{\dl_1}^\rho\rangle}\bigg|
& \lesssim \frac{\dl_1}{\dl}\big(|\al| + |\be|+\dl|\g| \big)\|w-v_1\|_{L^2(\T)}\\
& \lesssim \big(|\al| + |\be|+\dl |\g| \big)\|w-v_1 \|_{L^2(\T)}.
\end{split}
\label{JBY2}
\end{align}

\noi
Proceeding as above with 
\eqref{JB1},  \eqref{JB2}, 
and the mean value theorem, we have  
\begin{align}
\begin{split}
\|\II\|_{L^2(\T)} 
& \les 
\big(|\al| + |\be|+\dl |\g| \big) \\
& \hphantom{X}
\times  \big( \dl_1^{-1} |\ta_0| + \dl_1^{-2}(|x_0|+ |\dl-\dl_1|)
+ \exp(- c \dl_1^{-1})\big)
\|v_1 \|_{L^2(\T)}.
\end{split}
\label{JBY2a}
\end{align}

\noi
Thus, we obtain
\begin{align} 
\begin{split}
& \bigg| \frac{\langle \II,  (1-\dx^2) \partial_\dl Q_{\dl_1}^\rho\rangle}
{ \langle  \partial_{\dl} Q_{\dl_1}^\rho,(1-\dx^2) \partial_{\dl}Q_{\dl_1}^\rho\rangle}
 \bigg|\\
& \hphantom{X} \les 
\big(|\al| + |\be|+\dl |\g| \big)
 \big(  |\ta_0| + \dl_1^{-1}(|x_0|+ |\dl-\dl_1|)
+ \exp(- c \dl_1^{-1})\big)
\|v_1 \|_{L^2(\T)}.
\end{split}
\label{JBY3}
\end{align}

\noi
Lastly, we consider $\III$ in \eqref{JBY1}.
Recalling that $v_1 \in V_{\dl_1, 0, 0}$, 
it follows from  \eqref{JB1}
with the orthogonality of 
$ \dd_\dl Q_{\dl_1}^\rho$, $\dd_{x_0}  Q_{\dl_1}^\rho$, 
and $i  Q_{\dl_1}^\rho$ in $H^2(\T)$
 that 
\begin{align*}
\jb{\III, &  (1-\dx^2) \dd_\dl Q_{\dl_1}^\rho}  = 
\big\langle \big( \ds  P_{V_{\dl_1,0,0}}(\al -\wt \al, \be - \wt \be,\ga - \wt \g)\big)v_1 , 
(1-\dx^2) \dd_\dl Q_{\dl_1}^\rho\big\rangle
\\
 & =-(\al -\wt \al)  \jb{v_1,  (1-\dx^2) \partial^2_\dl Q_{\dl_1}^\rho}
  -(\be -\wt \be)  \jb{v_1,  (1-\dx^2)\dd_{x_0}  \partial_\dl Q_{\dl_1}^\rho}\\
 & \hphantom{X} -(\g -\wt \g)  \jb{v_1,  (1-\dx^2) i \partial_\dl Q_{\dl_1}^\rho}.
\end{align*}

\noi
Then, from
$\|v_1\|_{L^2(\T)} \les \eps_0 \ll1$, we have
\begin{align} 
\begin{split}
\frac{\langle \III, (1-\dx^2) \partial_\dl Q_{\dl_1}^\rho\rangle}{ \langle  \partial_{\dl} Q_{\dl_1}^\rho,
(1-\dx^2) \partial_{\dl}Q_{\dl_1}^\rho\rangle}
 & =  c_1 (\al -\wt  \al)
+ c_2 ( \be - \wt  \be)
+   c_3 \dl_1 ( \g -\wt  \g), 
\end{split}
\label{JBY4}
\end{align}

\noi
where $c_j = c_j(\eps_0) = O(\eps_0)$, $j = 1, 2, 3$.

Combining 
\eqref{JB3}, \eqref{eqn: gamma-1} - \eqref{JBY2}, 
\eqref{JBY3}, and \eqref{JBY4}, 
we have
\begin{align}
(1 + c_1)( \wt \al - \al) 
+ c_2 (\wt \be - \be)
+  c_3 \dl_1 (\wt \g -\g)
= 
O\big(A_{\dl_1,v_1,  \dl, x_0, \ta_0, w}( \al, \be, \ga, v)\big).
\label{JBY5}
\end{align}

\noi
By a similar computation, we also obtain
\begin{align}
c_4 (\wt \al - \al)
+ (1 + c_5)( \wt \be  - \be) 
+  c_6 \dl_1 (\wt \g - \g)
= 
O\big(A_{\dl_1,v_1,  \dl, x_0, \ta_0, w}( \al, \be, \ga, v)\big), 
\label{JBY6}
\end{align}

\noi
where $c_j = c_j(\eps_0) = O(\eps_0)$, $j = 4, 5, 6$.

As for 
the estimate~\eqref{GG1}, 
it follows from 
\eqref{JB2b} with \eqref{J9a} that 
\begin{equation}
  \begin{split}
\wt     \g&=
\frac{\langle \wt u ,(1-\dx^2)i  Q_{\dl_1}^\rho\rangle}
{ \langle  Q_{\dl_1}^\rho,(1-\dx^2) Q_{\dl_1}^\rho\rangle}
- \frac{\langle V_1 ,(1-\dx^2)i  Q_{\dl_1}^\rho\rangle}
{ \langle  Q_{\dl_1}^\rho,(1-\dx^2) Q_{\dl_1}^\rho\rangle}\\
\\
&=  
\frac{ 1}{\| Q_{\dl_1}^\rho\|^2_{H^1(\T)}}
\langle
\al  e^{i\ta_0}\partial_\dl Q_{\dl,x_0}^\rho+\be e^{i\ta_0}\partial_{x_0}
Q_{\dl,x_0}^\rho + \ga i e^{i\ta_0}Q_{\dl, x_0}^\rho \\
& \hphantom{XXXXXXX}
+ P_{V_{\dl, x_0,\ta_0}}v, 
(1-\dx^2)i  Q_{\dl_1}^\rho \rangle\\
& \hphantom{X}
+  \frac{\langle W - V_1 ,(1-\dx^2)i  Q_{\dl_1}^\rho\rangle}
{ \langle  Q_{\dl_1}^\rho,(1-\dx^2) Q_{\dl_1}^\rho\rangle}.
    \end{split}
\label{JB3c}
\end{equation}

\noi
Then, proceeding as before with the mean value theorem,
\eqref{JA7},  and \eqref{JA8}, 
the main contribution to \eqref{JB3c}
is given by 
\begin{equation}
  \begin{split}
&  \frac{ \jb{\g i e^{i\ta_0}Q_{\dl, x_0}^\rho ,
(1-\dx^2)i  Q_{\dl_1}^\rho}}
{ \langle  Q_{\dl_1}^\rho,(1-\dx^2) Q_{\dl_1}^\rho\rangle}\\
& \hphantom{XXX}
= \g \Big(1 + O\big(\ta_0^2 +  \dl_1^{-1} (|x_0| + |\dl - \dl_1|) + \exp(-c \dl_1^{-1})\big)\Big), 
    \end{split}
\label{JB3d}
\end{equation}

\noi
while 
the contribution to $\wt \g$
from the terms involving $\al$ and $\be$
can be bounded by
\begin{align}
\big(|\al|+|\be|\big)\big(\dl_1^{-1} |\ta_0| + \dl_1^{-2} (|x_0| + |\dl - \dl_1|) + \exp(-c \dl_1^{-1})\big).
\label{JB3e}
\end{align}

\noi
Proceeding as in \eqref{eqn: PVw-overlap}, we have
\begin{align}
\begin{split}
 \bigg|\frac{\langle P_{V_{\dl, x_0,\theta_0}}v 
  ,(1-\dx^2)i  Q_{\dl_1}^\rho\rangle}
{ \langle  Q_{\dl_1}^\rho,(1-\dx^2) Q_{\dl_1}^\rho\rangle}\bigg|
& = \bigg|\frac{\langle P_{V_{\dl, x_0,\theta_0}}v, 
i (1-\dx^2)(  Q_{\dl_1}^\rho -  e^{i\ta_0}  Q_{\dl, x_0}^\rho\rangle}
{ \langle  Q_{\dl_1}^\rho,(1-\dx^2) Q_{\dl_1}^\rho\rangle}\bigg|
\\
& 
\les \big(  |\ta_0| + \dl^{-1} (|x_0|+ |\dl-\dl_1|)
+ \exp(- c \dl_1^{-1})\big)\|v\|_{L^2(\T)}.
\end{split}
\label{JB3g}
\end{align}

By Cauchy-Schwarz inequality
with \eqref{JBY1} and 
\eqref{JB2a}, we have 
%
%
\begin{align} 
\begin{split}
\bigg|\frac{\langle \1, (1-\dx^2)i  Q_{\dl_1}^\rho\rangle}
{ \langle  Q_{\dl_1}^\rho,(1-\dx^2) Q_{\dl_1}^\rho\rangle}\bigg|
& \lesssim
\big( \dl^{-1}(|\al|+|\be|)+|\ga|\big)\|w - v_1\|_{L^2(\T)}.
\end{split}
\label{JBY7}
\end{align}

\noi
From 
\eqref{JBY2a}, we have 
\begin{align} 
\begin{split}
 \bigg|\frac{\langle \II, 
(1-\dx^2)i  Q_{\dl_1}^\rho\rangle}
{ \langle  Q_{\dl_1}^\rho,(1-\dx^2) Q_{\dl_1}^\rho\rangle}
 \bigg|
&  \les 
\big(|\al| + |\be|+\dl |\g| \big) \\
& \hphantom{X}
\times  \big(  \dl_1^{-1} |\ta_0| + \dl_1^{-2}(|x_0|+ |\dl-\dl_1|)
+ \exp(- c \dl_1^{-1})\big)
\|v_1 \|_{L^2(\T)}.
\end{split}
%
\label{JBY8}
\end{align}

\noi
Next, we consider $\III$ in \eqref{JBY1}.
Recalling that $v_1 \in V_{\dl_1, 0, 0}$, 
it follows from  \eqref{JB1}
with the orthogonality of 
$ \dd_\dl Q_{\dl_1}^\rho$, $\dd_{x_0}  Q_{\dl_1}^\rho$, 
and $i  Q_{\dl_1}^\rho$ in $H^2(\T)$
 that 
\begin{align*}
\jb{\III,   (1-\dx^2)i Q_{\dl_1}^\rho}  
 & =-(\al -\wt \al)  \jb{v_1,  (1-\dx^2) i  \partial_\dl Q_{\dl_1}^\rho}
  -(\be -\wt \be)  \jb{v_1,  (1-\dx^2) i  \dd_{x_0}  Q_{\dl_1}^\rho}\\
 & \hphantom{X} + (\g -\wt \g)  \jb{v_1,  (1-\dx^2)   Q_{\dl_1}^\rho}.
\end{align*}

\noi
Then, from
$\|v_1\|_{L^2(\T)} \les \eps_0 \ll1$, we have
\begin{align} 
\begin{split}
\frac{\langle \III, (1-\dx^2)i  Q_{\dl_1}^\rho\rangle}
{ \langle  Q_{\dl_1}^\rho,(1-\dx^2) Q_{\dl_1}^\rho\rangle}
 & =  c_7\dl_1^{-1} (\al -\wt \al)
+ c_8 \dl_1^{-1}(\be -\wt \be)
+  c_9  (\g -\wt \g), 
\end{split}
\label{JBY9}
\end{align}

\noi
where $c_j = c_j(\eps_0) = O(\eps_0)$, $j = 7, 8, 9$.

Hence, from \eqref{JB3c} - \eqref{JBY9}, 
we obtain 
\begin{align}
\begin{split}
c_7   \dl_1^{-1}(\wt \al - \al)
 +  c_8\dl_1^{-1} ( \wt \be  - \be) 
& +  (1+c_9)  (\wt \g - \g)\\
& = 
O\big(\dl_1^{-1} A_{\dl_1,v_1,  \dl, x_0, \ta_0, w}( \al, \be, \ga, v)\big).
\end{split} 
\label{JBY10}
\end{align}

\noi
By solving the system of linear equations
\eqref{JBY5}, \eqref{JBY6}, and \eqref{JBY10} for 
$\wt \al - \al$, $ \wt \be  - \be$, 
and $\wt \g - \g$, 
we obtain 
\eqref{eqn: alpha-alphatilde}, 
\eqref{eqn: beta-betatilde}, 
and 
\eqref{GG1}.

Finally, we turn to \eqref{eqn: v-w0}. 
By  \eqref{JB2b} and \eqref{J9a}, $\wt v$ is given by
\begin{equation}\label{JB4}
  \begin{split}
\wt  v
& = \wt u  -\wt \al  \partial_\dl Q_{\dl_1}^\rho-\wt \beta \partial_{x_0} Q_{\dl_1}^\rho 
- \wt \gamma i  Q_{\dl_1}^\rho- V_1\\
&= \al  e^{i\ta_0}\partial_\dl Q_{\dl,x_0}^\rho+\be e^{i\ta_0}\partial_{x_0}Q_{\dl,x_0}^\rho
 + \gamma i  e^{i\ta_0}Q_{\dl,x_0}^\rho\\
&\hphantom{X}
+ P_{V_{\dl, x_0,\theta_0}}v 
-\wt \al  \partial_\dl Q_{\dl_1}^\rho-\wt \beta \partial_{x_0} Q_{\dl_1}^\rho - \wt \gamma i  Q_{\dl_1}^\rho \\
&\hphantom{X}
+ W- V_1.
\end{split}
\end{equation}

For $v \in V_{\dl_1, 0, 0}(\T)$, 
from \eqref{eqn: Vdelta-def}, 
\eqref{JA3}, \eqref{JA5}, \eqref{JA8}, and computations similar to \eqref{eqn: gamma-1}, 
we have
\begin{equation}
  \begin{split}\label{JB5}
  P_{V_{\dl,x_0,\ta_0}}v
  & =v-\frac{\langle v, (1-\dx^2) e^{i\ta_0}\partial_\dl Q_{\dl,x_0}^\rho \rangle}
  {\|(1-\partial_{x}^2)\partial_\dl Q_{\dl,x_0}^\rho\|_{L^2(\T)}^2}
  (1-\dx^2)e^{i\ta_0}\partial_\dl Q_{\dl,x_0}^\rho\\
&\hphantom{X}-\frac{\langle v, (1-\dx^2)e^{i\ta_0}\partial_{x_0} Q_{\dl,x_0}^\rho \rangle}
{\|(1-\partial_{x}^2)\partial_{x_0} Q_{\dl,x_0}^\rho\|_{L^2(\T)}^2}(1-\dx^2)e^{i\ta_0}\partial_{x_0}Q_{\dl,x_0}^\rho\\
&\hphantom{X}-\frac{\jb{v,(1-\dx^2)i  e^{i\ta_0} Q_{\dl,x_0}^\rho}}
{\norm{(1-\dx^2)Q_{\dl,x_0}^\rho}_{L^2(\T)}^2} (1-\dx^2) ie^{i\ta_0}Q_{\dl,x_0}^\rho \\
  & =v-\frac{\langle v, (1-\dx^2)(e^{i\ta_0}\partial_\dl Q_{\dl,x_0}^\rho -  \partial_\dl Q_{\dl_1}^\rho) \rangle}
  {\|(1-\partial_{x}^2)\partial_\dl Q_{\dl,x_0}^\rho\|_{L^2(\T)}^2}(1-\dx^2)e^{i\ta_0}\partial_\dl Q_{\dl,x_0}^\rho\\
&\hphantom{X}-\frac{\langle v, (1-\dx^2)(e^{i\ta_0}\partial_{x_0} Q_{\dl,x_0}^\rho
 - \partial_{x_0}  Q_{\dl_1}^\rho) \rangle}
{\|(1-\partial_{x}^2)\partial_{x_0}Q_{\dl,x_0}^\rho\|_{L^2(\T)}^2}(1-\dx^2)e^{i\ta_0}
\partial_{x_0} Q_{\dl,x_0}^\rho\\
&\hphantom{X}-\frac{\jb{v,(1-\dx^2)i(e^{i\ta_0} Q_{\dl,x_0}^\rho -  Q_{\dl_1}^\rho)}}
{\norm{(1-\dx^2)Q_{\dl,x_0}^\rho}_{L^2(\T)}^2} (1-\dx^2) ie^{i\ta_0}Q_{\dl,x_0}^\rho \\
    &= v+O_{L^2(\T)}\Big( \big(\dl_1^{-1}( |x_0|+  |\dl-\dl_1|)
+     |\ta_0| 
    + \exp(- c \dl_1^{-1})\big)\|v\|_{L^2(\T)}\Big).
 \end{split}
\end{equation}

\noi
(Note that we could have obtained the same estimate integrating \eqref{JB2a} over the curve 
$(\alpha_s,\beta_s,\gamma_s) = (s(\delta-\delta_1), sx_0, s\ta_0), 
0 \leq s \leq 1$).

From \eqref{QQ1}, \eqref{JB3aa}, and \eqref{eqn: alpha-alphatilde}, 
we have
\begin{equation}\label{JB6}
  \begin{split}
\al &e^{i\ta_0}\partial_\dl Q_{\dl,x_0}^\rho-\wt \al  \partial_\dl Q_{\dl_1}^\rho\\
&   =
\al (\tau_{x_0} \rho) (e^{i\ta_0}\partial_\dl Q_{\dl,x_0}- \partial_\dl Q_{\dl_1})\\
& 
\quad  + \al ((\tau_{x_0} \rho) \partial_\dl Q_{\dl_1} -  \rho \partial_\dl Q_{\dl_1})
+ (\al -\wt \al) \partial_\dl Q_{\dl_1}^\rho\\
& = O_{L^2(\T)}\Big(  
\big( \dl_1^{-2} (|x_0|  + |\dl - \dl_1|) 
+\dl_1^{-1} |\ta_0|
+ \exp(- c \dl_1^{-1})\big)|\al|+ \dl_1^{-1} |\al - \wt \al|\Big)\\
&   = 
O_{L^2(\T)}\Big(\dl_1^{-1}
A_{\dl_1,v_1,  \dl, x_0, \ta_0, w}( \al, \be, \ga, v) \Big).
  \end{split}
\end{equation}

\noi
Similarly, 
from \eqref{QQ0} and \eqref{eqn: beta-betatilde}, we have
\begin{equation}\label{JB7}
  \begin{split}
\be  &e^{i\ta_0} \partial_{x_0}  Q_{\dl,x_0}^\rho-\wt \beta \partial_{x_0}  Q_{\dl_1}^\rho\\
& = \be \big(e^{i\ta_0}\partial_{x_0} \big((\tau_{x_0} \rho)Q_{\dl,x_0}\big)-\partial_{x_0} ((\tau_{x_0} \rho) Q_{\dl_1, 0})\big) \\
&\phantom{=\ }+ \be \big(\partial_{x_0} ((\tau_{x_0} \rho) Q_{\dl_1,0}) 
-
\dd_{x_0} (\tau_{x_0} (\rho Q_\dl))|_{x_0 = 0} \big)
+ (\be-\wt \beta)\partial_{x_0} Q_{\dl_1}^\rho\\
& =  O_{L^2(\T)}\Big( \big( \dl_1^{-2} ( |x_0|  + |\dl - \dl_1|)
+ \dl_1^{-1} |\ta_0|
+ \exp(- c \dl_1^{-1})\big)|\be|
+ \dl_1^{-1}|\be - \wt \be| \Big)\\
&   = 
O_{L^2(\T)}\Big(\dl_1^{-1}
A_{\dl_1, v_1,  \dl, x_0, \ta_0, w}( \al, \be, \ga, v) \Big).
  \end{split}
\end{equation}

\noi
From \eqref{GG1}, we have 
\begin{equation}\label{JB8}
  \begin{split}
\ga  i & e^{i\ta_0}  Q_{\dl,x_0}^\rho -\wt \ga i Q_{\dl_1}^\rho\\
& = \ga i (\tau_{x_0} \rho) (e^{i\ta_0}  Q_{\dl,x_0}-  Q_{\dl_1}) + \ga i ((\tau_{x_0} \rho) Q_{\dl_1} - \rho Q_{\dl_1})
+ (\ga-\wt \ga) i  Q_{\dl_1}^\rho\\
& =  O_{L^2(\T)}\Big( \big( \dl_1^{-1} ( |x_0|  + |\dl - \dl_1|)
+  |\ta_0|
+ \exp(- c \dl_1^{-1})\big)|\ga|
+ |\ga - \wt \ga| \Big)\\
&   = 
O_{L^2(\T)}\Big(\dl_1^{-1}
A_{\dl_1, v_1, \dl, x_0, \ta_0, w}( \al, \be, \ga, v) \Big).
  \end{split}
\end{equation}

Finally, from \eqref{JBY1}, 
\eqref{JB2a}, \eqref{JBY2a}, 
\eqref{eqn: alpha-alphatilde}, 
\eqref{eqn: beta-betatilde}, 
and 
\eqref{GG1}
with $\|v_1\|_{L^2(\T)} \les \eps_0 \ll1$, we have
\begin{align}
\begin{split}
\|W - V_1\|_{L^2(\T)}
&   \lesssim \big( \dl^{-1}(|\al|+|\be|)+|\ga|\big)\|w- v_1\|_{L^2(\T)} \\
& \hphantom{X}
+ 
\big(|\al| + |\be|+\dl |\g| \big) \\
& \hphantom{XX}
\times  \big( \dl_1^{-1} |\ta_0| + \dl_1^{-2} (|x_0|+ |\dl-\dl_1|)
+ \exp(- c \dl_1^{-1})\big)
\|v_1 \|_{L^2(\T)}\\
 & \hphantom{X}
+ 
\big( \dl_1^{-1}(|\al - \wt \al |+|\be- \wt \be |)+|\ga- \wt \g|\big)\|v_1\|_{L^2(\T)} \\
& \les \dl_1^{-1}
 A_{\dl_1, v_1, \dl, x_0, \ta_0, w}( \al, \be, \ga, v) .
\end{split}
\label{JB9}
\end{align}

\noi
Hence, we obtain
\eqref{eqn: v-w0}
from  \eqref{JB4} together with 
\eqref{JB5} - \eqref{JB9}. 
This concludes the proof of Lemma \ref{LEM:dF}.
\end{proof}

\subsection{Orthogonal coordinate system
in the finite-dimensional setting}

Given $N \in \NB$, 
let $\pi_N$ and $E_N$ be as in 
\eqref{EN0} and \eqref{EN1}, respectively.
Given $\eps$, $\dl_*$, $\dl^* > 0$, we define $U_\eps(\dl_*, \dl^*)$ by 
\begin{equation}
\begin{split}
U_{\eps}(\dl_*, \dl^*)=
\big\{u\in L^2(\T): 
\ &  \|u-e^{i \ta} Q_{\dl, x_0}^\rho\|_{L^2(\T)}<\eps \\
&  \text{for some }  \dl_* < \dl < \dl^* , \, x_0\in  \T,
\text{ and } \ta \in \R  \big\}.
\end{split}
\label{Ueps2}
\end{equation}

\noi
By modifying the proof of Proposition \ref{PROP:FOC}, we obtain the following proposition
on an orthogonal coordinate system
in the finite-dimensional setting.

\begin{proposition}\label{PROP:FOCN}
Let $N \in \NB$.
Given small $\eps_1 > 0$, 
there exist $N_0  = N_0 ( \eps_1) \in \NB$,  $\eps = \eps(\eps_1) > 0$,
$ \dl^*=  \dl^*(\eps_1)  > 0$, and $\dl_* = \dl_* (\eps_1, N)>0$
with
\begin{align*}
\lim_{N \to \infty} \dl_*(\eps_1, N) = 0
\end{align*}

\noi
\textup{(}for each fixed $\eps_1 >0$\textup{)}
 such that 
\begin{align*}
\begin{split}
  U_\eps & (   \dl_* ,   \dl^*)\cap E_N 
 \\ 
 & \subset \big\{u\in E_N:\,
 \|u-
e^{i\ta} \pi_N  Q_{\dl,x_0}^\rho \|_{L^2(\T)}<\eps_1, 
\  u-e^{i \ta} \pi_N  Q_{\dl, x_0}^\rho\in V_{\dl, x_0,\theta} (\T) \\
& \hspace{24mm}  \text{for some  }  0  < \dl <  \dl^* , \, x_0\in \T, \text{ and }\ta \in \R  
\big\}.
\end{split}
\end{align*}

\end{proposition}

Before proceeding to the proof of Proposition \ref{PROP:FOCN}, let us first discuss properties of truncated solitons.
Note that the frequency truncation operator
$\pi_N$ is parity-preserving;
 $\pi_N Q_\dl^\rho = \pi_N(\rho Q_\dl)$ is an even function for any $\dl > 0$. 
It also follows from \eqref{Q1} and \eqref{Q2} that 
$\pi_N \dd_\dl Q_\dl^\rho $ is an even function, 
while
 $\pi_N \dd_{x_0}  Q_\dl^\rho$ 
is an  odd function.
Hence, they are orthogonal in $H^k(\T)$, $k \in \Z_{\geq 0}$. Moreover, the operator $\pi_N$ also commutes with the pointwise conjugation, so $\pi_N Q_\dl^\rho$, $\pi_N  \dd_\dl Q_\dl^\rho$, 
and $\pi_N \dd_{x_0} Q_\dl^\rho$ are all real functions. Therefore, $\pi_N iQ_\dl^\rho$ is orthogonal\footnote{Recall that we view $H^k(\T)$
as a Hilbert space over reals.}  to both $\pi_N  \dd_\dl Q_\dl^\rho$ and $\pi_N \dd_{x_0}  Q_\dl^\rho$ in $H^k(\T)$, $k \in \Z_{\geq 0}$.
By a similar consideration centred at $x = x_0$, 
we also conclude that 
$\pi_N e^{i\ta} \dd_\dl Q_{\dl, x_0}^\rho$, 
 $\pi_N e^{i\ta} \dd_{x_0} Q_{\dl, x_0}^\rho$,  and 
 $\pi_N i  e^{i\ta} Q_{\dl,x_0}^\rho$ are pairwise
 orthogonal in $H^k(\T)$, $k \in \Z_{\geq 0}$.

In the following, we use $\F_\R$ to denote the Fourier transform of a function on the real line.
By \eqref{rhodef},  a change of variables, and the exponential decay of $Q$, we have
\begin{align*}
\ft { Q_\dl^\rho}(n) 
& = \dl^{-\frac 12} \int_{-\frac 12}^{\frac 12} 
\rho(x) Q(\dl^{-1} x) e^{-2\pi i n x} dx
= \dl^{\frac 12} \int_{\R}
\rho(\dl x)Q(x) e^{-2\pi i (\dl n) x} dx\\
& = \dl^{\frac 12} \F_\R(Q)(\dl n) 
+O(\exp(-c\dl^{-1}))\\
& \sim \dl^{\frac 12} \F_\R(Q)(\dl n)
\end{align*}

\noi
for $0 < \dl \ll 1$, provided that $\dl |n| \les 1$.
With  $A_1$ as in \eqref{JA2a}, 
an analogous computation yields
\begin{align}
\begin{split}
\ft {\dd_\dl Q_\dl^\rho}(n) 
& 
= \dl^{- \frac 12} \F_\R(A_1)(\dl n) 
 +O(\exp(-c\dl^{-1}))\\
& \sim  \dl^{- \frac 12} \F_\R(A_1)(\dl n) 
\end{split}
\label{JC2b}
\end{align}

\noi
for $0 < \dl \ll 1$, provided that $\dl |n| \les 1$.

Fix small $\g > 0$
and set 
$\dl_* = \dl_*(N) >0 $
such that 
\begin{align}
(\dl_*)^{1+ \g}  \ges N^{-1}
\label{JC2a}
\end{align} 

\noi
Then, 
by a Riemann sum approximation with \eqref{JC2b}, 
 we then have 
\begin{align*}
\begin{split}
\dl^{2k + 2} \| \pi_N  \dd_\dl Q_{\dl}^\rho\|_{\dot H^k(\T)}^2
& \ges   \sum_{\substack{|n|\leq N\\|n|\les \dl^{-1}}}  (\dl n)^{2k} 
\dl  |\F_\R(A_1)(\dl n) |^2 
\sim \| \pi^\R_{\min(\dl N, 1)} A_1\|_{\dot H^k(\R)}^2 \\
& \geq c_k >0
\end{split}
\end{align*}

\noi
for $\dl_* < \dl \ll 1$,
where $\pi^\R_N$ denotes the Dirichlet projection
onto frequencies $\{|\xi|\leq N\}$ for functions on the real line.
Since the estimate above holds independently of the base point $x_0 \in \T$, we have 
\begin{align}
 \| \pi_N  \dd_\dl Q_{\dl, x_0}^\rho\|_{ H^k(\T)}^2 \ges \dl^{- 2k - 2}
\label{JC3}
\end{align}

\noi
uniformly 
for $\dl_* < \dl \ll 1$ and $x_0 \in \T$. 
On the other hand, 
by integration by parts $2K$ times
together with the exponential decay of the ground state and \eqref{JC2a}, we have
\begin{align}
\begin{split}
\ft {\dd_\dl Q_\dl^\rho}(n) 
& 
= \dl^{- \frac 12} \int_{\R}
\rho(\dl x) A_1(x) e^{-2\pi i (\dl n) x} dx\\
& = O \big(\dl^{-2K - \frac 12} |n|^{-2K}\big).
\end{split}
\label{JC3a}
\end{align}

\noi
Then, 
with $\pi_N^\perp = \Id -  \pi_N$, 
it follows from \eqref{JC2a} and \eqref{JC3a} that 
\begin{align}
\begin{split}
 \| \pi_N^\perp \rho \dd_\dl Q_{\dl, x_0}\|_{ H^k(\T)}^2
& \les \dl^{K}
\end{split}
\label{JC3b}
\end{align}

\noi
for any $K \in \NB$, 
uniformly 
in  $\dl_* < \dl \ll 1$ and $x_0 \in \T$. 
Similarly, we have 
\begin{align}
\begin{split}
 \| \pi_N \dd_{x_0} Q_{\dl, x_0}^\rho\|_{ H^k(\T)}^2 & \ges \dl^{- 2k - 2}, \\
   \| \pi_N^\perp \dd_{x_0} Q_{\dl, x_0}^\rho\|_{ H^k(\T)}^2
&   \les \dl^{K}, \\
\| \pi_N  Q_{\delta,x_0}^\rho\|_{H^k(\T)}^2 & \ges \dl^{-2k}, \\
 \| \pi_N^\perp  Q_{\dl, x_0}^\rho\|_{ H^k(\T)}^2
&   \les \dl^{K}
 \end{split}
\label{JC4}
\end{align}

\noi
for any $k, K \in \NB$, 
uniformly 
in  $\dl_* < \dl \ll 1$ and $x_0 \in \T$.

\begin{proof}[Proof of Proposition \ref{PROP:FOCN}]
The proof of this proposition is based on a small modification of the proof of Proposition \ref{PROP:FOC}.
We only go over the main steps, indicating required modifications.
By a translation and a rotation, we may assume that $u\in E_N $ satisfies  
\[\|u-\rho Q_{\dl_0}\|_{L^2(\T)}<  \eps.\]

\noi
From the real line case discussed in the proof of Proposition \ref{PROP:FOC}
(see \eqref{J5a} above), 
we have 
\begin{align}
u 
=  e^{i\ta_1}Q_{\dl_1,x_1}+v,
\label{JC5a}
\end{align}

\noi
for some $(\dl_1,x_1,\ta_1)$ near $(\dl_0,0,0)$ and 
\begin{align}
v\in V^1_{\dl_1,x_1,\ta_1}(\R)
\quad 
\text{with} \quad \|v\|_{L^2(\R)} < \eps_0 \ll 1,
\label{JC5b}
\end{align} 
provided that   $\dl_0 = \dl_0(\eps_0) > 0$ is  sufficiently small.

Given $N \in \NB$, 
define  the map $F^N = F^N_{\dl_1, x_1,\ta_1} 
:\R_+\times \R\times(\R/2\pi\Z) \times V_{\dl_1,x_1,\ta_1}(\T)\cap E_N \to L^2(\T)$ by
\begin{equation}
\label{JC6}
F^N(\dl, x,\ta, v)=\pi_N e^{i\ta} Q_{\dl,x}^\rho+P_{V_{\dl,x,\ta}\cap E_N }v,
\end{equation}

\noi
where $P_{V_{\dl,x,\ta}\cap E_N }$ is the projection onto $V_{\dl,x,\ta}\cap E_N $  in $L^2(\T)$. 
From \eqref{JC5a} and \eqref{JC6}
with $v_1 = P_{V_{\dl_1,x_1,\ta_1}\cap E_N}(v|_{\T})$, we have 
\begin{equation}
\begin{split}
  u-F^N(\dl_1,x_1,\ta_1, v_1)
&   = 
     e^{i\ta_1}(Q_{\dl_1,x_1} - \pi_N  Q_{\dl_1,x_1}^\rho)
   + (v - v_1).
\end{split}
\label{JC7}
\end{equation}

\noi
From 
the orthogonality of 
$\pi_N e^{i\ta_1}\dd_{\dl} Q_{ \dl_1, x_1}^\rho$, 
$\pi_N e^{i\ta_1}\partial_{x_0} Q_{\dl_1,x_1}^\rho$,
and $\pi_N ie^{i\ta_1} Q_{ \dl_1, x_1}^\rho$ in $H^2(\T)$, we have 
\begin{align}
\begin{split}
v_1 
& = P_{V_{\dl_1,x_1,\ta_1}\cap E_N }(v|_{\T})\\
& = \pi_N v  
-\frac{\jb{\pi_N v, (1-\dx^2)e^{i\ta_1}\dd_{\dl} Q_{ \dl_1, x_1}^\rho}_{L^2(\T)}}{
\|(1-\dx^2) \pi_N\partial_{\dl}  Q_{\dl_1, x_1}^\rho\|^2_{L^2(\T)}}
(1-\dx^2)\pi_N e^{i\ta_1} \dd_{\dl} Q_{ \dl_1, x_1}^\rho\\
& \hphantom{XX}- \frac{\jb{\pi_N v, (1-\dx^2)e^{i\ta_1}\dd_{x_0} Q_{ \dl_1, x_1}^\rho}_{L^2(\T)}}{
\|(1-\dx^2)\pi_N \partial_{x_0} Q_{\dl_1, x_1}^\rho\|^2_{L^2(\T)}}
(1-\dx^2)\pi_N e^{i\ta_1}\dd_{x_0} Q_{ \dl_1, x_1}^\rho\\
&\hphantom{XX}- \frac{\jb{\pi_N v, (1-\dx^2)i e^{i\ta_1}Q_{ \dl_1, x_1}^\rho}_{L^2(\T)}}
{\|(1-\dx^2)\pi_N Q_{ \dl_1, x_1}^\rho\|_{L^2(\T)}^2} 
(1-\dx^2)\pi_N i e^{i\ta_1} Q_{ \dl_1, x_1}^\rho.
\end{split}
\label{JC7a}
\end{align}

\noi
Then, 
we can proceed as in the proof of Proposition \ref{PROP:FOC}
with \eqref{JC3}, \eqref{JC3b}, \eqref{JC4}, 
and \eqref{JC5b}
to estimate  
$\pi_N v - v_1$.
For example, 
the second term on the right-hand side of~\eqref{JC7a}
can be written as 
\begin{align}
\begin{split}
& \frac{\jb{ v, (1-\dx^2)e^{i\ta_1}\dd_{\dl} Q_{ \dl_1, x_1}^\rho}_{L^2(\T)}}{
\|(1-\dx^2) \pi_N \partial_{\dl}  Q_{\dl_1, x_1}^\rho\|^2_{L^2(\T)}}
(1-\dx^2)\pi_N e^{i\ta_1} \dd_{\dl} Q_{ \dl_1, x_1}^\rho\\
& -  \frac{\jb{ v, (1-\dx^2)\pi_N^\perp e^{i\ta_1}\dd_{\dl} Q_{ \dl_1, x_1}^\rho}_{L^2(\T)}}{
\|(1-\dx^2) \pi_N\partial_{\dl}  Q_{\dl_1, x_1}^\rho\|^2_{L^2(\T)}}
(1-\dx^2)\pi_N e^{i\ta_1} \dd_{\dl} Q_{ \dl_1, x_1}^\rho.
\end{split}
\label{JC7aa}
\end{align}

\noi
Thanks to \eqref{JC3}, \eqref{JC5b}, 
and the exponential decay of the ground state, 
the first term in~\eqref{JC7aa} is 
bounded as 
$O_{L^2(\T)} \big(\exp(-c\dl_1^{-1})\big)$.
On the other hand, 
from \eqref{JC3} and \eqref{JC3b}, 
we can bound  the second term in \eqref{JC7aa}
as  $O_{L^2(\T)} (\dl_1^K)$
for any $K \in \NB$.
By estimating the third and fourth terms
on the right-hand side of \eqref{JC7a}
in an analogous manner, we obtain
\begin{align}
\|\pi_N v - v_1\|_{L^2(\T)} = O(\dl_1^K) 
\label{JC7b}
\end{align}

\noi
for any $K \in \NB$.

From \eqref{JC5a},  $u \in E_N$, 
and \eqref{JC4}, 
we have
\begin{align}
\|\pi_N^\perp v \|_{L^2(\T)} = O(\dl_1^K) 
\label{JC7c}
\end{align}

\noi
for any $K \in \NB$.
Putting \eqref{JC7b}
and \eqref{JC7c} together, we obtain
\begin{align}
\|v - v_1\|_{L^2(\T)} = O(\dl_1^K) 
\label{JC8}
\end{align}

\noi
for  $\dl_* < \dl_1 \ll 1$, where $\dl_* = \dl_*(N)$ satisfies \eqref{JC2a}.
Finally, we conclude from \eqref{JC7}, 
\eqref{J5c}, \eqref{JC4},  and \eqref{JC8}
that 
\begin{equation}\label{JC11}
  \|u-F^N(\dl_1,x_1,\ta_1,v_1)\|_{L^2(\T)}= O(\dl_1^{K})
\end{equation}

\noi
for  $\dl_* < \dl_1 \ll 1$, where $\dl_* = \dl_*(N)$ satisfies \eqref{JC2a}.

By our choice of $\dl_* = \dl_*(N)$
such that \eqref{JC3} and  \eqref{JC4}  hold, 
we see that 
Lemmas~\ref{LEM:dF0}
and~\ref{LEM:dF} applied to $F^N = F^N_{\dl_1, 0, 0}$
hold uniformly 
for  $\dl_* < \dl_1 \ll 1$, 
since \eqref{JC3} and  \eqref{JC4} 
provide lower bounds on the denominators of the various terms appearing
in the proofs of these lemmas.
This allows us to apply the inverse function theorem (Lemma \ref{LEM: ift})
with  $R \sim  \dl_1^3$,  $\kk \sim \dl_1$, 
and $r \sim \dl_1^3$.
By taking $K > 3$ in \eqref{JC11}, 
we conclude that  $u$ lies in the image of $F^N = F^N_{\dl_1, x_1, \ta_1}$.
Lastly, we need to choose 
$\dl_* = \dl_* (\eps_1, N)>0$ and 
$\dl^* = \dl^* (\eps_1)>0$ sufficiently small such that  
$  \dl_*< \dl^* \les\eps_1^\frac{1}{3}$.
This concludes the proof of  Proposition~\ref{PROP:FOCN}.
\end{proof}

\subsection{A change-of-variable formula}
\label{SUBSEC:change}

Our main goal is to prove
the bound \eqref{eqn: to-bound}.
In the remaining part of the paper, 
we fix small $\dl^*>0$
and set 
\[U_\eps = U_\eps(\dl^*),\]

\noi 
where 
$ U_\eps(\dl^*)$ is defined in \eqref{Ueps}
for given small $\eps > 0$.
Recalling the low regularity of 
the Ornstein-Uhlenbeck loop in \eqref{bloop3}, 
we only work with functions $u \in 
H^s(\T)\setminus H^\frac 12 (\T)$, $s<\frac 12 $, 
for the $\mu$-integration.
Thus, with a slight abuse of notations, 
 we redefine $U_\eps = U_\eps(\dl^*)$
 in~\eqref{Ueps}
to mean $U_\eps \setminus H^\frac 12 (\T)$.
Namely, when we write 
$U_\eps = U_\eps(\dl^*)$ in the following, it is understood that 
we take an intersection with $\big(H^\frac 12 (\T)\big)^c$.

Given small $\dl>0$ and $ x_0 \in \T$, 
let $ V_{\dl,x_0,0} =  V_{\dl,x_0,0}(\T)$
be as in 
\eqref{eqn: Vdelta-def}
with $\ta = 0$.
In view of 
Proposition \ref{PROP:FOC}, 
the convention above: $U_\eps = U_\eps \cap \big(H^\frac 12 (\T)\big)^c$, 
and the fact that $Q_{\dl, x_0}^\rho \in H^1(\T)$, 
we redefine 
$ V_{\dl,x_0,0}$
to mean $ V_{\dl,x_0,0} \cap H^1(\T)$.
Namely, when we write 
$ V_{\dl,x_0,0}$
 in the following, it is understood that 
we take an intersection with $H^1(\T)$.
Now, let 
 $\mu^\perp_{\dl, x_0}$ denote the Gaussian measure with 
  $V_{\dl,x_0,0} \subset H^1(\T)$
    as its 
 Cameron-Martin space.
Then, we have 
 the following lemma on a change of variables for the $\mu$-integration over $U_\eps$.

\begin{lemma}\label{LEM:surface}
Fix $K > 0$ and sufficiently small $\dl^* > 0$.
Let $F(u)\ge 0 $ be a functional of $u$ on $\T$ 
which is continuous in some topology $H^s(\T)$, $s<\frac 12 $, 
such that $F\le C$ for some $C$ and  $F(u) = 0$ if $\|u\|_{L^2(\T)} > K$. 
Then,  there is a locally finite measure $\d \sigma(\dl)$ on $ (0,\dl^*)$ such that
\begin{align}
\begin{split}
 \int_{U_\eps}  F(u)\,\mu(\ds  u)
&  \le
\iiiint_{U_\eps}
  F\Big(e^{i\ta} (Q_{\dl, x_0}^\rho+v)\Big)\\
& \hphantom{XXX}
\times e^{-\frac{1}{2}\| Q_{\dl, x_0}^\rho\|^2_{ {H}^1(\T)}-\jb{ (1-\dx^2) Q_{\dl,x_0}^\rho,v}_{L^2(\T)}} 
\mu_{\dl, x_0}^\perp(\ds v)\d \sigma( \dl)\d x_0\d \ta, 
\end{split}
\label{LE0}
\end{align}

\noi
where
the domain of the integration on the right-hand side
is to be interpreted as 
\[
\Big\{(\dl, x_0, \ta, v) 
\in (0, \dl^*)\times \T \times (\R/2\pi \Z)
\times V_{\dl,x_0,0}(\T):
e^{i\ta} (Q_{\dl, x_0}^\rho+v) \in U_\eps\Big\}.
 \]
 
 \noi
Moreover, $\sigma$ is absolutely continuous 
with respect to the  Lebesgue measure $\ds  \dl$ on $(0, \dl^*)$, 
satisfying $|\frac{\d\sigma}{\d \dl}| \lesssim \delta^{-20}$.

\end{lemma}
The proof of Lemma \ref{LEM:surface} is based on a finite dimensional approximation and an application of the following lemma.

\begin{lemma}\label{LEM:meas}Let $M\subset \R^n$ be a closed submanifold of dimension $d$ and $\mathcal{N}$ be its normal bundle. 
Then, there is a neighborhood $U$ of $M$ such that $U$ is diffeomorphic to a subset $\M$
of $\mathcal{N}$ 
via the map $\phi:   \N  \to \R^n$ given by 
\begin{align}
\phi(x,v)= x+v
\label{CV1}
\end{align}

\noi
for $x\in M$ and  $v\in T_xM^\perp$. 
Furthermore,  the following estimate holds
for any non-negative measurable function $f: \R^n \to \R$ and any open set 
$V \subset \{(x,v): |v| \le 1\} \subset \mathcal{N}$\textup{:}
\begin{align}
\int_{\phi(V)} f(z) \d z \leq C_d \int_M \left( \int_{T_{x}M^\perp} f(x+v)  \mathbf{1}_{V}(x,v)\d v \right) \d \sigma(x).
\label{LE1}
\end{align}

\noi
Here, 
the measure $\ds \s$ is defined by 
\[\d\sigma(x) = \Big(1 + \sup_{k = 1, \dots, d}\|\nabla t_k(x)\|^d\Big) \d \omega(x),\]

\noi 
 where $\d \omega(x)$ is the surface measure on $M$
and  $ \{t_k(x)\}_{k = 1}^d = \{(t_k^1(x), \dots, t_k^n(x))\}_{k = 1}^d $ is  an orthonormal basis 
for the tangent space $T_x M$
with  
 the expression $\norm{\nabla t_k(x)}$  defined
 by 
 \[ \norm{\nabla t_k(x)}
 = \bigg(\sum_{j = 1}^n \sum_{i = 1}^d \Big|\frac{\dd}{\dd y_i} (t^j_k\circ \varphi^{-1})(y)\Big|^2\bigg)^\frac{1}{2}\]

\noi
for any coordinate chart $\varphi$ on a neighborhood of $x \in M$
such that $y = (y_1, \dots, y_d) = \varphi(x)$
and
$d (\varphi^{-1})(y)$ is an isometry with its image.
Note that the constant in \eqref{LE1} is independent
of  the dimension $n$ of the ambient space~$\R^n$.

\end{lemma}

\begin{proof}
Given $x \in M$, let $w_1(x), \dotsc, w_{n-d}(x)$ be an orthonormal basis of $T_{x}M^\perp$. Consider the map
\begin{align}
\psi: M_x \times \R_\alpha^{n-d} \to  \R^n, 
\qquad \psi(x, \alpha) = x + \sum_{j= 1}^{n-d} \alpha_j w_j(x),
\label{CV2}
\end{align}

\noi
where $\al = (\al_1, \dots, \al_{n-d})$.
Recalling that $\{w_j(x)\}_{j=1}^{n-d}$ is an orthonormal basis of  $T_{x}M^\perp$, 
it follows from the area formula (see \cite[Theorem 2 on p.\,99]{EG})
 with \eqref{CV1}, \eqref{CV2}, 
and $V \subset \{(x,v): |v| \le 1\}$
and then 
applying   a change of variables 
that 
\begin{align*}
& \int_{\phi(V)}  f(z) \d z \\
& \hphantom{X}
\le  \int_{\phi(V)} \#\bigg\{(x,\alpha) \in M \times \R^{n-d}: 
\Big(x,\sum_{j=1}^{n-d} \al_j w_j(x)\Big) \in V, \ \psi(x,\alpha) = z \bigg\}  f(z) \d z\\
& \hphantom{X}
 =  \iint_{M_x \times \{\al \in \R^{n-d}: |\al|\leq 1\}} f(\psi(x, \al)) \ind_V
 \bigg(\Big(x,\sum_{j=1}^{n-d} \alpha_j w_j(x)\Big)\bigg) |J_\psi (x,\alpha)|  \d \alpha \d \omega(x)\\
& \hphantom{X}
 \le  \int_M 
\left( \int_{T_{x}M^\perp} f(x+v)  \ind_{V}(x,v)\d v \right) 
\sup_{|\alpha|\le 1}
|J_\psi (x,\alpha)| \d \omega(x), 
\end{align*}

\noi
where $J_\psi$ is the determinant of the differential of the map $\psi$.
Hence, the bound \eqref{LE1} follows once we  prove 
\begin{align}
|J_\psi (x,\alpha)|  \le C_d \Big(1 + \sup_{ k = 1, \dots, d}\|\nabla t_k(x)\|\Big)^d 
\label{CV3}
\end{align}

\noi
for every $\al \in \R^{n-d}$ with $|\alpha| \le 1$, 
where  $ \{t_k(x)\}_{k = 1}^d $  is  an orthonormal basis 
for the tangent space $T_x M$.

 Recall that the tangent space of $M \times \R^{n-d}$ 
 at the point $(x,\alpha)$ is isomorphic to 
$T_xM \times \R^{n-d}$.
Then, denoting by $\{e_j\}_{j = 1}^{n-d}$ the standard basis of $\R^{n-d}$, 
it follows from \eqref{CV2} that\footnote{Hereafter, we suppress
the $x$-dependence of $w_j = w_j(x)$ and $t_k = t_k(x)$
when there is no confusion.}
\begin{align*}
\d\psi[t_k] &= t_k + \sum_{j = 1}^{n-d} \alpha_j \d w_j[t_k], \qquad k = 1, \dots, d, \\
\d \psi[e_j] & = w_j, \qquad j = 1, \dots, n - d.
\end{align*}

\noi
Hence, by taking 
$(t_1,\dots,t_d, e_1,\dotsc,e_{n-d})$ 
(and 
$(t_1,\dots,t_d, w_1,\dots,w_{n-d})$, respectively)
as the (orthonormal) basis of 
the domain $T_xM \times \R^{n-d}$ 
(and the codomain  $\R^n$, respectively), 
the matrix representation of $\d \psi$ is given by 
$$
A(x,\alpha) = 
\begin{pmatrix}
\Id_{d\times d} + B & 0 \\
D & \Id_{(n-d) \times (n-d)}
\end{pmatrix},
$$
where $B = B(\al) = \{B(\al)_{h, k}\}_{1\le h, k \le d}$ is given by 
\begin{align}
B(\al)_{h,k} = \sum_{j= 1}^{n-d} \alpha_j \jb{\d w_j[t_k] ,    t_h}_{\R^n}.
\label{CV4a}
\end{align}

\noi
Thus, we have
$$|J_\psi (x,\alpha)| = |\det A(x,\alpha)| = |\det (\Id_{d\times d} + B(\al))| 
\les_d1 + \sup_{1\le h,k\le d }|B(\al)_{h,k}|^d$$

\noi
for any  $\al \in \R^{n-d}$ with $|\alpha| \le 1$.
Therefore, 
the bound \eqref{CV3} (and hence \eqref{LE1}) follows once we prove
\begin{align}
\sup_{1\le h,k\le d}|B(\al)_{h,k}| \le \sup_{k = 1, \dots, d} \|\nabla t_k(x)\|, 
\label{CV4}
\end{align}

\noi
uniformly in 
 $\al \in \R^{n-d}$ with $|\alpha| \le 1$.

By differentiating 
the orthogonality relation $\jb{w_j(x),  t_k(x) }_{\R^n}= 0$, we obtain 
\begin{align}
\jb{\d w_j[\tau] ,  t_k }_{\R^n}= - \jb{w_j ,  \d t_k[\tau]}_{\R^n}
\label{CV5}
\end{align}

\noi
for any  $\tau \in T_xM$.
Thus, from \eqref{CV4a}, \eqref{CV5}, 
Cauchy-Schwarz inequality with $|\al | \le 1$,  
and the orthonormality of  $\{w_j\}_{j= 1}^{n-d}$, 
  we have 
\begin{align*}
|B(\al)_{h,k}| 
&= \bigg|\sum_{j=1}^{n-d} \alpha_j \jb{\d w_j[t_k] , t_h}_{\R^n}\bigg|\\
&= \bigg|\sum_{j=1}^{n-d} \alpha_j \jb{w_j,  \d t_h[t_k]}_{\R^n}\bigg|\\
&\le \bigg(\sum_{j=1}^{n-d} |\alpha_j|^2\bigg)^\frac12 
\bigg(\sum_{j=1}^{n-d} \jb{w_j ,  \d t_h[t_k]}_{\R^n}^2\bigg)^\frac12 \\
&\le \big|\d t_h[t_k]\big|\\
& \le \sup_{h = 1, \dots, d} \norm{\nabla t_h(x)},
\end{align*}

\noi
where we used Cauchy-Schwarz inequality once again in the last step.
This proves \eqref{CV4}
and hence concludes the proof of Lemma \ref{LEM:meas}.
\end{proof}

We now present the proof of Lemma \ref{LEM:surface}.

\begin{proof}[Proof of Lemma \ref{LEM:surface}]
Let $E_N$ 
and $U_\eps(\delta_*,\delta^*)$
as in \eqref{EN1} and \eqref{Ueps2}. By Proposition \ref{PROP:FOCN}, 
given any 
  $u\in  U_\eps(\delta_*,\delta^*)\cap E_N$, 
  there exist coordinates $(\dl,x_0,\ta) \in (0, \dl^*) \times \R \times (\R/2\pi \Z)$ and $v\in  V_{\dl,x_0, \ta} \cap E_N$ such that
\[u= e^{i\ta} \pi_N Q_{\dl, x_0}^\rho+v.\]

Given $(\dl, x_0,\ta)$, 
 we let $v_j=v_j(\dl, x_0,\ta)$, $j = 1, \dots, 4N-1$, denote an $H^1(\T)$-orthonormal basis\footnote{Once again, recall that we view $H^k(\T)$
as a Hilbert space over reals.} of $ V_{\dl, x_0,\ta}\cap E_N$.
Then, from 
\eqref{Q0} 
with $g = \{g_n\}_{ |n|\le N}$, we have 
\begin{align}
& \int_{U_\eps(\delta_*,\delta^*)}   F(\pi_N u)\mu(\ds  u)\notag \\
&\hphantom{X}
= \frac{1}{(2\pi)^{2N+1}}\int_{\C^{2N}} 
\ind_{{U_\eps(\delta_*,\delta^*)}}
\bigg(\sum_{ |n|\leq N} \frac{g_n}{\jb{n}} e^{2\pi i  nx}\bigg)
F\bigg(\sum_{ |n|\leq N} \frac{g_n}{\jb{n}} e^{2\pi i  nx}\bigg)e^{-\frac {|g|^2}2}\d  g\notag \\
\intertext{From Lemma \ref{LEM:meas} with $y = \{y_j\}_{j = 1}^{4N-1} \in \R^{4N-1} $, }
&\hphantom{X}
\les\frac{1}{(2\pi)^{(4N-1)/2}}\iint
\ind_{U_\eps(\delta_*,\delta^*)}\bigg(e^{i\ta}\pi_N Q_{\dl, x_0}^\rho+\sum_{j=1}^{4N-1}y_jv_j\bigg)\notag \\
&\hphantom{XXXX}
 \times 
F\bigg(e^{i\ta}\pi_N Q_{\dl, x_0}^\rho+\sum_{j=1}^{4N-1}y_jv_j\bigg) \notag \\
&\hphantom{XXXX}
 \times 
e^{-\frac{1}{2}\|e^{i\ta}\pi_N Q_{\dl, x_0}^\rho+\sum_{j=1}^{4N-1}y_jv_j\|^2_{{H}^1}}\,\d  y
\d \s_N(\dl, x_0, \ta)\notag \\
&\hphantom{X}
=\iint
\ind_{U_\eps(\delta_*,\delta^*)}\Big(e^{i\ta}\pi_N (Q_{\dl, x_0}^\rho+v)\Big)
 F\Big(e^{i\ta}\pi_N (Q_{\dl, x_0}^\rho+v)\Big)
\notag \\
&\hphantom{XXXX}
 \times 
e^{-\frac{1}{2}\|\pi_N Q_{\dl, x_0}^\rho\|^2_{ {H}^1(\T)}-\jb{ (1-\dx^2)\pi_N  Q_{\dl,x_0}^\rho,v}_{L^2(\T)}} 
\mu_{\dl, x_0}^\perp(\ds v)\d \sigma_N( \dl, x_0,\ta),
\label{eqn: surface}
\end{align}

\noi
where $\mu_{\dl, x_0}^\perp$ denotes  the Gaussian measure with 
  $V_{\dl,x_0,0} \subset H^1(\T)$
  as its 
 Cameron-Martin space.\footnote{In the last step
 of \eqref{eqn: surface}, 
 we used the decomposition $\mu^\perp_{\dl, x_0} = 
 \mu^\perp_{\dl, x_0, \leq N}\otimes  \mu^\perp_{\dl, x_0, > N}$, 
 where 
$  \mu^\perp_{\dl, x_0, \leq N}$ (and $\mu^\perp_{\dl, x_0, > N}$, respectively)
denotes  the Gaussian measure with 
  $V_{\dl,x_0,0} \cap E_N$ (and $V_{\dl,x_0,0} \cap \pi_N^\perp H^1(\T)$, respectively)
    as its 
 Cameron-Martin space.
 }
Here, the  measure
$\s_N$ is given by 
\[ d\s_N  (\dl, x_0,\ta) = 
\Big(1 + \sup_{k= 1, 2, 3} \norm{\nabla t_k(\dl, x_0,\ta)}^3\Big) \d \omega_N(\dl, x_0,\ta), \]

\noi
where $t_k = t_k(N, \dl, x_0, \ta)$, $k = 1, 2, 3$,  are the orthonormal vectors obtained by applying 
the Gram-Schmidt orthonormalization procedure in $H^1(\T)$ to 
 $\big\{\pi_N\partial_{\dl}(e^{i\ta}Q_{\dl, x_0}^\rho)$, 
 $\pi_N\partial_{x_0}(e^{i\ta}Q_{\dl, x_0}^\rho)$, 
 $\pi_N\partial_{\ta}(e^{i\ta}Q_{\dl, x_0}^\rho)\big\}$
  and 
the surface measure $\omega_N$ is given by 
\[ \d \omega_N (\dl, x_0,\ta)=|\gamma_N(\dl, x_0,\ta)|\ds \dl \ds  x_0\d\ta_0\]

\noi
with $\gamma_N(\dl, x_0,\ta)$ given by 
\begin{align}
& \gamma_N(\dl, x_0,\ta) \notag\\
&
=\det 
\begin{pmatrix}
\langle \pi_N\partial_{\dl}(e^{i\ta}Q_{\dl, x_0}^\rho),t_1\rangle_{H^1(\T)}
&  \langle \pi_N\partial_{x_0}(e^{i\ta}Q_{\dl, x_0}^\rho),t_1\rangle_{H^1(\T)} 
& \langle \pi_N\partial_{\ta}(e^{i\ta}Q_{\dl, x_0}^\rho),t_1\rangle_{H^1(\T)}\\
\langle \pi_N\partial_{\dl}(e^{i\ta}Q_{\dl, x_0}^\rho),t_2\rangle_{H^1(\T)}
&  \langle \pi_N\partial_{x_0}(e^{i\ta}Q_{\dl, x_0}^\rho),t_2\rangle_{H^1(\T)} 
& \langle \pi_N\partial_{\ta}(e^{i\ta}Q_{\dl, x_0}^\rho),t_2\rangle_{H^1(\T)}\\
\langle \pi_N\partial_{\dl}(e^{i\ta}Q_{\dl, x_0}^\rho),t_3\rangle_{H^1(\T)}
&  \langle \pi_N\partial_{x_0}(e^{i\ta}Q_{\dl, x_0}^\rho),t_3\rangle_{H^1(\T)} 
& \langle \pi_N\partial_{\ta}(e^{i\ta}Q_{\dl, x_0}^\rho),t_3\rangle_{H^1(\T)}
\end{pmatrix}\notag \\
&
=\det 
\begin{pmatrix}
\langle \partial_{\dl}(e^{i\ta}Q_{\dl, x_0}^\rho),t_1\rangle_{H^1(\T)}
&  \langle \partial_{x_0}(e^{i\ta}Q_{\dl, x_0}^\rho),t_1\rangle_{H^1(\T)} 
& \partial_{\ta}e^{i\ta}\langle Q_{\dl, x_0}^\rho,t_1\rangle_{H^1(\T)}\\
\langle \partial_{\dl}(e^{i\ta}Q_{\dl, x_0}^\rho),t_2\rangle_{H^1(\T)}
&  \langle \partial_{x_0}(e^{i\ta}Q_{\dl, x_0}^\rho),t_2\rangle_{H^1(\T)} 
& \partial_{\ta}e^{i\ta}\langle Q_{\dl, x_0}^\rho,t_2\rangle_{H^1(\T)}\\
\langle \partial_{\dl}(e^{i\ta}Q_{\dl, x_0}^\rho),t_3\rangle_{H^1(\T)}
&  \langle \partial_{x_0}(e^{i\ta}Q_{\dl, x_0}^\rho),t_3\rangle_{H^1(\T)} 
& \partial_{\ta}e^{i\ta}\langle Q_{\dl, x_0}^\rho,t_3\rangle_{H^1(\T)}
\end{pmatrix}.
\label{CV6}
\end{align}


Note that 
 $\{t_k\}_{k = 1}^3$ are chosen so that 
 $t_k(\dl, x_0, \ta)
 = \tau_{x_0}t_k(\dl, 0, \ta)$, 
 where $\tau_{x_0}$ is the translation defined in \eqref{trans1}.
Together with invariance under multiplication by a unitary complex number, 
we obtain
\begin{align}
  \gamma_N(\dl, x_0,\ta)=\gamma_N(\dl, 0, 0)
\label{CV7}
\end{align}

\noi
for any $x_0 \in \T$ and $\ta \in \R/(2\pi\Z)$.
%
%
%
From \eqref{CV7} and \eqref{CV6}
with \eqref{JA2}, \eqref{JA4}, and \eqref{JA7}, we have
\begin{align*}
\begin{split}
|\gamma_N(\dl,x_0, \ta)|
& = |\gamma_N(\dl, 0, 0)|
 \les\| \partial_\dl Q_\dl^\rho\|_{H^1(\T)}^3+\|\partial_{x_0}  Q_\dl^\rho\|_{H^1(\T)}^3
+ \|  Q_\dl^\rho\|_{H^1(\T)}^3\\
& \les \dl^{-6}.
\end{split}
\end{align*}

\noi
A computation analogous
to \eqref{JA2}, \eqref{JA4}, and \eqref{JA7}
together with \eqref{JC3} and \eqref{JC4} shows
\begin{align*}
\begin{split}
 \sup_{k = 1, 2, 3} \norm{\nabla t_k(\dl, x_0,\ta)}
& \les 
\sup_{\kk_1, \kk_2 \in \{\dl, x_0,\ta\}}
\| \partial^2_{\kk_1, \kk_2}(\pi_N  e^{i \ta}Q_{\dl, x_0}^\rho)\|_{H^1}\\
&  \les \delta^{-3}, 
\end{split}
\end{align*}

\noi
uniformly in $N \in \NB$.

Therefore, by taking the limit $N\to \infty$ in \eqref{eqn: surface}, 
the dominated convergence theorem yields
\begin{align*}
&  \int_{U_\eps(\delta_*,\delta^*)} F(u) \d \mu (u) \\
&\phantom{X}  \le
\iiiint_{U_\eps(\delta_*,\delta^*)}
 F\Big(e^{i\ta} (Q_{\dl, x_0}^\rho+v)\Big)
\\
&\phantom{XXXXXX} \times
e^{-\frac{1}{2}\| Q_{\dl, x_0}^\rho\|^2_{ {H}^1(\T)}-\jb{ (1-\dx^2) Q_{\dl,x_0}^\rho,v}_{L^2(\T)}} 
\mu_{\dl, x_0}^\perp(\ds v)\d \sigma( \dl)\d x_0\d \ta.
\end{align*}

\noi
Finally,
in view of Proposition \ref{PROP:FOCN}, 
by taking $\dl_*\to 0$, we obtain
 \eqref{LE0}.
\end{proof}

\subsection{Further reductions}
\label{SUBSEC:X1}

The forthcoming calculations aim to prove 
\eqref{eqn: to-bound}.
We first apply 
the change of variables (Lemma 
\ref{LEM:surface}) to the integral in \eqref{eqn: to-bound}.
In the remaining part of this section, 
we set $K = \|Q\|_{L^2(\R)}$
and use
$\jb{\,\cdot, \cdot\,}$ to denote the  inner product in $L^2(\T)$, 
unless otherwise specified. 
Note that the integrand
\begin{equation*}
\exp\bigg(\frac{1}{6}\int_\T |u|^6\ds  x\bigg)
\ind_{\{\|u\|_{L^2(\T)}\le K\}}
\end{equation*}

\noi
 is neither bounded nor continuous. 
 The lack of continuity is due to the sharp cutoff
$\ind_{\{\|u\|_{L^2(\T)}\le K\}}$.
Since the set of discontinuity has $\mu$-measure zero, 
we may start with a smooth cutoff and then pass to the limit.  
We can also replace this integrand  by a bounded one, as long as the bounds we obtain are uniform, at the cost of a simple approximation argument which we omit.

By applying 
Lemma \ref{LEM:surface} to the integral in \eqref{eqn: to-bound}
and  using also the translation  invariance of the surface measure
 (namely, the independence of all the quantities from $x_0$), 
we have
\begin{align}
\begin{split}
&\int_{U_\eps}\exp\bigg(\frac{1}{6}\int_\T |u|^6\,\ds  x\bigg)
\ind_{\{\|u\|_{L^2(\T)}\le K\}}\,\mu(\ds  u)\\
&\hphantom{X}
 \le\iiiint_{U_\eps} e^{F(v)}\ind_{\{\|Q_{\dl, x_0}^\rho+v\|_{L^2(\T)}\le K\}}
\,\mu_{\dl, x_0}^\perp(\ds v)d \sigma(\dl)\d x_0\d\ta\\
&\hphantom{X}
=2\pi\iint_{U_\eps}e^{F(v)}\ind_{\{\| Q_{\dl}^\rho+v\|_{L^2(\T)}\le K\}}
\,\mu_{\dl}^\perp(\ds v)\d \sigma(\dl),
\end{split}
\label{XX1}
\end{align}

\noi
where  $\mu_{\dl}^\perp = \,\mu_{\dl, 0}^\perp$ and 
\begin{align}
F(v)=\frac{1}{6}\int_{\T} | Q_{\dl}^\rho+v|^6\,\ds  x
-\frac{1}{2}\| Q_{\dl}^\rho\|^2_{{H}^1(\T)}-\jb{(1-\dx^2)  Q_{\dl}^\rho,v}.
\label{XX2}
\end{align}

A direct computation shows\footnote{In view of \eqref{inner1}, we did not need to use the real part
symbol in \eqref{XX3}.  We, however, chose to use the real part symbol
for clarity.}
\begin{equation}
\begin{split}
\frac{1}{6} & \int_\T  | Q_\dl^\rho+v|^6\,\ds  x\\
&= \frac{1}{6}\int ( Q_\dl^\rho)^6\,\ds  x+\jb{( Q_\dl^\rho)^5, \Re v}
+\langle ( Q_\dl^\rho)^4,\Re (v^2)+\tfrac{3}{2}|v|^2\rangle\\
&\hphantom{X}
+ \langle ( Q_\dl^\rho)^3, \tfrac{1}{3}\Re (v^3)+3|v|^2\Re v\rangle
+\langle ( Q_\dl^\rho)^2, |v|^2(\Re (v^2)+\tfrac{3}{2}|v|^2)\rangle\\
&\hphantom{X}
+ \langle Q_\dl^\rho,|v|^4\Re v\rangle+\frac{1}{6}\int_\T |v|^6\,\ds  x.
\end{split}
\label{XX3}
\end{equation}

\noi
From \eqref{XX2} and \eqref{XX3} with  \eqref{G3} and \eqref{QQ1}, we have
\begin{align}
\begin{split}
F(v) & = \frac{1}{6}\int_\T ( Q_\dl^\rho)^6\, \ds  x-\frac{1}{2}\| Q_\dl^\rho\|_{H^1(\T)}^2\\
&\hphantom{X}
 + \langle (\rho^5 - \rho) Q_\dl^5,  \Re v \rangle 
 + 2 \langle (\dd_x \rho) (\dd_x Q_\dl), \Re v \rangle + \langle (\dd_x^2 \rho) Q_\dl, \Re v\rangle \\
&\hphantom{X}
+(2\dl^{-2} -1)\langle  Q_\dl^\rho, \Re v\rangle
+ \langle ( Q_\dl^\rho)^4,2(\Re v)^2+\tfrac{1}{2}|v|^2\rangle\\
&\hphantom{X}
+ \langle ( Q_\dl^\rho)^3, \tfrac{1}{3}\Re (v^3)+3|v|^2\Re v\rangle
+\langle ( Q_\dl^\rho)^2, |v|^2(\Re (v^2)+\tfrac{3}{2}|v|^2)\rangle\\
&\hphantom{X}
+ \langle Q_\dl^\rho,|v|^4\Re v\rangle+\frac{1}{6}\int_\T |v|^6\,\ds  x,
\end{split}
\label{XX3a}
\end{align}

\noi
where we used $\Re (v^2) = 2(\Re v)^2 - |v|^2$
to get the seventh term on the right-hand side.

By the sharp Gagliardo-Nirenberg-Sobolev inequality
(Proposition \ref{THM:W}), 
we have $ \|Q_\dl\|_{L^6(\R)}^6 = 3 \|Q'_\dl\|_{L^2(\R)}^2$
and thus 
 we have
\begin{align}
\bigg|\frac{1}{6}\int_\T (Q_\dl^\rho)^6\, \ds  x-\frac{1}{2}
\| Q_\dl^\rho\|_{H^1(\T)}^2\bigg| \le C, 
\label{XX4}
\end{align}
uniformly in $0< \dl\le 1$, 
thanks to  the exponential decay of $Q$ (as in \eqref{corr1}). 
Moreover,  recalling from \eqref{rhodef} that $\rho \equiv 1$ on $[-\frac 18,\frac 18]$
and using  the exponential decay of $Q$ again with  $\| v\|_{L^2(\T)} \le 1$,  we obtain 
\begin{equation}
 \Big| \langle (\rho^5 - \rho) Q_\dl^5,  \Re v \rangle 
 + 2 \langle (\dd_x \rho) (\dd_x Q_\dl), \Re v \rangle + \langle (\dd_x^2 \rho) Q_\dl, \Re v\rangle 
 \Big| \le C,
\label{XX4b}
\end{equation}
uniformly in $0 \le \dl \le 1$.
By Young's inequality, 
the sum of the last four terms in \eqref{XX3a} is bounded by 
\begin{align}
\eta \int_\T ( Q_\dl^\rho)^4\big(2(\Re v)^2 + \tfrac{1}{2}|v|^2\big)\,\ds  x + C_\eta \int_\T |v|^6\,\ds  x
\label{XX5}
\end{align}
for any $0<\eta\ll1$  and a (large) constant $C_\eta> 0$.
Hence, from 
\eqref{XX1}, \eqref{XX3a}, \eqref{XX4}, \eqref{XX4b}, and \eqref{XX5}
together with Proposition \ref{PROP:FOC}, 
we have
\begin{align}
\begin{split}
\int_{U_\eps} &  \exp\left(\frac{1}{6}\int_{\T}|u|^6\,\ds  x\right)
\ind_{\{\|u\|_{L^2(\T)}\le K\}}\,\mu(\ds  u)\\
\les & \int_0^{\dl^*}  \int_{\{\|v\|_{L^2(\T)}\le \eps_1\}}e^{ G(v)}
e^{C_\eta \int |v|^6\,\ds  x}
\, 
\ind_{\{\| Q_{\dl}^\rho+v\|_{L^2(\T)}\le K\}}\,\mu_\dl^{\perp}(\ds v)\sigma(\ds \dl),
\end{split}
\label{XX5a}
\end{align}

\noi
where $G(v)$ is defined by 
\begin{align}G(v)=(2\dl^{-2}-1)\langle  Q_\dl^\rho,\Re v\rangle
+(1+\eta)
\langle (Q_\dl^\rho)^4,2(\Re v)^2+\tfrac{1}{2}|v|^2\rangle.
\label{XX6}
\end{align}

\noi
Here, $\eps_1 > 0$ is a small number to be chosen later (see Lemma \ref{lem: v6-int} below), 
which also appears in Proposition \ref{PROP:FOC},
determining
small $\eps = \eps(\eps_1) > 0$ and 
$ \dl^*=  \dl^*(\eps_1)  > 0$.

By a slight modification  of the argument presented in Subsection \ref{SUBSEC:1d}, we have the following
integrability result.

\begin{lemma}\label{lem: v6-int}
Given any  $C'_\eta>0$, there exists small $\eps_1>0$  such that
\begin{equation}
\int_{\{\|v\|_{L^2(\T)} \le \eps_1\}} \exp\bigg(C'_\eta \int_{\T}|v|^6\,\ds  x\bigg) \mu^\perp_\dl (\d v)<\infty, 
\label{XX7}
\end{equation}

\noi
uniformly in $0 < \dl \ll1$.
\end{lemma}
\begin{proof}
Let $W$ be a finite-dimensional subspace of $H^1(\T)$ of dimension $n$
with an orthonormal basis $\{w_1,\dots,w_n\} \subset H^2(\T)$ (with respect to the $H^1(\T)$-inner product). Define the projector $P_{W^\perp}$ by 
\begin{equation}
P_{W^\perp} u = u -\sum_{j=1}^n \jb{u,(1- \dx^2) w_j} w_j, 
\label{XX7a}
\end{equation}

\noi
where $\jb{\cdot, \cdot} = \jb{\cdot, \cdot}_{L^2(\T)}$.
On $H^1(\T)$, this is nothing but the usual $H^1$-orthogonal projection 
onto $W^\perp$.
The definition~\eqref{XX7a} allows us to extend 
$P_{W^\perp}$ to $L^2(\T)$.
Then, 
by the definition of $\mu^\perp_\dl = \mu_{\dl, 0}^\perp$, it suffices  to show that given $C'_\eta > 0$, there exist $M \gg 1$ and $\eps_1>0$, depending only on $n = \dim W$,  such that 
\begin{equation}
\int_{\{\|P_{W^\perp}u\|_{L^2(\T)} \le \eps_1\}} \exp\bigg(C'_\eta \int_{\T}|P_{W^\perp}u|^6\,\ds  x\bigg) \mu (\d u) \le M. 
\label{XX7w}
\end{equation}

\noi
Indeed, in view of \eqref{eqn: Vdelta-def}
and the definition of $\mu^\perp_\dl = \mu_{\dl, 0}^\perp$
with the Cameron-Martin space $V_{\dl, 0, 0}$, 
by  simply setting 
$w_1 = \frac{ \dd_\dl Q_\dl^\rho}{\| \dd_\dl Q_\dl^\rho\|_{H^1(\T)}}$, 
$w_2 = \frac{\dd_{x_0} Q_\dl^\rho}{\|\dd_{x_0}  Q_\dl^\rho\|_{H^1(\T)}}$, 
and $w_3 = \frac{i  Q_\dl^\rho}{\| Q_\dl^\rho\|_{H^1(\T)}}$ with $n = 3$, 
the desired bound~\eqref{XX7} follows from \eqref{XX7w}.

The proof of the inequality \eqref{XX7w} follows closely the argument presented in Subsection~\ref{SUBSEC:1d}. 
In the following, we point out the modifications required to obtain \eqref{XX7w}. First of all, we replace the definition \eqref{H2ab} of the set $E_k$ by 
\begin{align}
\begin{split}
E_k&=\big\{\|(P_{W^\perp}u)_{\geq 0} \|_{L^p}>\ld, \ldots, \|(P_{W^\perp}u)_{\ge k-1}\|_{L^p}>\ld, \|(P_{W^\perp}u)_{\ge k}\|_{L^p}\le \ld \big\} \\
&\subset \big\{\|(P_{W^\perp}u)_{\ge k-1}\|_{L^p}>\ld, \|(P_{W^\perp}u)_{\ge k}\|_{L^p}\le \ld \big\}.
\end{split}
\label{XX7w1}
\end{align}
Namely, in the definition \eqref{H2ab} of $E_k$, we replace $u$ with $P_{W^\perp}u$. 
Arguing as in \eqref{eqn: livius} with~\eqref{XX7w1}, 
it suffices to show that
\begin{equation}\label{XX7w2}
\begin{split}
\E\Big[  e^{C'_\eta \int_\T |P_{W^\perp} u(x)|^6\,\ds  x}, \,
& \|(P_{W^\perp} u)_{\ge k-1}\|_{L^p(\T)}>\ld, \\
& \|(P_{W^\perp} u)_{\ge k}\|_{L^p(\T)}\le \ld, \,\|P_{W^\perp} u \|_{L^2(\T)}\le \eps_1 \Big]
\end{split}
\end{equation}

\noi
is summable in $k \in \NB$.
Proceeding in the same way as in \eqref{eqn: ezekiel} and  \eqref{H2a}, 
we obtain the following analogue of \eqref{eqn: titus}, bounding \eqref{XX7w2}:
\begin{equation}\label{eqn: titusw}
\begin{split}
e^{C_6(  \eps) \ld^6 }
 \E\Big[& e^{C'_\eta(1+\dl_0)\int_\T |(P_{W^\perp}u)_{\le k-1}(x)|^6\, \ds  x},\\
 &  \|(P_{W^\perp}u)_{\ge k-1}\|_{L^p(\T)}>\ld ,\|P_{W^\perp}u\|_{L^2(\T)}\le \eps_1\Big], 
\end{split}
\end{equation}

\noi
where we used $\dl_0$ for the constant $\dl = \dl(6, \eps)$ in \eqref{H2a}
to avoid confusion with the dilation parameter~$\dl$. 
Note that from \eqref{GNS1} in Lemma \ref{LEM:periodicGN} with \eqref{XX7a}, we have 
\begin{align}
\begin{split}
\| ( & P_{W^\perp}u)_{\le k-1}  \|_{L^6(\T)}^6 \les \| (P_{W^\perp}u)_{\le k-1}\|_{H^1(\T)}^{2}\|P_{W^\perp}u\|_{L^2(\T)}^{4}  \\
& \les \Big(\|u_{\le k-1}\|_{H^1(\T)}^2 + \sum_{j=1}^n |\jb{u,(1-\dx^2) w_j}|^2 
\|(w_j)_{\le k-1}\|_{H^1(\T)}^2\Big)\|P_{W^\perp}u\|_{L^2(\T)}^{4} \\
&\le \Big(\|u_{\le k-1}\|_{H^1(\T)}^2 + \sum_{j=1}^n|\jb{u,(1-\dx^2) w_j}|^2\Big)\|P_{W^\perp}u\|_{L^2(\T)}^{4},
\end{split}
\label{XX7w3}
\end{align}
where the implicit constant depends only on $n = \dim W$.

Using the inequality \eqref{XX7w3} instead of~\eqref{GNS1}, and fixing $\ld = 1$, we obtain 
\begin{align*}
\eqref{eqn: titusw} \le  e^{C_6(\eps)} \E\bigg[
& \exp\bigg(CC'_\eta(1+\dl_0) \eps_1^4\Big(\|u_{\le k-1}\|_{H^1(\T)}^2 
+ \sum_{j=1}^n|\jb{u,(1-\dx^2) w_j}|^2\Big)\bigg), \\
& \|(P_{W^\perp}u)_{\ge k-1}\|_{L^p(\T)}>1\bigg]
\end{align*}

\noi
for some $C = C(n) > 0$.
Then, 
by H\"older's inequality, this implies the following analogue of \eqref{eqn: haddock}:
\begin{equation}
\begin{aligned}
\eqref{eqn: titusw} &\le e^{C_6(\eps)}\bigg\{\E\Big[\exp\Big((n+2)CC'_\eta(1+\dl_0) \eps_1^4\|u_{\le k-1}\|_{H^1(\T)}^2\Big)
\Big]\bigg\}^{\frac{1}{n+2}}\\
& \hphantom{\le}
\times \prod_{j=1}^n\bigg\{\E\Big[\exp\Big((n+2)CC'_\eta (1+\dl_0)\eps_1^4|\jb{u,(1-\dx^2) w_j}|^2\Big)
\Big]\bigg\}^{\frac{1}{n+2}} \\
& \hphantom{\le}
\times \Big\{\P\big(\|(P_{W^\perp}u)_{\ge k-1}\|_{L^6(\T)}> 1\big)\Big\}^{\frac{1}{n+2}}.
\end{aligned}  
\label{eqn: haddockw}
\end{equation}

We note 
from \eqref{Q0}
 that $\jb{u, (1-\dx^2) w_j}$ is a mean-zero Gaussian random variable with 
 $\E\big[|\jb{u,(1-\dx^2) w_j}|^2\big] = \| w_j \|_{H^1(\T)}^2 = 1$. In particular, 
 if $(n+2)CC'_\eta (1+\dl_0) \eps_1^4 < \frac 12$, 
it follows from~\eqref{H3a} that 
\begin{equation} \label{H4w}
\E\Big[\exp\Big((n+2)CC'_\eta (1+\dl_0) \eps_1^4|\jb{u,(1-\dx^2)w}|^2\Big)\Big] < \infty.
\end{equation}

\noi
Moreover, by Bernstein's inequality with  $\| w_j \|_{H^1(\T)}^2 = 1$, there exists $c = c(n) > 0$ such that
\begin{align} 
\begin{split}
\P\bigg(\|&  (\jb{u,(1-\dx^2)w_j} w)_{\ge k-1} \|_{L^6(\T)} > \frac1{n+1}\bigg) \\
&= \P\bigg(|\jb{u,(1-\dx^2)w_j}| > \frac{1}{(n+1) \| (w_j)_{\ge k-1} \|_{L^6(\T)}}\bigg) \\
&\le \P\bigg(|\jb{u,(1-\dx^2)w_j}| \gtrsim \frac{2^{\frac 23 k}}{\| w_j \|_{H^1(\T)}}\bigg) \\
&\les \exp(-c2^{\frac 43 k}), 
\end{split}
\label{XX7w4}
\end{align}

\noi
uniformly in $j = 1, \dots, n$.
Hence, 
from \eqref{XX7a}, 
 \eqref{eqn: vanguard}, and \eqref{XX7w4}, we obtain, for some $c' = c'(n) > 0$,
\begin{align}
\P& \big(\|(P_{W^\perp}u)_{\ge k-1}\|_{L^6(\T)}> 1\big) \notag\\
&\le \P\Big(\|u_{\ge k-1}\|_{L^6(\T)}> \frac 1{n+1}\Big) 
+ \sum_{j=1}^n \P\Big(\| (\jb{u,(1-\dx^2)w_j} w_j)_{\ge k-1} \|_{L^6(\T)} > \frac 1{n+1}\Big) \notag \\
&\les \exp(-c'2^{\frac 43 k})
 \label{H3w},
\end{align}

\noi
which is an analogue of \eqref{H3}. 
In addition, if $(n+2)CC'_\eta (1+\dl_0) \eps_1^4 < \frac 12$, 
then using \eqref{H3a}, we can repeat the computations in \eqref{H4} and obtain 
\begin{equation} \label{H5w}
\begin{split}
\E&  \Big[\exp\Big((n+2)CC'_\eta   (1+\dl_0) \eps_1^4\|u_{\le k-1}\|_{H^1(\T)}^2\Big)\Big] \\
& \le (1-2(n+2)CC'_\eta(1+\dl_0)\eps_1^4)^{-2^k}.
\end{split}
\end{equation}

Hence, by applying  \eqref{H5w}, \eqref{H4w}, and \eqref{H3w} to \eqref{eqn: haddockw}
and proceeding as in~\eqref{H4a}, we see that \eqref{XX7w2}
is summable in $k \in \NB$.
Therefore, we conclude that 
\begin{equation*}
\int_{\{\|P_{W^\perp}u\|_{L^2(\T)} \le \eps_1\}} \exp\bigg(C'_\eta \int_{\T}|P_{W^\perp}u|^6\,\ds  x\bigg) \mu (\d u) \les 1, 
\end{equation*}
where the implicit constant depends only on the dimension $n$ of $W$.
\end{proof}

%

By applying 
H\"older's inequality and Lemma \ref{lem: v6-int}
to
\eqref{XX5a}, 
we obtain
\begin{align*}
\begin{split}
& \int_{U_\eps}   \exp\left(\frac{1}{6}\int_{\T}|u|^6\,\ds  x\right)
\ind_{\{\|u\|_{L^2(\T)}\le K\}}\,\mu(\ds  u)\\
& \hphantom{X}
\les 
\int_0^{\dl^*}  \bigg(\int_{\{\|v\|_{L^2(\T)}\le \eps_1\}}
e^{(1+\eta)G(v)}
\ind_{\{\| Q_\dl^\rho+v\|_{L^2(\T)}\le K\}} \mu_\dl^\perp(\d v)\bigg)^{\frac{1}{1+\eta}}\sigma(\ds \dl).
\end{split}
\end{align*}

\noi
Our goal in the remaining part of this section is to bound
the inner integral on the right-hand side.

Next, we decompose the subspace\footnote{Recall the convention
$ V_{\dl,0,0} = V_{\dl,0,0}\cap H^1(\T)$ introduced
at the beginning of Subsection \ref{SUBSEC:change}.}
$V_{\dl,0,0} = V_{\dl,0,0}(\T) \subset H^1(\T)$ 
in \eqref{eqn: Vdelta-def}
as 
\begin{equation}\label{eqn: decomp}
V_{\dl,0,0}= V'
\oplus 
\mathrm{span}\{e\},
\end{equation}
where $\|e\|_{H^1(\T)}=1$ and $e$ is orthogonal in $H^1(\T)$ to 
\begin{align}
V'= \big\{w\in V_{\dl,0,0}: \langle w,  Q_\dl^\rho\rangle=0\big\}.
\label{XX8}
\end{align}

Denote by  $P^{H^1}_W$  the $H^1$-orthogonal projection onto 
a given  subspace $W \subset H^1(\T)$.
Then, by noting the orthogonality of $ Q_\dl^\rho$ with $\dd_{x_0}  Q_\dl^\rho$ and $i  Q_\dl^\rho$ in $L^2(\T)$
and by directly computing $\jb{ Q_\dl^\rho,  \dd_\dl Q_\dl^\rho}$ with \eqref{Q1} and integrating by parts, 
we  have
\begin{equation}
\label{eqn: e}
\begin{split}
e & = \frac{P^{H^1}_{V_{\delta,0,0}} (1-\dx^2)^{-1} Q_\delta^\rho}{\norm{P^{H^1}_{V_{\delta,0,0}} (1-\dx^2)^{-1}  Q_\delta^\rho}_{H^1(\T)}} \\
& = \frac{(1-\dx^2)^{-1} Q_\delta^\rho}
{\norm{(1-\dx^2)^{-1}  Q_\delta^\rho}_{H^1(\T)}} 
+ O_{L^2(\T)}(\exp(-c\delta^{-1})) 
\end{split}
\end{equation}

\noi
for $0 < \dl \ll1 $.
Corresponding to the decomposition \eqref{eqn: decomp} of the Cameron-Martin space
for $\mu_\dl^\perp$, we have the following decomposition of the measure:
\begin{equation*} 
\begin{split}
\ds \mu^\perp_\dl (v) &= \frac{1}{\sqrt{2\pi}}e^{-\frac{1}{2}g^2}\ds g\ds \mu^{\perp\perp}_\dl(w)
\end{split}
\end{equation*}

\noi
with 
\begin{align}
v= ge+w, \quad w \in V'.
\label{XX9}
\end{align}

\begin{lemma}\label{LEM:G1}
Let $G(v)$ be as in \eqref{XX6}. Then, we have 
\begin{align}
\begin{split}
& \int_{\{\|v\|_{L^2(\T)}\le \eps_1\}}
e^{(1+\eta)G(v)}
\ind_{\{\| Q_\dl^\rho+v\|_{L^2(\T)}\le K\}} \mu_\dl^\perp(\d v)\\
& \quad \lesssim \int \exp\left(-(1-\eta^2)\wt{H}_\dl(w)\right) \mu_\dl^{\perp\perp}(\d w), \label{XG1}
\end{split}
\end{align}

\noi
uniformly in $0 < \dl \ll1$, 
where $\wt{H}_\dl(w)$ is given by 
\begin{align*}
\wt{H}_\dl(w)= \dl^{-2}\int_{\T} |w|^2\,\ds  x- \left(1+5\eta\right)
\int_{\T} (Q_\dl^\rho)^4
\big(2(\Re w)^2+\tfrac{1}{2}|w|^2\big)
\ds  x.
\end{align*}

\end{lemma}

\begin{proof}

By expanding $\langle  Q_\dl^\rho+v,  Q_\dl^\rho+v\rangle$
with the decomposition \eqref{XX9}, we have
\begin{align*}
\begin{split}
\|Q\|_{L^2(\R)}^2=K^2&\ge \langle  Q_\dl^\rho+v,  Q_\dl^\rho+v\rangle\\
&= \| Q_\dl^\rho\|^2_{L^2(\T)}+2g\langle  Q_\dl^\rho, e\rangle + \|v\|^2_{L^2(\T)}\\
&= \| Q_\dl^\rho\|^2_{L^2(\T)}+2g\langle  Q_\dl^\rho, e\rangle+\|ge\|^2_{L^2(\T)}+\|w\|^2_{L^2(\T)}
+ 2 g\jb{e, w}.
\end{split}
\end{align*}

\noi
Together with \eqref{corr1}, we obtain
\begin{align}
\begin{split}
 2g & \langle    Q_\dl^\rho, e\rangle + \|v\|^2_{L^2(\T)}\\
&= 
2g\langle  Q_\dl^\rho, e\rangle+\|ge\|^2_{L^2(\T)}+\|w\|^2_{L^2(\T)}
+ 2 g\jb{e, w}\\
& \leq C \exp(-c\dl^{-1}).
\end{split}
\label{XG2}
\end{align}

\noi
We also note from \eqref{eqn: e} that 
\begin{align}
\langle  Q_\dl^\rho, e\rangle = \| Q_\dl^\rho\|_{H^{-1}(\T)} + 
O(\exp(-c\delta^{-1}))
 \sim \|Q_\dl\|_{H^{-1}(\R)} \sim \dl \sim \|e\|_{L^2(\T)}, 
\label{XG2a}
\end{align}

\noi
uniformly in $0 < \dl \ll 1$.

We also note that if $\|v\|_{L^2(\T)}\leq \eps \leq \eps_1$, 
then we have 
\begin{align}
\| ge \|_{L^2(\T)} + \|w\|_{L^2(\T)} \les \eps.
\label{XG2b}
\end{align}

\noi
Indeed, by observing
\begin{align*}
\Big|\|& Q_\dl^\rho\|_{L^2(\T)}^2  - \| Q_\dl^\rho + v \|_{L^2(\T)}^2\Big|\\
& = \Big|\| Q_\dl^\rho\|_{L^2(\T)}  - \| Q_\dl^\rho + v \|_{L^2(\T)}\Big|
\Big( \| Q_\dl^\rho\|_{L^2(\T)}  + \| Q_\dl^\rho + v \|_{L^2(\T)}\Big)\\
& \les \|v\|_{L^2(\T)} \Big( \| Q_\dl^\rho\|_{L^2(\T)}  + \|  v \|_{L^2(\T)}\Big)\\
& \les \eps, 
\end{align*}

\noi
we have
\begin{align}
\begin{split}
|2g\jb{ Q_\dl^\rho, e} |
& \leq \Big|2g\jb{ Q_\dl^\rho, e} +  \|v\|_{L^2(\T)}^2 \Big|  + \|v\|_{L^2(\T)}^2\\
&  = \Big|\| Q_\dl^\rho\|_{L^2(\T)}^2  - \|Q_\dl^\rho + v \|_{L^2(\T)}^2\Big| + O(\eps^2)\\
& \les \eps. 
\end{split}
\label{XG2c}
\end{align}

\noi
Then, from \eqref{XG2a} and \eqref{XG2c}, 
we conclude that 
$\| ge \|_{L^2(\T)}  \sim |g\jb{ Q_\dl^\rho, e} |\les \eps$.
By the triangle inequality with $w = v - ge$, 
we obtain 
$\|w\|_{L^2(\T)} \les \eps.$
This proves \eqref{XG2b}.

Now, from 
\eqref{XX6}, 
the decomposition \eqref{XX9}, 
Cauchy's inequality, and $\norm{Q_\delta^\rho}_{L^\infty(\T)} \sim \delta^{-\frac12}$, we have 
\begin{align}
\begin{split}
G(v)
& \le (2\dl^{-2}-1)\langle  Q_\dl^\rho,ge\rangle 
+ C_\eta \int_\T ( Q_\delta^\rho)^4 |ge|^2 \d x
 \\
 &\hphantom{XXXXX}
 + (1+2\eta)\int_{\T} ( Q_\dl^\rho)^4
\big(2(\Re w)^2+\tfrac{1}{2}|w|^2\big)\ds  x\\
& \le\dl^{-2}\big((2-\dl^2)\langle  Q_\dl^\rho,ge\rangle
 +  C'_\eta \norm{ge}_{L^2(\T)}^2\big)
 \\
 &\hphantom{XXXXX}
   + (1+2\eta)\int_{\T} ( Q_\dl^\rho)^4\big(2(\Re w)^2+\tfrac{1}{2}|w|^2\big)\ds x. 
\end{split}  
\label{eqn: G(v) - G(w)}
\end{align}

\noi
Noting that 
$(1+\eta) (1 + 2\eta) \le (1-\eta^2) (1+ 5\eta)$ for any $0 < \eta \ll 1$, in order to obtain \eqref{XG1}, 
we only need to bound the first term on the right-hand side of \eqref{eqn: G(v) - G(w)}.

\smallskip

\noi
$\bullet$ {\bf Case 1:}
$g  \geq 0$.
\\
\indent
In this case, from \eqref{XG2} and \eqref{XG2a} (which implies $\jb{ Q_\dl^\rho, e}>0$), we have 
\begin{align}
\|v\|^2_{L^2(\T)}  \le 
  g  \langle    Q_\dl^\rho ,  e\rangle + \|v\|^2_{L^2(\T)} 
  = 
O(\exp(-c\delta^{-1})).
\label{XG3}
\end{align}

\noi
Then, 
for sufficiently small $\dl > 0$, 
\eqref{XG3} shows that the hypothesis 
of \eqref{XG2b}
on the size of 
$\| v\|_{L^2(\T)}$
is satisfied 
with $\eps \sim 
\exp(-c\delta^{-1})$. 
Thus, from \eqref{XG2b}, 
we have 
%
%
\begin{align}
\| ge \|_{L^2(\T)} + \|w\|_{L^2(\T)}   = 
O(\exp(-c\delta^{-1})), 
\label{XG3a}
\end{align}

\noi
provided that $\dl = \dl(\eta, \eps_1) > 0$ is sufficiently small.
Hence, 
from 
 \eqref{XG3} and \eqref{XG3a} with \eqref{XX9}, 
we obtain, for some $C''_\eta > 0$,
\begin{align}
\begin{split}
(2-\dl^2)\langle Q_\dl^\rho,ge\rangle
 +  C'_\eta \norm{ge}_{L^2(\T)}^2
\leq 
- (1- \eta) \|w\|_{L^2(\T)}^2 + C''_\eta\exp(-c\delta^{-1}).
\end{split}
\label{XG4}
\end{align}

\noi
Therefore, 
from \eqref{eqn: G(v) - G(w)} and \eqref{XG4}, 
 the contribution to the left-hand side of \eqref{XG1} in this case
is bounded by 
\begin{align*}
\begin{split}
&  \lesssim 
\int \int_{
\{0 \leq g \les \exp(-c\delta^{-1})\}}\exp\left(-(1-\eta^2)\wt{H}_\dl(w)\right) e^{-\frac 12 g^2} dg \d \mu_\dl^{\perp\perp}( w)\\
&  \lesssim 
 \int  \exp\left(-(1-\eta^2)\wt{H}_\dl(w)\right) \mu_\dl^{\perp\perp}( \d w).
\end{split}
\end{align*}

\smallskip

\noi
$\bullet$ {\bf Case 2:}
$ g < 0$.\\
\indent
Under the condition $\| v\|_{L^2}\leq \eps_1$, 
it follows from \eqref{XG2}, 
Cauchy's inequality, \eqref{XG2a}
(which implies $\|ge\|_{L^2(\T)} \sim -g \jb{ Q_\dl^\rho, e}$), and \eqref{XG2b} that 
\begin{align*}
\begin{split}
2g\langle  Q_\dl^\rho, e\rangle
& \le -\big(1-\tfrac \eta 4\big)\|w\|_{L^2(\T)}^2
+ K_\eta\|ge\|^2_{L^2(\T)}
+ O(\exp(-c\dl^{-1}))\\
& \le -\big(1-\tfrac \eta 4\big)\|w\|_{L^2(\T)}^2
- \eps_1 K'_\eta g \jb{Q_\dl^\rho, e}
+ O(\exp(-c\dl^{-1})).
\end{split}
\end{align*}

\noi
By choosing $\eps_1 = \eps_1(\eta)> 0$ 
ands $\dl = \dl(\eta) > 0$ sufficiently small, 
we then have
\begin{align*}
\begin{split}
(2 - \dl^2) g\langle  Q_\dl^\rho, e\rangle
& \le -\big(1-\tfrac \eta 2\big)\|w\|_{L^2(\T)}^2
+ O(\exp(-c\dl^{-1})).
\end{split}
\end{align*}

\noi
Hence, we obtain
\begin{align*}
\begin{split}
(2-\dl^2) g\langle  Q_\dl^\rho, e\rangle
   +  C'_\eta \norm{ge}_{L^2(\T)}^2
& \le \big(2 -\dl^2 -  O(\eps_1)\big)  g\langle  Q_\dl^\rho, e\rangle
  \\
& \le -(1-\eta)\|w\|_{L^2(\T)}^2
+ O(\exp(-c\dl^{-1})).
\end{split}
\end{align*}

\noi
Proceeding as in Case 1, we also obtain \eqref{XG1} in this case.
\end{proof}

\subsection{Spectral analysis}
\label{SUBSEC:spec}

Given small $\dl, \eta > 0$, define an operator   $A = A(\dl, \eta)$ on $H^1(\T)$ by 
\begin{equation}\label{eqn: Adef}
Aw= P^{H^1}_{V'}(1-\dx^2)^{-1}\bigg(\dl^{-2}P^{H^1}_{V'}w-(1+5\eta)( Q_\dl^\rho)^4
\Big(2\Re (P^{H^1}_{V'}w)+\frac{1}{2}P^{H^1}_{V'}w\Big)\bigg), 
\end{equation}

\noi
where $V'$ is as in \eqref{XX8}.
Then, 
the integrand of the Gaussian integral 
on the right-hand side of~\eqref{XG1}  can be written as 
\begin{align}
\exp\left(-(1-\eta^2) \langle A w, w\rangle_{H^1(\T)}\right),
\label{AA1}
\end{align}

\noi
The main goal of this subsection is to establish the following integrability result.

\begin{proposition}\label{PROP: gauss-int}
Given any sufficiently small $\eta>0$, 
there exists a constant $c = c(\eta)>0$ such that
\begin{equation}\label{eqn: A-bound}
\int \exp\big(-(1-\eta^2)\langle Aw,w\rangle_{H^1(\T)}\big)\ds \mu^{\perp\perp}_\dl(w)\lesssim \exp(-c\dl^{-1}), 
\end{equation}

\noi
uniformly in $0 < \dl \ll1$.
\end{proposition}

From \eqref{eqn: Adef}
with 
$\norm{Q_\delta^\rho}_{L^\infty(\T)} \sim \delta^{-\frac12}$, 
we have 
\[| \langle Aw,w\rangle_{H^1(\T)}|\le C\dl^{-2}\|w\|_{L^2(\T)}^2.\]

\noi
Thus, by 
Rellich's lemma, we see that 
$A$ is a compact operator on $V' \subset H^1(\T)$
and thus the spectrum of $A$ consists of eigenvalues.
Recalling 
that 
$V' \subset H^1(\T)$ in~\eqref{XX8} is the 
 Cameron-Martin space
for $\mu_\dl^{\perp\perp}$, 
we see that 
evaluating the integral in \eqref{eqn: A-bound} is equivalent 
to estimating the product of  the eigenvalues of $\frac12\Id+(1-\eta^2)A$ on $V'$.
Thus, the rest of this subsection is devoted to studying
the eigenvalues of $A$.
We present the proof of Proposition~\ref{PROP: gauss-int} 
at the end of this subsection.

 Since $A$ preserves the subspace of $\Re H^1(\T)$ consisting of real-valued functions and 
 also the subspace $\Im H^1(\T)$ consisting of purely-imaginary-valued functions, 
in studying the spectrum of  $A$, 
 we split the analysis onto these two subspaces. 
 On the real subspace
$\Re H^1(\T)$, $A$ coincides with the operator
\begin{equation}\label{Areal}
A_1=P^{H^1}_{V_1(\dl)}(1-\dx^2)^{-1}\left(\dl^{-2}-\frac{5}{2}(1+5\eta)( Q_\dl^\rho)^4 \right)
P^{H^1}_{V_1(\dl)},
\end{equation}
where 
\begin{align}
\begin{split}
V_1 (\dl) = \big\{ w\in \Re H^1(\T): 
\, & \jb{w,(1-\dx^2) \partial_{\dl}Q_\delta^\rho} =0, \\
& \jb{w,(1-\dx^2)\partial_{x_0} Q_\delta^\rho}=0, \ \jb{w,  Q_\dl^\rho} = 0\big\}.
\end{split}
\label{VV1}
\end{align}

\noi
On the other hand, identifying $\Im H^1(\T)$ with $\Re H^1(\T)$ by the multiplication with $i$, 
we see that, on the imaginary subspace, $A$ coincides with the operator
\begin{equation}\label{Aimaginary}
A_2=P^{H^1}_{V_2(\dl)}(1-\dx^2)^{-1}\left(\dl^{-2}-\frac{1}{2}(1+5\eta)( Q_\dl^\rho)^4 \right)
P^{H^1}_{V_2(\dl)},
\end{equation}
where

\noi
\begin{align}
V_2(\dl)=\big\{w \in \Re H^1(\T) : \jb{w,(1-\dx^2) Q_\delta^\rho} =0 \big\}.
\label{VV2}
\end{align}

In the following, 
we  proceed to analyze the eigenvalues of $A_1$ and $A_2$. 
In Proposition~\ref{prop: e0-prop}, we provide a lower bound for the eigenvalues of $\frac12\Id+(1-\eta^2)A$.
Then, by  comparing $A_j$ with simple operators whose spectrum can be determined explicitly
(Lemma \ref{lem: eigenvalues}), 
we establish asymptotic bounds
on the eigenvalues
 in Proposition \ref{prop: eigen-asymp}.

\begin{proposition}\label{prop: e0-prop} 
There exists a constant $\eps_0>0$ such that for any  $ \dl , \eta>0$ sufficiently small, the smallest eigenvalue of $A=A(\dl, \eta)$ on $V'$, defined in \eqref{XX8}, is greater than $-\frac{1}{2}+\eps_0$.

\end{proposition}

\begin{proof}
Fix $0 <  \eps_0\ll 1$.
Suppose by contradiction that there exists a sequence $\{ \dl_n\}_{n  \in  \NB}$
of positive numbers tending to 0
and  $w_n^j\in V_j(\dl_n)$, $j = 1, 2$,  such that 
\begin{equation}
\begin{split}
\Big(-\frac{1}{2} & +\eps_0\Big) \int_{\T}(1-\dx^2)w_n^j \cdot w_n^j\ds  x
\\
&
\ge \langle A_j w_n^j, w_n^j\rangle_{H^1(\T)}
= \dl_n^{-2}\int_{\T} (w_n^j)^2\ds  x-c_j(1+5\eta)\int_{\T} ( Q_{\dl_n}^\rho)^4 (w_n^j)^2\ds  x
\end{split}
\label{XA2}
\end{equation}

\noi
for $j = 1, 2$, 
where $c_1 = \frac 52$ and $c_2=\frac12$.
Let $\wt A^j_n$ denote 
 the Schr\"odinger operator  given by
\begin{align}
\begin{split}
\wt A^j_n
& \hphantom{:}
=
P_{V_j(\dl_n)}T^j_nP_{V_j(\dl_n)}\\
 &  : = 
P_{V_j(\dl_n)}
\bigg(-\Big(\frac{1}{2}-\eps_0\Big)\dx^2 + \dl_n^{-2}+\Big(\frac{1}{2}-\eps_0\Big)-c_j(1+5\eta)
( Q_{\dl_n}^\rho)^4\bigg)P_{V_j(\dl_n)}, 
\end{split}
\label{XA2a}
\end{align}

\noi
where 
$P_{V_j(\dl_n)}$ 
denotes the projection onto $V_j(\dl_n)$
in $\Re L^2(\T)$. 
Then, by the min-max principle
(see \cite[Theorem XIII.1]{RS4}), 
the minimum eigenvalue 
 $\wt \ld_n^j$ is non-positive.
 We denote by 
$\wt w^j_n \in V_j(\dl_n)\subset L^2(\T)$    
an $L^2$-normalized eigenvector 
associated with this minimum eigenvalue $\wt \ld_n^j$:
\begin{align}
\wt A_n^j \wt w_n^j = \wt \ld^j_n \wt w^j_n.
\label{XA3}
\end{align}

\noi
Recalling $\|Q_\dl^\rho\|_{L^\infty(\T)} \sim \dl^{-\frac 12}$, 
we have 
\[ \jb{T_n^j w, w}  \ges -\dl_n^{-2}\|w\|_{L^2(\T)}^2\]

\noi
for any $w \in V_j(\dl_n)$, 
which shows that 
$\wt A^j_n$ is semi-bounded on $V_j(\dl_n)$
with a constant of order $\dl_n^{-2}$. 
In particular, we have 
\begin{align}
-C \dl_n^{-2} \leq \wt \ld_n^j \leq 0.
\label{XA3a}
\end{align}

\noi

When $j = 1$, 
the condition 
 $\langle  \wt w^1_n,  Q_{\dl_n}^\rho\rangle=0$
 in \eqref{VV1}
 together with 
 the positivity of $Q$ and the definition $Q_{\dl_n}^\rho =  \rho Q_{\dl_n}$ implies that $\wt w^1_n(x_n)=0$ for some $x_n\in \T$. 
We define a sequence $\{v^1_n\}_{n \in \NB}$ of $L^2$-normalized functions on $\R$ by
\begin{align}
v^1_n(x)=\begin{cases}\dl_n^{\frac 12 }
\wt w_n^1(\dl_nx+x_n-\frac{1}{2}),& 
\quad  \text{for }
|x|\le \frac{1}{2}\dl_n^{-1},\\
0,&\quad  \text{for } |x|>\frac{1}{2}\dl_n^{-1}, 
\end{cases}
\label{XA4}
\end{align}

\noi
where the addition  is understood mod 1. 
When $j = 2$, it follows from the continuity of $\wt w_n^2$
and $\|\wt w_n^2\|_{L^2(\T)} = 1$ that 
$|\wt w^2(x_n)| \leq 1$ for some $x_n \in \T$.
Then, we define a sequence $\{v^2_n\}_{n \in \NB}$ of 
functions on $\R$ by
\begin{align}
v^2_n(x)=\begin{cases}\dl_n^{\frac 12 }
\wt w_n^2(\dl_nx+x_n-\frac{1}{2}),& 
\quad  \text{for }
|x|\le \frac{1}{2}\dl_n^{-1},\\
0,&\quad  \text{for } |x|>\frac{1}{2}\dl_n^{-1} + 
\dl_n^{\frac 12 } |\wt w^2(x_n)|, 
\end{cases}
\label{XA44}
\end{align}

\noi
and by linear interpolation
for  $\frac{1}{2}\dl_n^{-1} < |x|
\le \frac{1}{2}\dl_n^{-1} + 
\dl_n^{\frac 12 } |w(x_n)|$.
In both cases, 
 we have $v_n^j \in H^1(\R)$.
In the remaining part of the proof,  we drop the superscript $j$ 
for simplicity of notations, when there is no confusion.

Define $\wt T_n^j$ by 
\begin{align}
\begin{split}
\wt {T}_n^j
&=-\Big(\frac{1}{2}-\eps_0\Big)\dx^2+ 1
+ 
\Big(\frac{1}{2}-\eps_0\Big)\dl_n^2
 -c_j(1+5\eta)
\big(\tau_{- \dl_n^{-1}(x_n - \frac 12)}\rho(\dl_n \cdot) Q\big)^4, 
\end{split}
\label{XA4a}
\end{align}

\noi
where  $\tau_{x_0}$ denotes the translation by $x_0$
as in \eqref{trans1}.
Then, 
from \eqref{XA2a}, \eqref{XA4}, \eqref{XA44}, and \eqref{XA4a}, 
a direct computation shows that 
\begin{align}
\wt T_n^j v_n (x) = \dl_n^{\frac 52}T_n^j \wt w_n(y)
\label{XA4b}
\end{align}

\noi
with $y = \dl_n x + x_n - \frac 12$, 
 as long as $|x|\le \frac{1}{2}\dl_n^{-1}$.

Define $\wt V_j\subset \Re H^1(\R)$, $j = 1, 2 $, by 
\begin{align}\label{XV1}
\begin{split}& \wt V_1 = \big\{v\in \Re H^1(\R): \langle v, \partial_{\dl} 
Q\rangle_{\dot H^1(\R)} =0 , \langle v,\partial_{x_0} 
Q \rangle_{\dot H^1(\R)}=0, \\
& \hspace{37.3mm}\langle v, Q \rangle_{L^2(\R)} =0\big\},\\
\end{split}\\
& \wt V_2 = \big\{v\in \Re H^1(\R): \jb{v,Q}_{\dot H^1(\R)}=0\big\}.
 \label{XV2}
\end{align}

\noi
Noting that $\jb{v,w}_{\dot H^1(\R)} = \jb{v,-\dx^2 w}_{L^2(\R)}$
and that $Q$ is a smooth function, we can view $\wt V_j$ as a subspace of $\Re L^2(\R)$.
We now define the operator $\wt T_j$, $j = 1, 2$,  on $\Re L^2(\R)$ by 
\begin{align}
\begin{split}
\wt {T}_j
& \hphantom{:}=
P_{\wt{V}_j}^{L^2(\R)}S_jP_{\wt{V}_j }^{L^2(\R)}\\
&:=
P_{\wt{V}_j }^{L^2(\R)}\bigg(-\Big(\frac{1}{2}-\eps_0\Big)\dx^2+ 1 -c_j(1+5\eta)
Q^4\bigg)P_{\wt{V}_j }^{L^2(\R)}. 
\end{split}
\label{XA5}
\end{align}
%
%

\noi
Here, 
$P_{\wt{V}_j}^{L^2(\R)}$
denotes the projection onto $\wt{V}_j$
in $\Re L^2(\R)$.
Then, from \eqref{XA3}, \eqref{XA4a},  \eqref{XA4b}, and \eqref{XA5}
along with 
 the smoothness and exponential decay of $Q$ and its derivatives, 
 we have 
\begin{equation}
\jb{\wt{T}_j \tau_{\dl_n^{-1}(x_n - \frac 12)}v_n - \ld_n \tau_{\dl_n^{-1}(x_n - \frac 12)} v_n, \tau_{\dl_n^{-1}(x_n - \frac 12)} v_n}_{L^2(\R)} \too 0
\label{XA6}
\end{equation}

\noi
as $n \to\infty$, 
where $\ld_n=\dl_n^2\wt \ld_n\le 0$. 
Moreover, 
\eqref{VV1}, \eqref{VV2}, 
\eqref{XA4}, \eqref{XA44},
\eqref{XV1}, and \eqref{XV2}
along with \eqref{rhodef} and 
 the exponential decay of $Q$ and its derivatives
 once again, we have  
 \begin{align}
 \|P_{\wt{V}_j}^{L^2(\R)} \tau_{\dl_n^{-1}(x_n - \frac 12)}v_n \|_{L^2} \gtrsim 1.
 \label{VV3}
 \end{align}
Note that from \eqref{XA3a}, we have
\[-C \le \ld_n \le 0\]

\noi
for any $n \in \NB$.
Thus, passing to a subsequence, we may assume that $\ld_n\to \ld$
for some $\ld\le 0$ and thus from \eqref{XA6}, we obtain
\begin{equation}\label{eqn: T-convergence}
\jb{\wt{T}_j \tau_{\dl_n^{-1}(x_n - \frac 12)} v_n - \ld \tau_{\dl_n^{-1}(x_n - \frac 12)} v_n, \tau_{\dl_n^{-1}(x_n - \frac 12)} v_n}_{L^2(\R)} \too 0.
\end{equation}

Noting that $\wt T_j$ is semi-bounded, 
let $\ld_0$ be the minimum 
of the spectrum of $\wt{T}_j$ on $\wt V_j$. 
By the min-max principle
(see \cite[Theorem XIII.1]{RS4}), 
we have
\begin{equation*}
\jb{\wt{T}_j v - \ld v, v}_{L^2(\R)} \ge (\ld_0 - \ld) \|P_{\wt{V}_j}^{L^2(\R)} v\|_{L^2(\R)}^2
\end{equation*}
for every $v \in H^1(\R)$.
Therefore, from \eqref{eqn: T-convergence}
and \eqref{VV3}, we obtain that $\ld_0 \le \ld \le 0$.

Since $Q$ is a Schwartz function, 
 the essential spectrum of the Schr\"odinger operator $S_j$ in~\eqref{XA5} 
 is equal to $[1,\infty)$;
see \cite[Theorem V-5.7]{Kato}.
Moreover, 
$\wt T_j$ in \eqref{XA5} differs from $S_j$  by a finite rank 
 perturbation
and thus its essential spectrum is $[1,\infty)$;
see \cite[Theorem IV-5.35]{Kato}.
In particular, $\ld_0 \le 0$ does {\it not} lie in the essential spectrum of $\wt T_j$.
Namely, $\ld_0$ belongs to the discrete spectrum of $\wt T_j$.
Hence, there exists $v \in H^2(\R)$ such that 
\begin{equation}\label{eqn: T-lambda}
\wt T_jv=\ld_0 v.
\end{equation}
In order to derive a contradiction, we invoke the following lemma.


\begin{lemma}\label{lem: B-lwr-bd}
For any sufficiently small $\eps_0 > 0$ and $\eta > 0$, 
the following statements hold.
The operator
\begin{equation*}
B_1(\eps_0, \eta) = -\Big(\frac{1}{2}-\eps_0\Big)\dx^2+1-\frac{5}{2}(1+ 5\eta)Q^4,
\end{equation*}
viewed as an operator on $L^2(\R)$, 
is strictly positive on $\wt V_1 $ defined in \eqref{XV1}.
Similarly, the operator 
\begin{equation*}
B_2(\eps_0, \eta) = -
\Big(\frac{1}{2}-\eps_0\Big)\dx^2 + 1-\frac{1}{2}(1+5\eta)Q^4, 
\end{equation*}

\noi
viewed as an operator on $L^2(\R)$, 
is strictly positive on $\wt V_2$ defined in \eqref{XV2}.

\end{lemma}

This lemma  shows 
that 
 \eqref{eqn: T-lambda} can  not hold for  $\ld\le 0$.
Therefore,  we arrive at a contradiction to \eqref{XA2}.
This concludes the proof of 
Proposition \ref{prop: e0-prop} 
(modulo the proof of Lemma \ref{lem: B-lwr-bd} which we  present below).
\end{proof}

We now present 
the proof of Lemma \ref{lem: B-lwr-bd}. 

\begin{proof}[Proof of Lemma \ref{lem: B-lwr-bd}]
In the following, we only prove the strict positivity
of $B_j := B_j(0, 0)$ on $\wt V_j$, $j = 1, 2$, 
when $\eps_0 = \eta = 0$.
Namely, we show that there exists $\ta > 0$ such that 
\begin{equation}
\jb{B_j v, v}_{L^2(\R)} \ge \theta \|v\|_{L^2(\R)}^2
\label{XH3}
\end{equation}

\noi
for any $v\in \wt V_j$.
Then, by writing 
\[ B_j(\eps_0, \eta) = (1 - 2\eps_0) B_j(0, 0)
+ 2\eps_0\big( 1- c_j Q^4\big)- 5c_j \eta Q^4\]

\noi
with $c_1 = \frac 52$ and $c_2=\frac12$, 
the strict positivity of 
$B_j(\eps, \eta)$  for  sufficiently small $\eps_0, \eta > 0$
 follows from \eqref{XH3}.

Consider the Hamiltonian $H(u) =H_\R(u)$:
$$H(u) = \frac12 \int_\R |u'|^2 \d x-\frac16 \int_\R |u|^6 \d x .$$

\noi
By the sharp Gagliardo-Nirenberg-Sobolev  inequality
(Proposition \ref{THM:W}), 
we know that $H: H^1(\R) \to \R$ has a global minimum at $Q$, when restricted to the manifold 
$\norm{u}_{L^2(\R)} = \norm{Q}_{L^2(\R)}$. 
In view of \eqref{elliptic2}, 
we see that $u = Q$ and $\ld = -1$
satisfy the following 
 Lagrange multiplier problem:
\begin{align}
d H|_u(v) = \ld d G|_u(v)
\label{LG1} 
\end{align}

\noi
for any $v \in H^1(\T)$, 
where   $G(u) = \|u \|_{L^2(\R)}^2 - \|Q\|_{L^2(\R)}^2$.
By a direct computation, 
the second variation of  $H(u) + G(u)$ at $Q$ in the direction $v$ is given by 
$ 2\jb{B_1\Re v,\Re v}_{L^2(\R)} + 2\jb{B_2\Im v,\Im v}_{L^2(\R)}$, 
while the constraint $G(u) = 0$ gives $\jb{v, Q}_{L^2(\R)} = 0$.
Then, 
from the second derivative test for 
constrained minima, 
we obtain that 
\begin{align}
&  \jb{B_1w,w}_{L^2(\R)} \ge 0 \quad \text{on } \big\{w \in \Re H^1(\R): \jb{w,Q}_{L^2(\R)} =0\big\}, 
\label{XH4}\\
& \jb{B_2w,w}_{L^2(\R)} \ge 0 \quad \text{for every } w \in \Re H^1(\R). \notag
\end{align}

From \eqref{elliptic}, we have  $B_2 Q = 0$.
Also by differentiating \eqref{elliptic2}, 
we obtain $B_1 Q' = 0$.
Moreover, 
by the elementary theory of Schr\"odinger operators,\footnote{If there is a linearly independent solution $v$ to $B_j v = 0$, then by considering the Wronskian, we obtain that $v' (x)\not\to 0 $ as $|x| \to \infty$.  This in particular implies $v \notin L^2 (\R) \cup \dot H^1(\R)$.} 
 the kernel of  the operator $B_j$, $j = 1, 2$,  has dimension 1. 
In particular, when $j = 2$, $B_2 Q = 0$ implies
\eqref{XH3} for any $v \in \wt V_2$. 
Hereafter, we use the fact that 
the  essential spectrum of $B_j$ is $[1,\infty)$
and that the spectrum of $B_j$ on $(-\infty, 1)$
consists only of isolated eigenvalues of finite multiplicities.\footnote{This claim follows from 
Weyl's criterion (\cite[Theorem VII.12]{RS1}).}

In the following, we focus  on $B_1$.
With $Q^\perp = (\text{span} (Q))^\perp$, 
let $P_{Q^\perp}$
denote the $L^2(\R)$-projection onto $Q^\perp$ given by 
\[P_{Q^\perp} v= v- \frac{\langle Q, v\rangle_{L^2(\R)}}{\|Q\|_{L^2(\R)}^2}Q.\]

\noi
Then, from \eqref{XH4}, we have
\begin{equation}\label{eqn: PBP}
\big\langle P_{Q^\perp}B_1P_{Q^\perp} v,v\big\rangle_{L^2(\R)} \ge 0.
\end{equation}

\noi
Obviously, 
the quadratic form \eqref{eqn: PBP} vanishes for $v=Q$. 
On $Q^\perp$, we have $P_{Q^\perp} B_1 v=0$ if and only if 
\begin{equation}\label{eqn: Bvnull}
B_1v\in \mathrm{ker}P_{Q^\perp}=\mathrm{span}(Q).
\end{equation}

\noi
Recalling that 
 the restriction of $B_1$ to $(Q')^\perp$ is invertible, 
 we see that the condition \eqref{eqn: Bvnull} holds at most for a two-dimensional space. 
Recall from Remark \ref{REM:ortho} that  $Q' \, (= - \dd_{x_0} Q)$
is orthogonal to $Q$ in $L^2(\R)$.
Moreover, by differentiating  
 \eqref{elliptic2}
 in $\dl$, we have 
\[B_1 \dd_\dl Q=2Q\in \mathrm{ker}P_{Q^\perp},\]

\noi
while $\dd_\dl Q \perp Q$ in $L^2(\R)$.
Hence, we have 
 $P_{Q^\perp}B_1P_{Q^\perp}v=0$ 
 on the three-dimensional subspace $\text{span}\{Q,\partial_{x_0} Q, \partial_{\dl}Q\}$.
 As a result, it follows from \eqref{XV1} that there exists $\ta > 0$ such that 
\begin{equation*}
\langle B_1 v, v\rangle_{L^2(\R)} \ge \theta \|v\|_{L^2(\R)}^2
\end{equation*}
for any $v\in \wt {V}_1$.
\end{proof}

\begin{remark} \label{REM:Leo2}
\rm 
The spectral analysis of the operators $B_1$ and $ B_2$ defined in Lemma \ref{lem: B-lwr-bd} resembles closely the analysis of similar Schr\"odinger operators that is at the basis of the results in  \cite{ns1,ns2,ns3,ns4,kns1,kns2, ds12, kns3,kns4,kns5}. 
For example, focusing on the operator $B_1(0,0)$, we have the following picture: 
\begin{itemize}
\item $B_1(0,0)$ has one negative eigenvalue, whose existence can be inferred from 
$$\jb{B_1(0,0) Q, Q}_{L^2(\R)} = - 2 \int Q^6 < 0.$$

\smallskip

\item the eigenvalue $0 \in \sigma(B_1(0,0))$ has multiplicity $2$, corresponding to the directions on the tangent space,  $\mathrm{span}\{\partial_\dl Q, \partial_{x_0} Q\}$,  to the soliton manifold .

\smallskip
\item On the resulting space of codimension 3, corresponding to the orthogonal 
complement to the eigenfunctions described above, $B_1(0,0)$ is strictly positive. 
\end{itemize}

\smallskip

\noi
Compare this result with \cite[Lemma 3.2]{nsbook} and \cite[Lemma 3.7]{ds12}, which depict a very similar picture in their respective cases. 
We would like to point out that, in the current work, we make use of the orthogonal coordinate system around the soliton manifold to remove the complications coming from the $0$ eigenvalues, while the papers cited above deal with the problem differently. 
In particular, we require the strict positivity (see \eqref{eqn: lambdan-uprbd}) of the operators 
$\frac 12 \Id + A_j$ defined in \eqref{Areal} and \eqref{Aimaginary} in order to guarantee that  the denominator appearing in \eqref{eqn: AN-bound}
does not vanish.
Moreover, in the same expression, we crucially use the asymptotics of the positive eigenvalues provided by \eqref{eqn: lambdan-lwrbd}, while this finer analysis does not seem to be present in the aforementioned works.

\end{remark}

Next, we establish asymptotic bounds
on the eigenvalues of $A_j$, $j = 1, 2$.
We achieve this goal by comparing $A_j$
to a simpler operator whose spectrum is studied
in Lemma \ref{lem: eigenvalues} below.

\begin{proposition}\label{prop: eigen-asymp}
Let $j = 1, 2 $. The spectrum of the operators $A_j$
defined  in \eqref{Areal} or \eqref{Aimaginary} consists of a countable collection of eigenvalues. 
Denoting by $\{-\ld_n^{j}\}_{n \in \NB}$
 the negative eigenvalues of $A_j$
 and by $\{\mu_n^{j}\}_{n \in \NB}$ the positive eigenvalues, we have
 \begin{align}
& 1-2\ld_n^{j}\ge 2\eps_0>0,\label{eqn: lambdan-uprbd}\\
& \ld_n^j\lesssim \frac{1}{n^2}, 
 \qquad \text{and}\qquad
\mu_n^j\gtrsim \frac{\dl^{-2}}{n^2}.
 \label{eqn: lambdan-lwrbd}
\end{align}

\end{proposition}

\begin{proof}
The first bound  \eqref{eqn: lambdan-uprbd} is just a restatement of Proposition \ref{prop: e0-prop}.
In the following, we establish  the asymptotic behavior \eqref{eqn: lambdan-lwrbd}
of $\lambda_n^j$ and $\mu_n^j$.

For $j = 1, 2$, define an operator $T_j$ by 
\begin{align}
T_j =(1-\dx^2)^{-1}\big(\dl^{-2}-c_j(1+5\eta)( Q_\dl^\rho)^4 \big),
\label{XT1}
\end{align}

\noi
where $c_1=\frac 52$ and $c_2=\frac 12$.
Since $T_j$ is
 the composition of a bounded multiplication operator and $(1-\dx^2)^{-1}$, which is a compact operator on $H^1(\T)$, it is compact and thus has a countable sequence of eigenvalues accumulating only at zero.  We label the negative eigenvalues 
 as $-\bar\ld_1^j\le -\bar\ld_2^j \le \cdots \leq 0$, while the positive eigenvalues are labeled as $\bar\mu_1^j\ge \bar\mu_2^j\ge \cdots \geq 0$. Since $A_j$ is the composition of $T_j$ 
 with a projection, $A_j$ is also compact.
 Then, by the min-max principle, we have 
\begin{align}
-\bar{\ld}_n^j\le -\ld_n^j \le -\bar{\ld}_{n+3}^j
\qquad 
\text{and}\qquad 
\bar{\mu}_n^j \le \mu_n^j \le \bar{\mu}_{n+3}^j.
\label{XT1a}
\end{align}

\noi
Hence, it suffices to prove the asymptotic bounds
 \eqref{eqn: lambdan-lwrbd} 
 for 
 the eigenvalues $- \bar \lambda_n^j$ and $\bar \mu_n^j$
 of $T_j$.

We estimate 
the eigenvalues $\bar \ld_n^j$ and  $\bar \mu_n^j$ for large $n$ by comparing $T_j$
with an operator with piecewise constant coefficients, whose spectrum can be computed explicitly. 
Fix 
 $a> 0$  such that $10Q^4(a)\le 1$. 
 Given small $\dl > 0$, define the function $S = S_\dl$ on $[-\frac 12,\frac 12)$ by
setting
\begin{align}
S(x)=\begin{cases}
 \dl^{-2} (2 - 5\|Q\|_{L^\infty(\R)}^4), &  \text{for } |x|\le a\dl, \\
\frac 12 \dl^{-2}, & \text{for } a\dl < |x| \leq \frac 12.
\end{cases}.
\label{XT2}
\end{align}

\noi
Recalling the explicit formula
for the ground state 
$ Q(x) = \frac{6^\frac14}{\cosh^\frac 12 (2^\frac 32 x)}$, 
we have $\|Q\|_{L^\infty(\R)}^4 = Q^4(0) = 6$.
Hence, for any sufficiently small $\eta > 0$, we have 
$S\le \dl^{-2}-\frac 52(1+\eta)( Q_\delta^\rho)^4$
pointwise.
Therefore, from \eqref{XT1}, we have
 \[ (1-\dx^2)^{-1}S\le T_j\] 
 
 \noi
 for $j = 1, 2$ in the sense of operators on $H^1(\T)$. 
 By the min-max principle, it then suffices to establish the asymptotic behavior 
 \eqref{eqn: lambdan-lwrbd} for the eigenvalues 
  of $(1-\dx^2)^{-1}S$. This is an explicit computation which we carry out in the following lemma.

\begin{lemma}\label{lem: eigenvalues}
Given  $0 < \dl \ll 1$, 
define the operator $R = R_\dl$
by \[ R =(1-\partial_{x}^2)^{-1}S\]
 on $H^1(\T)$, 
where $S$ is as in \eqref{XT2}. Then,  $R$ has a countable sequence of eigenvalues 
$\{-\wt {\ld}_n\}_{n \in \NB}\cup\{ \wt {\mu}_n\}_{n \in \NB}$ with 
$\wt {\ld}_n$, $\wt {\mu}_n\ge 0$. Moreover, we have the following asymptotics for the eigenvalues\textup{:}
\begin{align}
\wt{\ld}_n\sim \frac{1}{n^2}
\qquad \text{and}\qquad  \wt {\mu}_n\sim \frac{\dl^{-2}}{n^2}, 
\label{XT3}
\end{align}

\noi
where the implicit constants are independent of 
$0 < \dl \ll 1$.

\end{lemma}

We first complete the proof of Proposition 
\ref{prop: eigen-asymp}, assuming Lemma \ref{lem: eigenvalues}.
It follows from the asymptotics \eqref{XT3}
and the min-max principle  that 
\begin{align}
\bar{\ld}_n \les \frac{1}{n^2}
\qquad \text{and}\qquad  \bar {\mu}_n\ges \frac{\dl^{-2}}{n^2}
\label{XT4}
\end{align}

\noi
for  the eigenvalues $- \bar \lambda_n^j$ and $\bar \mu_n^j$
 of $T_j$.
Then, the  asymptotic bounds \eqref{eqn: lambdan-lwrbd}
 for the eigenvalues 
of  $- \lambda_n^j$ and $\mu_n^j$ of $A_j$
follows from \eqref{XT1a} and \eqref{XT4}.
\end{proof}

We now present 
the proof of Lemma \ref{lem: eigenvalues}.

\begin{proof}[Proof of Lemma \ref{lem: eigenvalues}]
Noting that $S$ defined in \eqref{XT2} is a bounded operator, 
we see that 
 $R=(1-\partial^2_{x})^{-1}S$ is a compact operator on $H^1(\T)$ and thus has a countable sequence of eigenvalues $\{\ld_n\}_{n \in \NB}$ accumulating only at 0
 with the associated 
 eigenfunctions $\{f_n\}_{n \in \NB}$ forming 
 a complete orthonormal system in $H^1(\T)$.
Moreover, since the operator $R$ has even and odd functions as invariant subspaces, 
we  assume that any such eigenfunction $f_n$ is  even or odd.
The eigenvalue equation:
\begin{align}
(1-\partial_{x}^2)^{-1}S f_n=\ld_n f_n
\label{XY0}
\end{align}
shows that the second derivative of any eigenfunction is piecewise continuous. 
With $C_0 = 5\|Q\|_{L^\infty(\R)}^4-2 = 
5\cdot 6-2 = 28$,  
define
\begin{equation}
\begin{split}
A_0 &= \ld^{-1} \dl^{-2} C_0+1, \\
B_0 &= \tfrac 12 \ld^{-1} \dl^{-2} - 1.
\end{split}
\label{XY1}
\end{equation}

\noi
Then, by dropping the subscript $n$, 
we can rewrite the eigenvalue equation \eqref{XY0} as
\begin{align}
\begin{array}{ll}
f''=A_0f,  & \text{for } x\in [-a\dl,a\dl],  \rule[-2.5mm]{0pt}{0pt} \\ 
f''=-B_0f ,\quad & \text{for } x\in [-\tfrac 12, -a\dl)\cup(a\dl,\tfrac 12]. 
\end{array}
\label{XY2}
\end{align}

\noi
Let $A = \sqrt{|A_0|}$ and $B = \sqrt{|B_0|}$. 
In the following, we carry out case-by-case analysis, 
depending on the signs of $A_0$ and $B_0$.
It follows from \eqref{XY1} that $\ld > 0$ implies $A_0 > B_0$, while 
$\ld < 0$ implies $B_0 < 0$.
This in particular implies that the case  $A_0 < 0 < B_0$
can not happen for any parameter choice.

\medskip

\noi
$\bullet$ 
\textbf{Case 1}: $A_0,B_0 < 0$ and $f$ is even.
\quad Without loss of generality, we may assume that 
\begin{align*}
f(0)=1, \quad 
f'(0)=0,\quad \text{and}\quad 
f'(\tfrac 12 )=0.
\end{align*}

\noi
From solving \eqref{XY2} on $[-a\dl, a\dl]$, we have
\[f(x)=\cos(Ax).\]

\noi
On the other hand, 
the general solution to \eqref{XY2} on $(a\dl,\frac 12]$ is given by 
\[f(x)=\al \cosh(B(x-a\dl))+\beta \sinh (B(x-a\dl)), \qquad x\in (a\dl, \tfrac 12].\]

\noi
Thus, with the notations:
\begin{align*}
g(x+):=\lim_{\eps\downarrow 0} g(x+\eps)
\qquad \text{and}\qquad
g(x-):=\lim_{\eps\downarrow 0} g(x-\eps), 
\end{align*}

\noi
 we have the following ``transmission conditions''
 at $x=a\dl$:
\begin{align*}
f(a\dl-)&=\cos(A a\dl) =f(a\dl+)=\al, \\
f'(a\dl-)&=-A\sin(Aa\dl)=f'(a\dl+)=\beta B.
\end{align*}

\noi
Hence, on $(a\dl,\frac 12]$, we have
\[f(x)= \cos(Aa\dl)\cosh(B(x-a\dl))-\frac{A}{B}\sin(Aa\dl)\sinh(B(x-a\dl)).\]
To enforce the periodicity condition $f'(\frac 12) = 0$, we need
\begin{align*}
0&=f'(\tfrac 12 )\\
&=B\cos(Aa\dl)\sinh(B(2^{-1}-a\dl))-A\sin(Aa\dl)\cosh(B(2^{-1}-a\dl)).
\end{align*}

\noi
By rearranging, we obtain the condition
\begin{equation}
\frac{A}{B}\tan(Aa\dl) \coth(B(2^{-1}-a\dl))=1 \label{C1}.
\end{equation}

Our goal is to show that for every $k \in \Z_{\ge 0}$, there exists exactly one $A_k^1$ such that \eqref{C1} is satisfied and 
\begin{equation} \label{C2}
\max\Big(k \pi - \frac \pi 2,0\Big) <  A^1_k a \dl  < k \pi + \frac \pi 2. 
\end{equation}

\noi
From \eqref{XY1}, we have 
\begin{equation} \label{C3}
 B = \sqrt{1 + \frac{A^2 + 1}{2C_0}}.
\end{equation}

\noi
Thus, 
we can rewrite \eqref{C1} as $F_1(A) = 1$, 
where $F_1(A)$ is given by 
\begin{equation*}
F_1(A) = \frac{A}{\sqrt{1 + \frac{A^2 + 1}{2C_0}}} \tan(Aa\dl)
  \coth\Bigg(\sqrt{1 + \frac{A^2 + 1}{2C_0}}(2^{-1}-a\dl)\Bigg).
\end{equation*}

\noi
The existence of $A^1_k$ satisfying $F_1(A^1_k) = 1$ follows from continuity
of $F_1$, together with the limits 
$$ F_1(0) = 0,\quad F_1\Big(\frac{(2k-1)\pi}{ 2a \dl} +\Big) = - \infty, \quad 
\text{and} \quad 
F_1\Big(\frac{(2k+1)\pi}{ 2a \dl} -\Big) =  \infty. $$
In order to show uniqueness, it suffices to show that for every $A$ with $F_1(A) = 1$, 
we have  $F_1'(A) > 0$. 
We first note that, for $0\le A \ll \delta^{-1}$, we have 
$$ F_1(A) \les A^2\dl. $$

\noi
Thus, in order to have $F_1(A) = 1$, we must have $A \gtrsim \delta^{-\frac 12}$. 
Moreover, for $A \gtrsim \delta^{-\frac12}$, we have 
\begin{equation*} 
F_1(A) = \sqrt{2C_0} \tan(Aa \dl)(1 + O(A^{-1})),
\end{equation*}
from which we obtain $A \gtrsim \dl^{-1}$.
Namely, $F_1(A)= 1$ implies 
$A \gtrsim \dl^{-1}$.

From \eqref{C3}, we have
\begin{align}
\frac AB = \sqrt{2C_0} \, \big(1 + O(A^{-2})\big)\qquad \text{and}\qquad
\frac BA = \frac{1}{\sqrt{2C_0}}  \big(1 + O(A^{-2})\big)
\label{C4b}
\end{align}

\noi
for $A \gg 1$.
Therefore, whenever $F_1(x) = 1$, we have 
\begin{equation} \label{C5}
\begin{aligned}
F_1'(A) 
&=  
\frac{A}{\sqrt{1 + \frac{A^2 + 1}{2C_0}}} 
\coth\Bigg(\sqrt{1 + \frac{A^2 + 1}{2C_0}}(2^{-1}-a\dl)\Bigg) a\dl \big(1 + (\tan(Aa\dl))^2\big)\\
&\phantom{=} 
+ \frac{1}{\sqrt{1 + \frac{A^2 + 1}{2C_0}}} 
 \coth\Bigg(\sqrt{1 + \frac{A^2 + 1}{2C_0}}(2^{-1}-a\dl)\Bigg) \tan(Aa\dl)\\
&\phantom{=} 
- \frac{A^2}{2C_0\big(1 + \frac{A^2 + 1}{2C_0}\big)^\frac32} 
\coth\Bigg(\sqrt{1 + \frac{A^2 + 1}{2C_0}}(2^{-1}-a\dl)\Bigg) \tan(Aa\dl)\\
&\phantom{=}
 +
 \frac{A}{\sqrt{1 + \frac{A^2 + 1}{2C_0}}} 
\frac {d}{dA} \Bigg(\coth\Bigg(\sqrt{1 + \frac{A^2 + 1}{2C_0}}(2^{-1}-a\dl)\Bigg)\Bigg) \tan(Aa\dl)\\
&= a\dl\bigg(\sqrt{2C_0} + \frac{1}{\sqrt{2C_0}}\bigg)   + O(A^{-2})
> 0
\end{aligned}
\end{equation}

\noi
for $A \ges \dl^{-1}$.
Note that we used the relation \eqref{C1}
in handling  the first, second, and third terms after the first equality in \eqref{C5}.

\smallskip

\noi
$\bullet$ 
\textbf{Case 2:} $A_0,B_0 < 0$ and  $f$ is odd. 
\quad 
Without loss of generality, we may assume that 
\begin{align*}
f(0)=0,\quad
f'(0)=A, \quad \text{and}\quad 
f(\tfrac 12)=0.
\end{align*}

\noi
By solving the eigenfunction equation \eqref{XY2} and solving for the transmission conditions
as in Case 1, we have 
\begin{equation*}
f(x) = \sin(A\alpha \dl) \cosh( B(x - a \dl)) + \frac A B \cos(Aa\dl)\sinh( B(x - a \dl))
\end{equation*}

\noi
 for $ x \in (a \dl, \frac 12]$.
Therefore, the periodicity condition $f(\tfrac 12)=0$ becomes
\begin{equation} \label{C6}
F_2(A) : = \frac A B \cot(A \alpha \dl) \tanh(B(2^{-1} - a \dl)) = - 1, 
\end{equation}

\noi
where $B = B(A)$ is as in \eqref{C3}.
As in Case 1, we want to show that for every $k \in \Z_{\ge 0}$, there exists exactly one value $A^2_k$ that satisfies \eqref{C6} with 
\begin{equation}\label{C7}
k\pi < A^2_k a \dl < (k+1) \pi.
\end{equation}

\noi
From \eqref{C6}, 
we have  $\cot(Aa\dl) < 0$, which implies $A a \dl > \frac{\pi}2$.
Then, using \eqref{C4b} and \eqref{C7} with $A \ges \dl^{-1}$, 
we can proceed as in Case 1
and  show existence and uniqueness of such  $A^2_k$, $k \in \Z_{\ge 0}$.
Note that in this case, we show $F_2'(A) < 0$ whenever $F_2(A) = -1$.

\smallskip

\noi
$\bullet$ \textbf{Case 3:}  $B_0 < 0 <A_0$ and  $f$ is even.
\quad By solving the eigenfunction equation \eqref{XY2}, we obtain 
$$ f(x) = \cosh(Aa\dl) \cosh(B(x-a\dl)) + \frac A B \sinh(Aa\dl) \sinh(B(x-a\dl)) 
$$

\noi
for $x \in ( a\dl, \frac 12]$.
By imposing the periodicity condition $f'(\frac 12 ) = 0$, we obtain 
\begin{equation*}
\frac A B  \tanh(Aa\dl) \coth(B(2^{-1}-a\dl)) = - 1.
\end{equation*}
However, since $A, B, a > 0$, we see that this condition can never be satisfied
for $0 < \dl \ll 1$.

\smallskip

\noi
$\bullet$ \textbf{Case 4:}  $B_0< 0<A_0$ and $f$ is odd. 
\quad In this case, we have 
$$ f(x) = \sinh(Aa\dl) \cosh(B(x-a\dl)) + \frac A B \cosh(Aa\dl) \sinh(B(x-a\dl)) 
$$

\noi
for $x \in ( a\dl, \frac 12]$.
By imposing the periodicity condition $f(\frac 12) = 0$, we obtain
\begin{equation*}
\frac A B  \coth(Aa\dl) \tanh(B(2^{-1}-a\dl)) = - 1.
\end{equation*}
Since $A, B, a > 0$, we once again see that  this condition can not be satisfied
for $0 < \dl \ll 1$.

\smallskip

\noi
$\bullet$ \textbf{Case 5:} $A_0, B_0 > 0$ and $f$ is even. 
\quad 
By solving the eigenfunction equation \eqref{XY2}, we have 
$$ f(x) = \cosh(Aa\dl) \cos(B(x-a\dl)) + \frac A B \sinh(Aa\dl) \sin(B(x-a\dl))$$

\noi
for $x \in ( a\dl, \frac 12]$.
Then, by imposing the periodicity condition $f'(\frac 12 ) = 0$, we obtain
\begin{equation} \label{C8}
B\tan(B(2^{-1}-a\dl)) - A \tanh(Aa\dl) = 0.
\end{equation}

\noi
By writing
\begin{equation}\label{C10}
A = \sqrt{2(B^2+1)C_0 + 1},
\end{equation}
we can rewrite \eqref{C8} as $G_1(B) = 0$, where
\begin{equation*}
\begin{split}
G_1(B) 
& = B\tan(B(2^{-1}-a\dl)) \\
& \hphantom{=}- \sqrt{2(B^2+1)C_0 + 1} \tanh\big(\sqrt{2(B^2+1)C_0 + 1}\, a\dl\big) .
\end{split}
\end{equation*}

\noi
As in  Case 1, we need to  show that for every $k \in \Z_{\ge 0}$, there exists exactly one value $B^1_k$ 
such that $G(B^1_k) = 0$ with 
\begin{equation}\label{C9}
\max \Big( k\pi - \frac \pi 2, 0 \Big)  < B^1_k(2^{-1} - a \dl) < k \pi + \frac \pi 2.
\end{equation}

\noi
Existence follows from continuity of $G_1$, together with the limits 
$$ G_1(0) < 0, \quad 
G_1\Big(\frac{2k-1}{1-2a\dl} + \Big) = - \infty, 
\quad 
\text{and}
\quad
G_1\Big(\frac{2k+1}{1-2a\dl} - \Big) =  \infty. $$

\noi
As  for uniqueness of $B^1_k$, 
it suffices to  show 
$G_1'(B) > 0$ whenever
 $G_1(B) = 0$.
Using~\eqref{C8}, we have 
\begin{align*}
G_1'(B) &= \tan(B(2^{-1}-a\dl)) + B(2^{-1}-a\dl) \big(1 + (\tan(B(2^{-1}-a\dl)))^2\big) \\
&\phantom{=} 
- \frac{2C_0B}{\sqrt{2(B^2+1)C_0 + 1}}\tanh\big(\sqrt{2(B^2+1)C_0 + 1}\, a\dl\big)\\
&\phantom{=} - 2 a \dl C_0B
+ 2a\dl C_0 B \big(\tanh(\sqrt{2(B^2+1)C_0 + 1}\, a\dl)\big)^2 \\
&= \tan(B(2^{-1}-a\dl)) + B(2^{-1}-a\dl) \big(1 + (\tan(B(2^{-1}-a\dl)))^2\big) \\
&\phantom{=} - \frac{2C_0B^2}{{2(B^2+1)C_0 + 1}}\tan(B(2^{-1}-a\dl)) +O(B \dl)\\
&\ge  B(2^{-1}-a\dl) + O(B \dl) > 0
\end{align*}

\noi
for $0 < \dl \ll 1$.

\smallskip

\noi
$\bullet$ 
\textbf{Case 6:} $A_0, B_0 > 0$ and  $f$ is odd. 
\quad 
By solving the eigenfunction equation \eqref{XY2}, we have that 
$$ f(x) = \sinh(Aa\dl) \cos(B(x-a\dl)) + \frac A B \cosh(Aa\dl) \sin(B(x-a\dl))
$$

\noi
for $x \in ( a\dl, \frac 12]$.
Then, by imposing the periodicity condition $f(\frac 12) = 0$, we obtain
\begin{equation} 
B\cot(B(2^{-1}-a\dl))\tanh(Aa\dl) + A  = 0.
\label{C11a}
\end{equation}

\noi
In view of \eqref{C10}, define $G_2(B)$ by 
\begin{align*}
\begin{split}
G_2(B) & = B\cot(B(2^{-1}-a\dl))\tanh\big(\sqrt{2(B^2+1)C_0 + 1}\,a\dl \big)\\
& \hphantom{=}
 + \sqrt{2(B^2+1)C_0 + 1}.
 \end{split} 
\end{align*}

\noi
As in the previous cases, 
 we  show that for every $k \in \Z_{\ge 0}$, there exists exactly one value $B^2_k$ 
 such that $G_2(B_k^2) = 0$ and 
\begin{equation}\label{C12}
k\pi < B^2_k(2^{-1} - a \dl) < (k+1) \pi .
\end{equation}

\noi
Noting 
$G_2(0) > 0$, 
existence follows from the same argument as in the previous cases.
As for uniqueness, 
we show that $G_2'(B) < 0$
 whenever $G_2(B) = 0$.
Using~\eqref{C11a}
(which implies $\cot(B(2^{-1}-a\dl)) < 0$) 
and 
noting that $\tanh\big(\sqrt{2(B^2+1)C_0 + 1}\,a\dl\big)$
 is strictly increasing in $B$, 
we have
\begin{align*}
G_2'(B) 
&\le \cot(B(2^{-1}-a\dl))\tanh\big(\sqrt{2(B^2+1)C_0 + 1}\, a\dl\big) \\
&\phantom{\le} - (2^{-1}-a\dl)B \tanh\big(\sqrt{2(B^2+1)C_0 + 1}\,a\dl\big)
\big(1 + (\cot(B(2^{-1}-a\dl)))^2\big)\\
&\phantom{\le} + \frac{2C_0B}{\sqrt{2(B^2+1)C_0 + 1}}\\
&= - \frac{\sqrt{2(B^2+1)C_0 + 1}}{B} 
+ \frac{2C_0B}{\sqrt{2(B^2+1)C_0 + 1}}\\
&\phantom{=} 
- (2^{-1}-a\dl)B\tanh\big(\sqrt{2(B^2+1)C_0 + 1}\,a\dl\big)
\big(1 + (\cot(B(2^{-1}-a\dl)))^2\big)\\
&\le 
- (2^{-1}-a\dl)B\tanh\big(\sqrt{2(B^2+1)C_0 + 1}\,a\dl\big)
\big(1 + (\cot(B(2^{-1}-a\dl)))^2\big)\\
&<0
\end{align*}

\noi
for $0 < \dl \ll 1$.

\medskip

\noi
{\bf Conclusion:}
Fix $0 < \dl \ll 1$.
Putting all the cases together, 
 we conclude  that 
  every $\lambda$ 
such that in \eqref{XY1}, we have $A_0 = -(A_k^1)^2$, $A_0 = -(A_k^2)^2$, $B_0 = (B_k^1)^2$, or $B_0 = (B_k^2)^2$ 
for some $k \in \Z_{\ge 0}$
corresponds to  an eigenvalue for $R$.
Moreover,  this exhausts all the possibilities for the eigenvalues of $R$, with possible exceptions of $\lambda$ such that $A_0 = 0$ or $B_0 = 0$. 
By inverting the formulas of $A_0= -(A_k^j)^2$ and  $B_0= (B_k^j)^2$ in \eqref{XY1} for  $\lambda$
together with  \eqref{C2}, \eqref{C7}, \eqref{C9},  \eqref{C12}, 
and 
$A_k^j \ges \dl^{-1}$, $j = 1, 2$, 
we 
obtain the asymptotics:
$$\wt{\ld}_n\sim \frac{1}{n^2}\qquad \text{and}\qquad \wt{\mu}_n\sim \frac{\dl^{-2}}{n^2}. $$

\noi
This completes the proof of Lemma \ref{lem: eigenvalues}.
\end{proof}

Finally, 
we conclude this subsection by presenting the proof of Proposition \ref{PROP: gauss-int}.

\begin{proof}[Proof of Proposition \ref{PROP: gauss-int}]
Given  $N\in \mathbb{N}$, let $W_N$ be the subspace spanned by the eigenvectors 
of $A$ corresponding to 
the eigenvalues $\ld_n^j$ and $\mu_n^j$, $1\le n\le N$, $j = 1, 2$. 
Define $A_N=AP_N$, where $P_N = P_{W_N}^{H^1}$ is the projection onto $W_N$
in $H^1(\T)$. Then,  by Proposition~\ref{prop: eigen-asymp}
with $\d \al \d\be = \prod_{j = 1}^2 \prod_{n = 1}^N \d \al_n^j\d \be_n^j$
and \eqref{H3a}, we have
\begin{align}
& \int   \exp(-(1-\eta^2) \langle A_N w,w\rangle_{H^1(\T)})\ds \mu_\dl^{\perp\perp}(w) 
\notag \\
&\hphantom{X}
=\frac{1}{(2\pi)^{2N}}\int_{\R^{4N}}
 \exp\bigg(-\frac{1}{2}\sum_{j=1}^2\sum_{n=1}^N(1-2(1-\eta^2)\ld_n^j)(\al_n^j)^2
\notag \\
& \hphantom{XXXXXXXXXXXl} -\frac{1}{2}\sum_{j=1}^2\sum_{n=1}^N (1+2(1-\eta^2)\mu_n^j)(\be_n^j)^2\bigg)\,\ds  \al\ds  \be  
\notag \\
&\hphantom{X}
=\Bigg(\prod_{j=1}^2\prod_{n=1}^N \frac{1}{\sqrt{1-2(1-\eta^2)\ld_n^j}}\Bigg)
\Bigg(\prod_{j=1}^2\prod_{n=1}^N \frac{1}{\sqrt{1+2(1-\eta^2)\mu_n^j}}\Bigg)
\notag \\
&\hphantom{X}
\le C_\eta \exp\bigg(-\frac{1}{2}\sum_{j=1}^2\sum_{\dl^{-1}\ll n\le N}
\log\big(1+2(1-\eta^2)\mu_n^j\big)\bigg) 
\notag \\
&\hphantom{X}
\le C_\eta  \exp\bigg(-c\sum_{n\gg \dl^{-1}} \frac{\dl^{-2}}{n^2}\bigg)
\lesssim \exp(-c'\dl^{-1}), 
\label{eqn: AN-bound}
\end{align}

\noi
where we used
$\log (1+x) \geq x - \frac {x^2}2 \ges x$
for $|x|\ll 1$ in the penultimate step.

For $v, w \in L^2(\T)$, we have 
\begin{align*}
 \sup_{N \in \NB}|\langle A_N v, w\rangle_{H^1(\T)}|
+ |\langle Av,w\rangle_{H^1(\T)}|\le C\dl^{-2}\|v\|_{L^2(\T)}\|w\|_{L^2(\T)}.
\end{align*}

\noi
Then, by a density argument, 
we see that $\jb{A_N w, w}$ converges
to $\jb{A w, w}$ as $N \to \infty$ 
for each $w \in L^2(\T)$.
Hence, the estimate \eqref{eqn: A-bound} follows from \eqref{eqn: AN-bound}
and  Fatou's lemma. 
\end{proof}

\subsection{Proof of Theorem \ref{THM:OPT}}
\label{SUBSEC:pf}

We now put all the steps together and present the proof of  of Theorem \ref{THM:OPT}.
Let 
$K = \|Q\|_{L^2(\R)}$.
By Lemma \ref{LEM:u-Q}, we have 
\begin{align}
\begin{split}
Z_{6,K=\|Q\|_{L^2(\R)}}
& =\E\Big[e^{\frac{1}{p}\int_\T |u|^6\,\ds  x}\mathbf{1}_{\{\|u\|_{L^2(\T)}\le K\}}\Big]\\
& \leq 
  \E\Big[e^{\frac{1}{6}\int_\T |u(x)|^6\,\ds  x}, 
  S_\g\Big]
  + \E\Big[e^{\frac{1}{6}\int_\T |u(x)|^6\,\ds  x}
\ind_{\{\|u\|_{L^2(\T)} \leq  K\}}
,  U_\eps(\dl^*) \Big]\\
& =: \1(\g) + \II(\eps, \dl^*), 
\end{split}
\label{ZZ1}
\end{align}

\noi
where $S_\g$ and  $U_\eps(\dl^*)$ are as in \eqref{G1} and  \eqref{Ueps}.
From \eqref{E-upper}, 
we have
\begin{align}
\1(\g) < \infty
\label{ZZ2}
\end{align}

\noi
 for any $\g > 0$.
As for $\II(\eps, \dl^*)$, 
it follows from \eqref{XX5a}, 
H\"older's inequality, and Lemma \ref{lem: v6-int}
that 
\begin{align*}
\begin{split}
\II(\eps, \dl^*)& 
\les 
\int_0^{\dl^*}  
\bigg(\int_{\{\|v\|_{L^2(\T)} \le \eps_1\}} \exp\bigg(C'_\eta \int_{\T}|v|^6\,\ds  x\bigg) \mu^\perp_\dl (\d v)
\bigg)^\frac{\eta}{1+\eta}\\
& 
\hphantom{=}
\times \bigg(\int_{\{\|v\|_{L^2(\T)}\le \eps_1\}}
e^{(1+\eta)G(v)}
\ind_{\{\| Q_\dl^\rho+v\|_{L^2(\T)}\le K\}} \mu_\dl^\perp(\d v)\bigg)^{\frac{1}{1+\eta}}\sigma(\ds \dl)\\
& \les 
\int_0^{\dl^*}  
 \bigg(\int_{\{\|v\|_{L^2(\T)}\le \eps_1\}}
e^{(1+\eta)G(v)}
\ind_{\{\| Q_\dl^\rho +v\|_{L^2(\T)}\le K\}} \mu_\dl^\perp(\d v)\bigg)^{\frac{1}{1+\eta}}\sigma(\ds \dl), 
\end{split}
\end{align*}

\noi
provided that $\eps_1 = \eps_1(\eta) > 0$ is sufficiently small.
Note that in obtaining \eqref{XX5a}, we applied
Proposition \ref{PROP:FOC} which gives
an orthogonal coordinate system in a neighborhood of the soliton manifold.
Then, from 
Lemma \ref{LEM:G1}, \eqref{AA1}, and Proposition \ref{PROP: gauss-int}, 
we obtain
\begin{align}
\begin{split}
\II(\eps, \dl^*)
& \les 
\int_0^{\dl^*}  
 \bigg(
 \int \exp\big(-(1-\eta^2)\langle Aw,w\rangle_{H^1(\T)}\big)\ds \mu^{\perp\perp}_\dl(w)
\bigg)^{\frac{1}{1+\eta}}\sigma(\ds \dl)\\
& \les 
\int_0^{\dl^*}  \exp(-c\dl^{-1})
\sigma(\ds \dl) < \infty, 
\end{split}
\label{ZZ3}
\end{align}

\noi
provided that $\eta > 0$ is sufficiently small, 
where we used 
Lemma \ref{LEM:surface} in the last step.
Therefore, from \eqref{ZZ1}, \eqref{ZZ2}, and \eqref{ZZ3}, 
we conclude that 
\[Z_{6,K=\|Q\|_{L^2(\R)}}< \infty. \]

For readers' convenience, 
we go over how we choose the parameters.
We first choose $\eta > 0 $ in \eqref{XX5} sufficiently small
such that  Proposition \ref{PROP: gauss-int}
holds.
Next, we fix 
small $\eps_1 = \eps_1(\eta) > 0$ 
such that  Lemma \ref{lem: v6-int} holds.
Then, Proposition \ref{PROP:FOC}
determines 
small  $\eps = \eps(\eps_1)>0 $ and $\dl^* = \dl^*(\eps_1)> 0 $.
See also Lemma \ref{LEM:G1}, 
where we need smallness of $\dl = \dl(\eta, \eps_1) > 0$.
Finally, we fix $\g =\g(\eps) > 0$
by applying Lemma \ref{LEM:u-Q}.
This completes the proof of  Theorem \ref{THM:OPT}.

\begin{remark}\label{REM:mean}
\rm 
In this section, we presented the proof of Theorem \ref{THM:OPT}, 
where the base Gaussian process is given by 
the  Ornstein-Uhlenbeck loop  in~\eqref{Q0}.
When the base Gaussian process is given by 
the mean-zero  Brownian loop 
in \eqref{bloop}, 
 the same but simpler argument gives 
Theorem \ref{THM:OPT}.
For example, 
in the case of 
the mean-zero  Brownian loop, 
we can omit~\eqref{G1x}
and 
the reduction to the mean-zero case at the beginning of the proof of Lemma~\ref{LEM:u-Q}.
In the proof of Proposition \ref{PROP:FOC}, 
we introduced
$V_{\dl,x_0,\theta}^{\g_0}$
in \eqref{eqn: Vdelta0-def}
in order to use a scaling argument in the non-homogeneous setting.
 In the  case of 
the mean-zero  Brownian loop, 
we can simply use 
$V_{\dl,x_0,\theta}^{0}$, i.e.~$\g_0 = 0$.
The rest of the proof remains essentially the same.

\end{remark}

\begin{ack}\rm
The authors would like to 
 thank  Nikolay Tzvetkov for some interesting discussions.
 The authors would like to express their gratitude  to the anonymous referee for the helpful comments which improved the quality of the paper.
 T.O.~was supported by the European Research Council (grant no.~637995 ``ProbDynDispEq''
and grant no.~864138 ``SingStochDispDyn"). 
P.S.~was partially supported by NSF grant DMS-1811093.
L.T.~was supported by the European Research Council (grant no.~637995 ``ProbDynDispEq'')
and  by the Deutsche
Forschungsgemeinschaft (DFG, German Research Foundation) under Germany's Excellence
Strategy-EXC-2047/1-390685813, through the Collaborative Research Centre (CRC) 1060. 
L.T.~was also supported by the Maxwell Institute Graduate School in Analysis and its Applications, a Centre for Doctoral Training funded by the UK Engineering and Physical Sciences Research Council (grant EP/L016508/01), the Scottish Funding Council, Heriot-Watt University, and the University of Edinburgh.
\end{ack}


\begin{thebibliography}{99}	

%
\bibitem{Agrawal}
G.P. Agrawal, 
{\it Nonlinear Fiber Optics}, 5th ed. San Francisco, CA, USA: Academic, 2012.
%


\bibitem{BG}
N.~Barashkov, M.~Gubinelli, 
{\it  A variational method for $\Phi^4_3$}, 
Duke Math. J.
 169 (2020), no. 17, 3339--3415.

\bibitem{BO}
\'A.~B\'enyi, T.~Oh,
{\it The Sobolev inequality on the torus revisited}, Publ. Math. Debrecen 83 (2013), no. 3, 359--374. 



\bibitem{BOP1}
\'A.~B\'enyi, T.~Oh, O.~Pocovnicu,
{\it  Wiener randomization on unbounded domains and an application to almost sure well-posedness of NLS}, Excursions in harmonic analysis. Vol. 4, 3--25, Appl. Numer. Harmon. Anal., Birkh\"auser/Springer, Cham, 2015. 

\bibitem{BOP2}
\'A.~B\'enyi, T.~Oh, O.~Pocovnicu, 
{\it On the probabilistic Cauchy theory of the cubic nonlinear Schr\"odinger equation on $\R^d$, $d \geq 3$},  Trans. Amer. Math. Soc. Ser. B 2 (2015), 1--50. 




\bibitem{BOP3}
\'A.~B\'enyi, T.~Oh, O.~Pocovnicu, 
{\it  Higher order expansions for the probabilistic local Cauchy theory of the cubic nonlinear Schr\"odinger equation on $\R^3$}, Trans. Amer. Math. Soc. Ser. B 6 (2019), 114--160. 


\bibitem{BOP4}
\'A.~B\'enyi, T.~Oh, O.~Pocovnicu,
{\it  On the probabilistic Cauchy theory for nonlinear dispersive PDEs,} 
 Landscapes of Time-Frequency Analysis. 1--32, Appl. Numer. Harmon. Anal., Birkh\"auser/Springer, Cham, 2019. 


\bibitem{BO93}
J.~Bourgain, 
{\it Fourier transform restriction phenomena for certain lattice subsets and applications to nonlinear evolution equations. I. Schr\"odinger equations,} Geom. Funct. Anal. 3 (1993), no. 2, 107--156. 


\bibitem{BO94}
J.~Bourgain, 
{\it Periodic nonlinear Schr\"odinger equation and invariant measures}, 
Comm. Math. Phys. 166 (1994), no. 1, 1--26.

\bibitem{BO96}
J.~Bourgain, 
{\it Invariant measures for the 2D-defocusing nonlinear Schr\"odinger equation,} Comm. Math. Phys. 176 (1996), no. 2, 421--445. 

\bibitem{BO97}
J.~Bourgain, 
{\it Invariant measures for the Gross-Piatevskii equation,} 
J. Math. Pures Appl. 76 (1997), no. 8, 649--702. 


\bibitem{BBulut} 
J.~Bourgain,  A.~Bulut, 
{\it Almost sure global well posedness for the radial nonlinear Schr\"odinger equation on the unit ball I: the 2D case.,}
Ann. Inst. H. Poincar\'e Anal. Non Lin\'eaire 31 (2014), no. 6, 1267--1288.

\bibitem{Bre}
J.~Brereton, 
{\it Invariant measure construction at a fixed mass},
 Nonlinearity 32 (2019), no. 2, 496--558.



\bibitem{Bring}
B.~Bringmann, 
{\it Stable blowup for the focusing energy critical nonlinear wave equation under random perturbations,}
Comm. Partial Differential Equations 45 (2020), no. 12, 1755--1777.


\bibitem{Bring3}
B.~Bringmann, 
{\it 
Invariant Gibbs measures for the three-dimensional wave equation with a Hartree nonlinearity II: Dynamics}
arXiv:2009.04616 [math.AP].

\bibitem{BS}
D.~Brydges, G.~Slade, 
{\it Statistical mechanics of the 2-dimensional focusing nonlinear Schr\"odinger equation,}
 Comm. Math. Phys. 182 (1996), no. 2, 485--504. 



\bibitem{BT1}
N.~Burq, N.~Tzvetkov, 
{\it Random data Cauchy theory for supercritical wave equations. I. Local theory,}
 Invent. Math. 173 (2008), no. 3, 449--475.


\bibitem{BT3}
N.~Burq, N.~Tzvetkov,
{\it Probabilistic well-posedness for the cubic wave equation},
J. Eur. Math. Soc. 16 (2014), no. 1, 1--30.



\bibitem{CFL} 
E.~Carlen, J.~Fr\"ohlich, J.~Lebowitz, 
{\it Exponential relaxation to equilibrium for a one-dimensional focusing non-linear Schr\"odinger equation with noise},
Comm. Math. Phys. 342 (2016), no. 1, 303--332. 



\bibitem{CK}
A.~Chapouto, N.~Kishimoto,
{\it  Invariance of the Gibbs measures for the periodic generalized KdV equations}, 
arXiv:2104.07382 [math.AP].



\bibitem{CO}
J.~Colliander, T.~Oh, 
{\it  Almost sure well-posedness of the cubic nonlinear Schr\"odinger equation below $L^2(\T)$}, Duke Math. J. 161 (2012), no. 3, 367--414. 



\bibitem{CR}
J.~Colliander, P.~Rapha\"el, 
{\it Rough blowup solutions to the $L^2$ critical NLS,}
Math. Ann. 345 (2009), no. 2, 307--366. 

\bibitem{DP}
G.~Da Prato, 
{\it An introduction to infinite-dimensional analysis,}
Universitext. Springer-Verlag, Berlin, 2006. x+209 pp.



\bibitem{DPD2}
G.~Da Prato, A.~Debussche, 
{\it Strong solutions to the stochastic quantization equations,} Ann. Probab. 31 (2003), no. 4, 1900--1916.



\bibitem{DZ}
 G.~Da Prato, J.~Zabczyk, 
 {\it Stochastic equations in infinite dimensions,} Second edition. Encyclopedia of Mathematics and its Applications, 152. Cambridge University Press, Cambridge, 2014. xviii+493 pp.


\bibitem{DNY2}
Y.~Deng, A.~Nahmod, H.~Yue,
{\it 
Invariant Gibbs measures and global strong solutions for nonlinear Schr\"odinger equations in dimension two},
arXiv:1910.08492 [math.AP].

\bibitem{DNY3}
Y.~Deng, A.~Nahmod, H.~Yue,
{\it 
Random tensors, propagation of randomness, and nonlinear dispersive equations}, 
arXiv:2006.09285 [math.AP].



\bibitem{DNY4}
Y.~Deng, A.~Nahmod, H.~Yue, 
{\it Invariant Gibbs measure and global strong solutions for the Hartree NLS equation in dimension three,}
 J. Math. Phys. 62 (2021), no. 3, 031514, 39 pp. 



\bibitem{Dodson}
B.~Dodson,
{\it A determination of the blowup solutions to the focusing, quintic NLS with mass equal to the mass of the soliton}, arXiv:2104.11690 [math.AP].


\bibitem{ds12} 
R.~Donninger, B.~Sch\"orkhuber, \emph{Stable self-similar blow up for energy subcritical wave equations}, Dyn. Partial Differ. Equ. 9 (2012), no. 1, 63--87. 


\bibitem{driver} 
B.~Driver, 
{\it Analysis tools with applications,} 
lecture  notes, 2003. \url{http://www.math.ucsd.edu/~bdriver/240-01-02/Lecture_Notes/anal.pdf}


\bibitem{EJS}
W.~E, A.~Jentzen, H.~Shen, 
{\it Renormalized powers of Ornstein-Uhlenbeck processes and well-posedness of stochastic Ginzburg-Landau equations,}
Nonlinear Anal. 142 (2016), 152--193. 


\bibitem{EG}
L.C.~Evans, R.~Gariepy, 
{\it Measure theory and fine properties of functions,} Studies in Advanced Mathematics. CRC Press, Boca Raton, FL, 1992. viii+268 pp.



\bibitem{FM}
C.~Fan, D.~Mendelson,
{\it Construction of $L^2$ log-log blowup solutions for the mass critical nonlinear Schr\"odinger equation}, 
arXiv:2010.07821 [math.AP].

\bibitem{fernique} 
X.~Fernique, 
{\it Regularit\'e des trajectoires des fonctions al\'eatoires gaussiennes,}
\'Ecole d'\'Et\'e de Probabilit\'es de Saint-Flour, IV-1974,  1--96. Lecture Notes in Math., Vol. 480, Springer, Berlin, 1975.


\bibitem{frank} 
R.~Frank, \emph{Grounds states of semilinear equations}, 
lecture notes from \emph{Current topics in Mathematical Physics}, Luminy, 2013.
\url{http://www.mathematik.uni-muenchen.de/~frank/luminy140202.pdf}


\bibitem{GIP}
M.~Gubinelli, P.~Imkeller, N.~Perkowski, 
{\it Paracontrolled distributions and singular PDEs,} Forum Math. Pi 3 (2015), e6, 75 pp. 



\bibitem{GKO2}
M.~Gubinelli, H.~Koch, T.~Oh, 
{\it Paracontrolled approach to the three-dimensional stochastic nonlinear wave equation with quadratic nonlinearity}, to appear in J. Eur. Math. Soc. 


\bibitem{HK}
T.~Hmidi, S.~Keraani, 
{\it Blowup theory for the critical nonlinear Schr\"odinger equations revisited,}
 Int. Math. Res. Not. 2005, 2815--2828.



\bibitem{KMM}
S.~Kamvissis, K.~McLaughlin, P.~Miller, 
{\it Semiclassical soliton ensembles for the focusing nonlinear Schr\"odinger equation}, 
Annals of Mathematics Studies, 154. Princeton University Press, Princeton, NJ, 2003. xii+265 pp.

\bibitem{Kato}
T.~Kato, {\it Perturbation theory for linear operators,} Reprint of the 1980 edition. Classics in Mathematics. Springer-Verlag, Berlin, 1995. xxii+619 pp.

\bibitem{KV}
R.~Killip, M.~Vi\c{s}an, 
{\it Nonlinear Schr\"odinger equations at critical regularity,}
 Evolution equations, 325--437, Clay Math. Proc., 17, Amer. Math. Soc., Providence, RI, 2013.


\bibitem{Kishimoto}
N.~Kishimoto, 
{\it Remark on the periodic mass critical nonlinear Schr\"odinger equation,}
 Proc. Amer. Math. Soc. 142 (2014), no. 8, 2649--2660. 



\bibitem{kns1}
J.~Krieger, K.~Nakanishi, W.~Schlag, \emph{Global dynamics above the ground state energy for the one-dimensional NLKG equation}, Math. Z. 272 (2012), no. 1-2, 297--316.

\bibitem{kns2}
J.~Krieger, K.~Nakanishi, W.~Schlag, \emph{Global dynamics of the nonradial energy-critical wave equation above the ground state energy}, Discrete Contin. Dyn. Syst. 33 (2013), no. 6, 2423--2450.

\bibitem{kns3}
J.~Krieger, K.~Nakanishi, W.~Schlag, \emph{Global dynamics away from the ground state for the energy-critical nonlinear wave equation}, Amer. J. Math. 135 (2013), no. 4, 935--965.

\bibitem{kns4}
J.~Krieger, K.~Nakanishi, W.~Schlag, \emph{Threshold phenomenon for the quintic wave equation in three dimensions}, Comm. Math. Phys. 327 (2014), no. 1, 309--332. 

\bibitem{kns5}
J.~Krieger, K.~Nakanishi, W.~Schlag, \emph{Center-stable manifold of the ground state in the energy space for the critical wave equation}, Math. Ann. 361 (2015), no. 1-2, 1--50.
 



\bibitem{kwong} 
M.K.~Kwong, 
{\it Uniqueness of positive solutions of $\Dl u - u + u^p = 0$ in $\mathbf{R}^n$},
Arch. Rational Mech. Anal. 105 (1989), no. 3, 243--266. 



\bibitem{LMW}
J.~Lebowitz, P.~Mounaix, W.-M.~Wang, 
{\it Approach to equilibrium for the stochastic NLS,} Comm. Math. Phys. 321 (2013), no. 1, 69--84.


\bibitem{LRS} 
J.~Lebowitz, H.~Rose, E.~Speer, 
{\it Statistical mechanics of the nonlinear Schr\"odinger equation,}
 J. Statist. Phys. 50 (1988), no. 3-4, 657--687. 


\bibitem{LZ}
D.~Li, X.~Zhang, 
{\it On the rigidity of minimal mass solutions to the focusing mass-critical NLS for rough initial data,}
Electron. J. Differential Equations 2009, No. 78, 19 pp.





\bibitem{LM}
J.~L\"uhrmann, D.~Mendelson,
{\it Random data Cauchy theory for nonlinear wave equations of power-type on $\mathbb{R}^3$},
 Comm. Partial Differential Equations 39 (2014), no. 12, 2262--2283.


\bibitem{MM1}
Y.~Martel, F.~Merle, 
{\it Blow up in finite time and dynamics of blow up solutions for the 
$L^2$-critical generalized KdV equation,} J. Amer. Math. Soc. 15 (2002), no. 3, 617--664. 


\bibitem{MM2}
Y.~Martel, F.~Merle, 
{\it 
Nonexistence of blow-up solution with minimal $L^2$-mass for the critical gKdV equation,}
Duke Math. J. 115 (2002), no. 2, 385--408. 





\bibitem{McKean}
H.P.~McKean, 
{\it Statistical mechanics of nonlinear wave equations. IV. Cubic Schr\"odinger,} 
 Comm. Math. Phys. 168 (1995), no. 3, 479--491. 
 {\it Erratum: Statistical mechanics of nonlinear wave equations. IV. Cubic Schr\"odinger}, Comm. Math. Phys. 173 (1995), no. 3, 675.
%


\bibitem{Merle93}
F.~Merle, 
{\it Determination of blow-up solutions with minimal mass for nonlinear Schr\"odinger equations with critical power,}
Duke Math. J. 69 (1993), no. 2, 427--454.

\bibitem{Merle}
F.~Merle, 
{\it Existence of blow-up solutions in the energy space for the critical generalized KdV equation},
J. Amer. Math. Soc. 14 (2001), no. 3, 555--578. 


\bibitem{MR1}
F.~Merle, P.~Rapha\"el, 
{\it Sharp upper bound on the blow-up rate for the critical nonlinear Schr\"odinger equation,}
 Geom. Funct. Anal. 13 (2003), no. 3, 591--642.

\bibitem{MR2}
F.~Merle, P.~Rapha\"el, 
{\it On universality of blow-up profile for $L^2$ critical nonlinear Schr\"odinger equation}, Invent. Math. 156 (2004), no. 3, 565--672. 


\bibitem{MR3}
F.~Merle, P.~Rapha\"el, 
{\it The blow-up dynamic and upper bound on the blow-up rate for critical nonlinear Schr\"odinger equation}, Ann. of Math.  161 (2005), no. 1, 157--222. 


\bibitem{MR4}
F.~Merle, P.~Rapha\"el, 
{\it On a sharp lower bound on the blow-up rate for the $L^2$ critical nonlinear Schr\"odinger equation}, J. Amer. Math. Soc. 19 (2006), no. 1, 37--90.

\bibitem{MRS}
F.~Merle, P.~Rapha\"el, J.~Szeftel, 
{\it The instability of Bourgain-Wang solutions for the $L^2$ critical NLS},  Amer. J. Math. 135 (2013), no. 4, 967--1017.


\bibitem{nagy}
B.V.Sz.~Nagy, 
{\it \"Uber Integralgleichungen zwischen einer Funktion und ihrer Ableitung,}
 Acta Univ. Szeged. Sect. Sci. Math
  10 (1941), 64--74. 




 
\bibitem{ns1}
K.~Nakanishi, W.~Schlag, \emph{Global dynamics above the ground state energy for the focusing nonlinear Klein-Gordon equation}, J. Differential Equations 250 (2011), no. 5, 2299--2333.

\bibitem{nsbook}
K.~Nakanishi, W.~Schlag, \emph{Invariant manifolds and dispersive Hamiltonian evolution equations}, Zurich Lectures in Advanced Mathematics. European Mathematical Society (EMS), Z\"urich, 2011. vi+253 pp.

\bibitem{ns2}
K.~Nakanishi, W.~Schlag, \emph{Global dynamics above the ground state energy for the cubic NLS equation in 3D}, Calc. Var. Partial Differential Equations 44 (2012), no. 1-2, 1--45. 

\bibitem{ns3}
K.~Nakanishi, W.~Schlag, \emph{Invariant manifolds around soliton manifolds for the nonlinear Klein-Gordon equation}, SIAM J. Math. Anal. 44 (2012), no. 2, 1175--1210.

\bibitem{ns4}
K.~Nakanishi, W.~Schlag, \emph{Global dynamics above the ground state for the nonlinear Klein-Gordon equation without a radial assumption} Arch. Ration. Mech. Anal. 203 (2012), no. 3, 809--851.


\bibitem{OgT}
T.~Ogawa, Y.~Tsutsumi, 
{\it Blow-up of solutions for the nonlinear Schr\"odinger equation with quartic potential and periodic boundary condition,} Functional-analytic methods for partial differential equations (Tokyo, 1989), 236--251, Lecture Notes in Math., 1450, Springer, Berlin, 1990.



\bibitem{OOP}
T.~Oh, M.~Okamoto, O.~Pocovnicu,
{\it  On the probabilistic well-posedness of the nonlinear Schr\"odinger equations with non-algebraic nonlinearities}, Discrete Contin. Dyn. Syst. A. 39 (2019), no. 6, 3479--3520. 

\bibitem{OOT}
T.~Oh, M.~Okamoto, L.~Tolomeo, 
{\it Focusing $\Phi^4_3$-model with a Hartree-type nonlinearity},
arXiv:2009.03251 [math.PR].



\bibitem{OOT2}
T.~Oh, M.~Okamoto, L.~Tolomeo, 
{\it Stochastic quantization of the $\Phi^3_3$-model},
arXiv:2108.06777 [math.PR].



\bibitem{OP}
T.~Oh, O.~Pocovnicu,
{\it  Probabilistic global well-posedness of the energy-critical defocusing quintic nonlinear wave equation on $\R^3$}, J. Math. Pures Appl. 105 (2016), 342--366. 



\bibitem{OQ}
T.~Oh, J.~Quastel, 
{\it On invariant Gibbs measures conditioned on mass and momentum,}
 J. Math. Soc. Japan 65 (2013), no. 1, 13--35.


\bibitem{OQV}
T.~Oh, J.~Quastel, B.~Valk\'o,
{\it  Interpolation of Gibbs measures and white noise for Hamiltonian PDE}, J. Math. Pures Appl. 97 (2012), no. 4, 391--410. 





\bibitem{ORSW}
T.~Oh, T.~Robert, P.~Sosoe, Y.~Wang,
{\it On the two-dimensional hyperbolic stochastic sine-Gordon equation}, 
Stoch. Partial Differ. Equ. Anal. Comput. 9 (2021), 1--32. 


\bibitem{OST}
T.~Oh, K.~Seong, L.~Tolomeo,
{\it  A remark on Gibbs measures with log-correlated Gaussian fields},
arXiv:2012.06729 [math.PR].

\bibitem{PW}
G.~Parisi, Y.S.~Wu, 
{\it Perturbation theory without gauge fixing,}
Sci. Sinica 24 (1981), no. 4, 483--496. 


\bibitem{Perel}
G.~Perelman, 
{\it On the formation of singularities in solutions of the critical nonlinear Schr\"odinger equation}, 
Ann. Henri Poincar\'e 2 (2001), no. 4, 605--673. 


\bibitem{PR}
F.~Planchon, P.~Rapha\"el, 
{\it Existence and stability of the log-log blow-up dynamics for the $L^2$-critical nonlinear Schr\"odinger equation in a domain},
Ann. Henri Poincar\'e 8 (2007), no. 6, 1177--1219. 

\bibitem{Poc}
O.~Pocovnicu, 
{\it Almost sure global well-posedness for the energy-critical defocusing nonlinear wave equation on 
$\R^d$, $ d= 4$ and $5$}, 
 J. Eur. Math. Soc. (JEMS) 19 (2017), no. 8, 2521--2575. 


\bibitem{Rap}
P.~Rapha\"el, 
{\it Stability of the log-log bound for blow up solutions to the critical non linear Schr\"odinger equation,} Math. Ann. 331 (2005), no. 3, 577--609. 

\bibitem{RS1}
M.~Reed, B.~Simon, 
{\it Methods of modern mathematical physics. I. Functional analysis,} Second edition. Academic Press, Inc. [Harcourt Brace Jovanovich, Publishers], New York, 1980. xv+400 pp. ISBN:


\bibitem{RS4}
M.~Reed, B.~Simon, 
{\it Methods of modern mathematical physics. IV. Analysis of operators,} Academic Press [Harcourt Brace Jovanovich, Publishers], New York-London, 1978. xv+396 pp.



\bibitem{Rider}
B.~Rider, 
{\it 
On the $\infty$-volume limit
of the focusing cubic Schr\"odinger equation},
Comm. Pure Appl. Math. 55 (2002), no. 10, 1231--1248.


\bibitem{RSS}
S.~Ryang, T.~Saito, K.~Shigemoto, 
{\it Canonical stochastic quantization}, Progr. Theoret. Phys. 73 (1985),
no. 5, 1295--1298.


\bibitem{SS}
C.~Sulem, P.-L.~Sulem, 
{\it The nonlinear Schr\"odinger equation. Self-focusing and wave collapse,} Applied Mathematical Sciences, 139. Springer-Verlag, New York, 1999. xvi+350 pp. 

\bibitem{Tao}
T.~Tao, 
{\it Nonlinear dispersive equations. Local and global analysis.}
 CBMS Regional Conference Series in Mathematics, 106. Published for the Conference Board of the Mathematical Sciences, Washington, DC; by the American Mathematical Society, Providence, RI, 2006. xvi+373 pp. 


\bibitem{TW}
L.~Tolomeo, H.~Weber, 
{\it Phase transition for invariant measures of the focusing Schr\"odinger equation}, 
 in preparation.


\bibitem{Tzv1} 
N.~Tzvetkov, 
{\it Invariant measures for the nonlinear Schr\"odinger equation on the disc}, 
 Dyn. Partial Differ. Equ. 3 (2006), no. 2, 111--160.


\bibitem{Tzv2} 
N.~Tzvetkov, 
\emph{Invariant measures for the defocusing nonlinear Schr\"odinger equation}, 
 Ann. Inst. Fourier (Grenoble) 58 (2008), no. 7, 2543--2604.



\bibitem{weinstein} 
M.~Weinstein, \emph{Nonlinear Schr\"odinger equations and sharp interpolation inequalities}, 
 Comm. Math. Phys. 87 (1983), no. 4, 567--576.



\vspace{1cm}




 


\end{thebibliography}
\end{document}